\documentclass[12pt]{amsart}

\usepackage{amssymb, amsmath, amsthm}

\numberwithin{equation}{section}

\newcommand{\op}{\mathrm{Op}}
\newcommand{\R}{\mathbb{R}}
\newcommand{\C}{\mathbb{C}}
\newcommand{\Z}{\mathbb{Z}}
\newcommand{\N}{\mathbb{N}}
\newcommand{\K}{\mathbb{K}}
\newcommand{\calA}{\mathcal{A}}
\newcommand{\calF}{\mathcal{F}}
\newcommand{\calFi}{\mathcal{F}^{-1}}

\newcommand{\calS}{\mathcal{S}}

\newcommand{\II}{I\hspace{-3pt}I}
\newcommand{\III}{I\hspace{-3pt}I\hspace{-3pt}I}
\newcommand{\supp}{\mathrm{supp}\,}

\newcommand{\andd}{\quad\mathrm{and}\quad}

\newtheorem{theorem}{Theorem}[section]
\newtheorem{proposition}[theorem]{Proposition}
\newtheorem{corollary}[theorem]{Corollary}
\newtheorem{lemma}[theorem]{Lemma}
\newtheorem{definition}[theorem]{Definition}

\theoremstyle{definition}
\newtheorem{remark}[theorem]{Remark}

\begin{document}


\title[Bilinear pseudo-differential operators]
{Bilinear pseudo-differential operators
with exotic class symbols of limited smoothness}

\author[T. Kato]{Tomoya Kato}

\address
{Division of Pure and Applied Science, 
Faculty of Science and Technology, Gunma University, 
Kiryu, Gunma 376-8515, Japan}

\email{t.katou@gunma-u.ac.jp}

\date{\today}



\thanks{This work was supported by 
the association for the
advancement of Science and Technology, 
Gunma University.}

\begin{abstract}
We consider bilinear pseudo-differential operators
with symbols in the bilinear H\"ormander class,
$BS_{\rho, \rho}^m$, $m \in \R$, $0 \leq \rho < 1$.
The aim of this paper is to discuss smoothness conditions 
for symbols to assure the boundedness
from $L^2 \times L^2$ to $h^1$
and from $L^2 \times bmo$ to $L^2$.
\end{abstract}

\maketitle



\section{Introduction}

Let $m \in \R$
and $0 \leq \delta \leq \rho \leq 1$.
The bilinear H\"ormander symbol class
$BS^m_{\rho,\delta}=
BS^m_{\rho,\delta}(\R^n)$
consists of all functions
$\sigma(x,\xi,\eta) \in C^{\infty}((\R^n)^{3})$
such that
\[
|\partial^{\alpha}_x\partial^{\beta}_{\xi}
\partial^{\gamma}_{\eta}\sigma(x,\xi,\eta)|
\le C_{\alpha,\beta,\gamma}
(1+|\xi|+|\eta|)^{m+\delta|\alpha|-\rho(|\beta|+|\gamma|)}
\]
for all multi-indices
$\alpha,\beta,\gamma \in \N_0^n
=\{0, 1, 2, \dots \}^n$. 
For a bounded measurable function 
$\sigma = \sigma (x,\xi,\eta)$ on $(\R^n)^3$,
the bilinear pseudo-differential operator
$T_{\sigma}$ is defined by
\[
T_{\sigma}(f,g)(x)
=\frac{1}{(2\pi)^{2n}}
\int_{ (\R^n)^{2} }e^{i x \cdot(\xi+\eta)}
\sigma(x,\xi,\eta)\widehat{f}(\xi)
\widehat{g}(\eta)\, d\xi d\eta
\]
for $f, g \in \calS(\R^n)$. 
To mention about the boundedness 
of the operator $T_{\sigma}$,
we will use the following terminology with a slight abuse.
Let $X,Y,Z$ be function spaces. 
If there exists a constant $C$ such that the estimate
\begin{equation*}
\| T_\sigma (f,g) \|_{Z} 
\leq C
\| f \|_{X} \| g \|_{Y}
\end{equation*}
holds
for all $f \in \calS \cap X$ and $g \in \calS \cap Y$, 
then we simply say that
the operator $T_\sigma$ 
is bounded from $X \times Y$ to $Z$.
If $\calA$ is a class of symbols,  
we denote by $\op(\calA)$
the class of all bilinear 
operators $T_{\sigma}$ 
such that $\sigma \in \calA$. 
If the operator $T_{\sigma}$ is bounded from $X \times Y$ to $Z$ 
for all $\sigma \in \calA$, 
then we write 
$\op(\calA) 
\subset B (X \times Y \to Z)$.

The boundedness of the bilinear pseudo-differential operators
with symbols in the bilinear H\"ormander class
has been studied by a lot of researches.
In the case $\rho = 1$ and $\delta < 1$,
since the bilinear operator $T_{\sigma}$ for 
$\sigma \in BS_{1,\delta}^{0}$
becomes a bilinear Carder\'on--Zygmund operator
in the sense of Grafakos--Torres \cite{grafakos torres 2002},
it is bounded from $L^p \times L^q$ to $L^r$,
$1 < p,q < \infty$,
$1/p+1/q=1/r$.
See Coifman--Meyer \cite{coifman meyer 1978}, 
B\'enyi--Torres \cite{benyi torres 2003}, and 
B\'enyi--Maldonado--Naibo--Torres 
\cite{benyi maldonado naibo torres 2010}.
Moreover, in the case $\rho=\delta=1$,
B\'enyi--Torres \cite{benyi torres 2003}
and Koezuka--Tomita \cite{koezuka tomita 2018} proved 
that the bilinear operator $T_{\sigma}$ for 
$\sigma \in BS_{1,1}^{0}$ is 
bounded from $L^{p}_{s} \times L^{q}_{s}$ to $L^{r}_{s}$,
$1 < p,q < \infty$,
$1/p+1/q=1/r$, $s > 0$,
where $L^{p}_{s}$ is the $L^{p}$-Sobolev space
equipped with the norm $\| f \|_{L^p_s} = \| (I-\Delta)^{s/2} f \|_{L^p}$.

In this paper, we are interested in the case $0 \leq \rho = \delta < 1$.
Also, we only consider the boundedness of the operator $T_{\sigma}$
from $L^2 \times L^2$ to $L^1$
and its dual, the boundedness from $L^2 \times L^\infty$ to $L^2$,
which are often understood as bases of the bilinear case.
Now, in the case $0 \leq \rho = \delta < 1$,
it is known that the class $BS_{\rho,\rho}^{0}$
does not provide the boundedness on $L^p \times L^q$.
This was pointed out by B\'enyi--Torres \cite{benyi torres 2004}.
Then, Miyachi--Tomita 
\cite{miyachi tomita 2013 IUMJ, miyachi tomita 2018 AIFG} 
proved that
\begin{equation} \label{miyachitomitaL1}
\op(BS_{\rho,\rho}^{-(1-\rho)n/2} (\R^n))
\subset B(L^2 \times L^2 \to L^1)
\end{equation}
and
\begin{equation} \label{miyachitomitaL2}
\op(BS_{\rho,\rho}^{-(1-\rho)n/2} (\R^n))
\subset B(L^2 \times L^\infty \to L^2),
\end{equation}
where, if $\rho=0$,
$L^1$ in \eqref{miyachitomitaL1} 
is replaced by the local Hardy space $h^1$
and $L^{\infty}$ in \eqref{miyachitomitaL2}
is replaced by the local $BMO$ space $bmo$.
(See also \cite{BBMNT 2013, MRS 2014}
for the preceding results in the case $m < -(1-\rho)n/2$.)
Also, in the case $\rho=0$,
these boundedness were further extended 
to sharper ones
by using $L^2$-based amalgam spaces 
in \cite{KMT-arxiv, KMT-arxiv-2}.
Now, by interpolation between 
\eqref{miyachitomitaL1} and \eqref{miyachitomitaL2}, 
it holds that
\begin{equation} \label{miyachitomitaLr}
\op(BS_{\rho,\rho}^{-(1-\rho)n/2} (\R^n))
\subset B(L^p \times L^q \to L^r)
\end{equation}
for $1 \leq r \leq 2\leq p,q \leq \infty$
and $1/p+1/q=1/r$.
In \cite{miyachi tomita 2013 IUMJ},
it was also proved that the order $m=-(1-\rho)n/2$
is critical to have \eqref{miyachitomitaL1}, 
\eqref{miyachitomitaL2}, and \eqref{miyachitomitaLr}.

We shall have a look to regularity conditions for symbols
to assure the boundedness of 
the bilinear pseudo-differential operator.
Before that, let us recall the linear case.
The linear H\"ormander symbol class,
$S_{\rho, \delta}^{m}=S_{\rho, \delta}^{m} (\R^n)$, 
is defined by the set of
all functions $\sigma \in C^\infty (\R^{n} \times \R^{n})$ such that
\begin{equation*}
| \partial_x^\alpha \partial_\xi^\beta \sigma (x,\xi) | 
\leq C_{\alpha, \beta} (1+|\xi|)^{m + \delta |\alpha| - \rho |\beta|}
\end{equation*}
for all $\alpha, \beta \in \N_0^n$.
For a bounded measurable function 
$\sigma = \sigma (x,\xi)$ on $(\R^n)^2$,
the linear pseudo-differential operator
$\sigma (X,D)$ is defined by
\begin{equation*}
\sigma (X,D) f (x) 
=
\frac{1}{( 2 \pi )^{n}} 
\int_{\mathbb{R}^n} e^{i x \cdot \xi} \sigma (x , \xi ) \widehat f (\xi ) \, d\xi
\end{equation*}
for $f \in \calS(\R^n)$.
The celebrated Calder\'on--Vaillancourt theorem in \cite{calderon vaillancourt 1972}
states that if symbols are in $S^0_{\rho,\rho}$, $0 \leq \rho <1$,
the linear pseudo-differential operator 
is bounded on $L^2$.
However, this theorem requires 
much smoothness to symbols.
Then, this was relaxed to,
roughly speaking, 
the smoothness of symbols up to 
$n/2$ for each variable $x$ and $\xi$
by Cordes \cite{cordes 1975}, 
Coifman--Meyer \cite{coifman meyer 1978}, 
Miyachi \cite{miyachi 1987}, 
Muramatu \cite{muramatu 1987}, 
Sugimoto \cite{sugimoto 1988 JMSJ}, 
and Boulkhemair \cite{boulkhemair 1995},
for $\rho = 0$, and 
Marschall \cite{marschall 1987} and
Sugimoto \cite{sugimoto 1988 TMJ},
for $0 < \rho < 1$.
Now, we shall consider the bilinear case.
For $\rho = 0$,
it was shown by \cite{KMT-arxiv, KMT-arxiv-2} that
the smoothness of symbols up to 
$n/2$ for each variable $x$, $\xi$, and $\eta$
assures the boundedness \eqref{miyachitomitaLr}.
See also
Herbert--Naibo \cite{herbert naibo 2014, herbert naibo 2016}
for the preceding results.
For $0 < \rho < 1$,
the author cannot find related results.

The purpose of this paper is to improve 
the boundedness \eqref{miyachitomitaLr}
for $0< \rho < 1$ in two ways.
Firstly, we show that
the target space $L^1$
for $r = 1$
can be replaced by $h^1$
and the domain space $L^\infty$ 
for $p = \infty$ or $q = \infty$ 
can be replaced by $bmo$
when $0 < \rho < 1$ as well as $\rho = 0$.
Secondly, we determine a smoothness condition 
of symbols for $0 < \rho < 1$
to obtain the boundedness \eqref{miyachitomitaLr}
as for $\rho=0$.

We shall state the main theorem of this paper. 
Before that, we define a Besov type symbol class.
Let $\{ \psi_{k} \}_{k \in \N_0}$ and $\{ \Psi_j \}_{j\in \N_0}$
be Littlewood--Paley partitions of unity 
on $\R^n$ and $(\R^n)^2$, respectively.
For $j \in \N_0$, 
$\boldsymbol{k} = (k_0,k_1,k_2) \in (\N_0)^3$,
$\boldsymbol{s} = (s_0,s_1,s_2) \in [0,\infty)^3$, and
$\sigma = \sigma(x,\xi,\eta) \in L^\infty((\R^n)^{3})$,
write $\boldsymbol{k} \cdot \boldsymbol{s} = k_0s_0 + k_1s_1 + k_2s_2$,
\begin{equation*}
\Delta_{\boldsymbol{k}} \sigma(x,\xi,\eta) =
\psi_{k_0}(D_x)
\psi_{k_1}(D_\xi)
\psi_{k_2}(D_\eta)
\sigma (x,\xi,\eta),
\end{equation*}
and
\begin{equation*}
\sigma_j^\rho(x,\xi,\eta)
=
\sigma (2^{-j\rho}x,2^{j\rho}\xi, 2^{j\rho}\eta)
\Psi_j (2^{j\rho}\xi, 2^{j\rho}\eta).
\end{equation*}
Then, for $m \in \R$, we denote by
$BS^{m}_{\rho,\rho}(\boldsymbol{s}; \R^{n})$
the set of all $\sigma \in L^\infty ((\R^n)^{3})$ such that
\begin{equation*}
\| \sigma \|_{BS_{\rho,\rho}^m (\boldsymbol{s};\R^{n})}
=
\sup_{ j\in\N_0} \sum_{ \boldsymbol{k} \in (\N_0)^3 }
2^{-jm + \boldsymbol{k} \cdot \boldsymbol{s}}
\| \Delta_{\boldsymbol{k}} [ \sigma_j^\rho ] \|_{ L^{2}_{ul} ((\R^{n})^3) }
< \infty,
\end{equation*}
where $L^{2}_{ul}$ is 
the uniformly local $L^2$ space
(see Section \ref{basicnotation}).
By using this symbol class, 
the main theorem of this paper is given as follows.

\begin{theorem} \label{main-thm}
Let $0 \leq \rho <1$, $m=-(1-\rho)n/2$, and
$\boldsymbol{s}=(s_0,s_1,s_2) \in [0,\infty)^3$ satisfy
$s_0 > n/2$ and $s_{1}, s_{2} \geq n/2$.
Then, if $\sigma \in BS_{\rho,\rho}^m (\boldsymbol{s};\R^{n})$,
the bilinear pseudo-differential operator $T_\sigma$ 
is bounded from $L^2 (\R^n) \times L^2 (\R^n) $ to $h^1 (\R^n) $ and 
is bounded from $L^2 (\R^n) \times bmo (\R^n) $ to $L^2 (\R^n) $.
\end{theorem}

We write
\begin{equation*}
X^p=
\left\{
\begin{array}{ll}
\vspace{3pt}
h^1, & \textrm{if } p=1, \\
\vspace{3pt}
L^p, & \textrm{if }  1<p<\infty, \\
bmo, & \textrm{if } p=\infty.
\end{array}
\right.
\end{equation*}
Then, by virtue of Theorem \ref{main-thm},
we have the following boundedness related to \eqref{miyachitomitaLr}
for the bilinear H\"ormander class
with limited smoothness.

\begin{corollary} \label{main-cor}
Let $0 \leq \rho <1$ and 
$1 \leq r \leq 2\leq p,q \leq \infty$ 
satisfy $1/p+1/q=1/r$.
Then, if
$\sigma \in C^{([n/2]+1,[n/2]+1,[n/2]+1)} 
((\R^n)^{3})$ 
satisfies that
\begin{equation*}
| \partial_{x}^{\alpha} \partial_{\xi}^{\beta} 
\partial_{\eta}^{\gamma} \sigma (x , \xi ,\eta) | 
\leq C_{\alpha, \beta, \gamma} 
\big( 1 + |\xi| + |\eta| \big)^
{-(1-\rho)n/2 + \rho ( |\alpha| - |\beta| + |\gamma| ) }
\end{equation*}
for all $\alpha,\beta,\gamma \in \N_0^n$
with
$|\alpha|,|\beta|,|\gamma| \leq [n/2]+1$,
the bilinear pseudo-differential operator $T_\sigma$ is bounded
from $X^p (\R^n) \times X^q (\R^n) $ to $X^r (\R^n) $.
\end{corollary}

We end this section by explaining the plan of this paper.
In Section \ref{secpreli}, we introduce basic notations, function spaces
and their properties which will be used throughout this paper.
In Section \ref{secmaintheorem}, 
we display two key theorems, 
Theorems \ref{l2l2h1bdd} and \ref{l2bmol2bdd}, 
and then we derive Theorem \ref{main-thm} 
and Corollary \ref{main-cor} from them.
In Section \ref{seclemma}, 
we prepare several lemmas 
which will be used to prove
Theorems \ref{l2l2h1bdd} and \ref{l2bmol2bdd}.
After we decompose symbols into 
easy forms to handle in Section \ref{secchangeform},
we actually prove Theorems \ref{l2l2h1bdd} and \ref{l2bmol2bdd} 
in Sections \ref{secl2l2h1} and \ref{secl2bmol2},
respectively.
In Section \ref{secsharp},
we consider the sharpness of indices
stated in Theorem \ref{main-thm}.
In Appendix \ref{appexistdecom},
we consider the existence of the decomposition 
used in Section \ref{secchangeform}.
In Appendix \ref{appl2l2l1bdd},
we give a small remark on
the critical case $s_0 = n/2$.


\section*{Acknowledgments}
The author sincerely express very deep thanks 
to Professor Akihiko Miyachi and Professor Naohito Tomita,
who discussed with him many times
and always gave him many remarkable comments.
Especially, they gave the author an unpublished memo
considering the boundedness from $X^p \times X^q$ to $X^r$
of bilinear Fourier multiplier operators with symbols in the exotic class.
Thanks to their fruitful helps, the author could complete this paper.


\section{Preliminaries} \label{secpreli}

\subsection{Basic notations} \label{basicnotation}

We collect notations which will be used throughout this paper.
We denote by $\R$, $\N$, $\Z$ and $\N_0$
the sets of reals, natural numbers, integers and nonnegative integers, respectively. 
We denote by $Q$ the $n$-dimensional 
unit cube $[-1/2,1/2)^n$.
A disjoint union of translations $\tau + Q$, $\tau \in \Z^n$, 
generates the Euclidean space $\R^n$.
This implies integral of a function on $\R^{n}$ can be written as 
\begin{equation}\label{cubicdiscretization}
\int_{\R^n} f(x)\, dx = 
\sum_{\tau \in \Z^{n}} \int_{Q} f(x+ \tau)\, dx. 
\end{equation}
For $x \in \R^n$ and $R>0$, 
we denote $Q(x,R)$ by 
the closed cube $x + [-R,R]^{n}$
and $B_R = B(0,R) $ by the closed ball 
$\{ x \in \R^n : | x | \leq R \}$.
We write the characteristic function 
on the set $\Omega$ as $\mathbf{1}_{\Omega}$.
For $1 \leq p \leq \infty$, $p^\prime$ is 
the conjugate number of $p$ defined by $1/p + 1/p^\prime =1$.
We write 
$[s] = \max\{ n \in \Z : n \leq s \}$ for $s \in \R$.

For two nonnegative functions $A(x)$ and $B(x)$ defined 
on a set $X$, 
we write $A(x) \lesssim B(x)$ for $x\in X$ to mean that 
there exists a positive constant $C$ such that 
$A(x) \le CB(x)$ for all $x\in X$. 
We often omit to mention the set $X$ when it is 
obviously recognized.  
Also $A(x) \approx B(x)$ means that
$A(x) \lesssim B(x)$ and $B(x) \lesssim A(x)$.

We denote the Schwartz space of rapidly 
decreasing smooth functions
on $\R^d$ 
by $\calS (\R^d)$ 
and its dual,
the space of tempered distributions, 
by $\calS^\prime(\R^d)$. 
The Fourier transform and the inverse 
Fourier transform of $f \in \calS(\R^d)$ are given by
\begin{align*}
\mathcal{F} f  (\xi) 
&= \widehat {f} (\xi) 
= \int_{\R^d}  e^{-i x \cdot \xi} f(x) \, dx, 
\\
\mathcal{F}^{-1} f (x) 
&= \check f (x)
= \frac{1}{(2\pi)^d} \int_{\R^d}  e^{i x \cdot \xi } f( \xi ) \, d\xi,
\end{align*}
respectively.
We sometimes write $\calF[f]$ and $\mathcal{F}^{-1}[f]$
when the form of $f$ is complicated. 
Furthermore, we sometimes deal with the partial Fourier transform 
of a Schwartz function $f (x,\xi,\eta)$, $x,\xi,\eta \in\R^n$.
We denote the partial Fourier transform 
with respect to the $x$, $\xi$, and $\eta$ variables by 
$\calF_0$, $\calF_1$, and $\calF_2$, respectively.
We also write $\calF_{1,2}=\calF_1 \calF_2$.
For $m \in \calS^\prime (\R^n)$, 
the Fourier multiplier operator is given by
\begin{equation*}
m(D) f 
=
\mathcal{F}^{-1} \big[ m \widehat{f}\, \big]
=
\left( \mathcal{F}^{-1} m \right) \ast f.
\end{equation*}
We also use the notation $(m(D)f)(x)=m(D_x)f(x)$ 
when we indicate which variable is considered.

For $m \in \N$, we denote by $C ( (\R^{n})^{m} )$
the set of all bounded and 
uniformly continuous functions on 
$(\R^{n})^{m}$.
For $N_i \in \N_0$, $i = 1,\dots, m$, we define
\begin{align*}
&
C^{(N_{1},\dots,N_{m})} ( (\R^{n})^{m} )
\\&=
\left\{ f : 
\partial^{ \alpha_{1} }_{x_1} \dots 
\partial^{ \alpha_{m} }_{x_{m} } f ( x_{1} ,\dots , x_{m} )
\in C ( (\R^{n})^{m} ) 
\;\;\textrm{for}\;\;
| \alpha_i | \leq N_{i}, 
i=1,\dots,m 
\right\}.
\end{align*}

For a measurable subset $E \subset \R^d$, 
the Lebesgue space $L^p (E)$, $1 \leq p\leq \infty$, 
is the set of all those 
measurable functions $f$ on $E$ such that 
$\| f \|_{L^p(E)} = 
( \int_{E} | f(x) |^p \, dx )^{1/p} 
< \infty $
if $1 \leq p < \infty$ 
or 
$\| f \|_{L^\infty (E)} 
= 
\operatorname{ess \, sup}_{x\in E} |f(x)|
< \infty$ if $p = \infty$. 
We also use the notation 
$\| f \|_{L^p(E)} = \| f(x) \|_{L^p_{x}(E)} $ 
when we indicate the variable explicitly. 
The uniformly local $L^2$ space, denoted by $L^{2}_{ul}$,
is the set of 
all measurable functions $f$ on $\R^d$ such that
\begin{equation*}
\| f \|_{L^{2}_{ul} (\R^d) } 
= \sup_{\nu \in \Z^d}
\left( \int_{[-1/2, 1/2)^d} 
\big| f(x+\nu) 
\big|^2 \, dx 
\right)^{1/2}
< \infty.
\end{equation*}

Let $\K$ be a countable set.
For $1 \leq q \leq \infty$, we denote by $\ell^q = \ell^q (\K)$ 
the set of all complex number sequences 
$\{ a_k \}_{ k \in \K }$ such that
$ \| \{ a_k \}_{ k\in\K } \|_{ \ell^q } 
= ( \sum_{ k \in \K } | a_k |^q )^{ 1/q } < \infty $
if $1 \leq q < \infty$
or 
$\| \{ a_k \}_{ k\in\K } \|_{ \ell^\infty } 
= \sup_{k \in \K} |a_k| < \infty$ 
if $q = \infty$.
For the sake of simplicity, we will write $ \| a_k \|_{ \ell^q } $
instead of the more correct notation $ \| \{ a_k \}_{ k\in\K } \|_{ \ell^q } $.
Moreover, we use the notation 
$\| a_{k} \|_{\ell^{q}} = \| a_{k} \|_{\ell^{q}_{k} (\K)} $ 
when we indicate the variable.

Let $X,Y,Z$ be function spaces.  
We denote the mixed norm by
\begin{equation*}
\| f (x,y,z) \|_{ X_x Y_y Z_z } 
= \Big\| \big\| \| f (x,y,z) 
\|_{ X_x } \big\|_{ Y_y } \Big\|_{ Z_z }.
\end{equation*} 
(Here pay special attention to the order 
of taking norms.)  
We shall use these mixed norms for  
$X, Y, Z$ being $L^p$ or $\ell^p$.

We end this subsection by stating the Schur lemma.
See, e.g., \cite[Appendix A]{grafakos 2014 modern}.

\begin{lemma} \label{schur}
Let $\{A_{j,k}\}_{j,k\in\N_{0}}$ be 
a sequence of nonnegative numbers satisfying that
$\| A_{j,k} \|_{ \ell^{1}_{k} (\N_{0}) \ell^{\infty}_{j} (\N_{0}) } \leq 1$
and
$\| A_{j,k} \|_{ \ell^{1}_{j} (\N_{0}) \ell^{\infty}_{k} (\N_{0}) } \leq 1$.
Then, we have
\begin{equation*}
\sum_{j,k\in\N_{0}} A_{j,k} b_j c_k 
\leq
\| b_j \|_{\ell^2 (\N_{0})}
\| c_k \|_{\ell^2 (\N_{0})}
\end{equation*}
for all sequences of nonnegative numbers 
$\{ b_j \}_{j\in\N_{0}}$ and $\{ c_k \}_{k\in\N_{0}}$.
\end{lemma}


\subsection{Besov spaces}

We recall the definition of the Besov space.

Let $\phi \in \calS (\R^d)$ 
satisfy that $\phi = 1$ on $\{ \xi \in \R^d : | \xi | \leq 1 \}$ 
and $\supp \phi \subset \{ \xi \in \R^d : | \xi | \leq 2 \}$. 
We set $\phi_k = \phi (\cdot / 2^k )$, $k\in\N_0$.
We write $\psi = \phi - \phi (2\cdot)$
and set $\psi_0 = \phi$ and $\psi_k = \psi (\cdot / 2^k )$, $k\in\N$.
Then, $\supp \psi \subset \{ \xi \in \R^d : 1/2 \leq | \xi | \leq 2 \}$ 
and $\sum_{k\in\N_{0}} \psi_{k} \equiv 1$.
We denote the Fourier multiplier operators 
$\psi_k (D) $ by $\Delta_k$, $k \in \N_{0}$.
We call $\{ \psi_{k} \}_{k \in \N_{0}}$ 
a Littlewood--Paley partition of unity.
For $1 \leq p,q \leq \infty$ and $s \in \R$, 
the Besov space $B_{p,q}^s (\R^d) $ 
consists of all $f \in\calS^\prime(\R^d)$
such that
\begin{equation*}
\| f \|_{B^s_{p,q}(\R^d)} 
= \left\| 2^{ks} \| \Delta_k f \|_{L^p(\R^d)} \right\|_{ \ell^q_{k} (\N_{0}) } < \infty.
\end{equation*}
We will sometimes write $\Delta_k[f]$ when the form of $f$ is complicated. 
It is known that Besov spaces are independent of the choice of
the Littlewood--Paley partition of unity. 
See \cite{triebel 1983} 
for more properties of Besov spaces.


\subsection{Local Hardy space $h^1$ and spaces $bmo$ and $BMO$}

We recall the definition of the local Hardy space $h^1(\R^n)$ 
and the spaces $bmo(\R^n)$ and $BMO(\R^n)$.

Let $\phi \in \calS(\R^n)$ be such that
$\int_{\R^n}\phi(x)\, dx \neq 0$. 
Then, the local Hardy space $h^1(\R^n)$ 
consists of
all $f \in \calS'(\R^n)$ such that
$ \|f\|_{h^1}=\|\sup_{0<t<1}|\phi_t*f|\|_{L^1}
<\infty$,
where $\phi_t(x)=t^{-n}\phi(x/t)$.
It is known that $h^1(\R^n)$
does not depend on the choice of the function $\phi$,
and that $h^1(\R^n) \hookrightarrow L^1(\R^n)$. 

The space $bmo(\R^n)$ consists of
all locally integrable functions $f$ on $\R^n$ such that
\[
\|f\|_{bmo}
=\sup_{|Q| \le 1}\frac{1}{|Q|}
\int_{Q}|f(x)-f_Q|\, dx
+\sup_{|Q|\geq1}\frac{1}{|Q|}
\int_Q |f(x)|\, dx
<\infty,
\]
where $f_Q=|Q|^{-1}\int_Q f$,
and $Q$ ranges over the cubes in $\R^n$.

The space $BMO(\R^n)$ consists of
all locally integrable functions $f$ on $\R^n$ such that
\[
\|f\|_{BMO}
=\sup_{Q}\frac{1}{|Q|}
\int_{Q}|f(x)-f_Q|\, dx
<\infty,
\]
where the supremum is taken over all cubes in $\R^n$.

It is known that
the dual space of $h^1(\R^n)$ is $bmo(\R^n)$
and that the embeddings
$L^\infty(\R^n) \hookrightarrow bmo(\R^n) \hookrightarrow BMO(\R^n)$
hold.
See Goldberg \cite{goldberg 1979} for more properties.

We end this subsection by stating the following lemma.
This was mentioned in the memo by Miyachi--Tomita.
Hence, although this lemma is not the author's contribution,
let the author give a proof below for the reader's convenience.

\begin{lemma} \label{foumulopbmo}
Let $a \geq 0$.
Suppose that 
$\varphi \in \calS(\R^n)$
and $\psi \in \calS(\R^n)$ satisfies $\psi(0)=0$.
Then, the following hold for 
any $f \in \calS (\R^n)$.
\begin{enumerate}
\setlength{\itemindent}{0pt} 
\setlength{\itemsep}{3pt} 

\item
$\| \varphi(D/2^a) f \|_{L^\infty (\R^n)} \lesssim (1+a) \| f \|_{bmo (\R^n)}$.
\vspace{3pt}

\item
$\| \psi(D/2^a) f \|_{L^\infty (\R^n)} \lesssim \| f \|_{BMO (\R^n)}$.
\end{enumerate}
Here, the implicit constants above
are independent of $a \geq 0$.
\end{lemma}

\begin{proof}
We first consider the assertion (1). We have
\begin{align*}
\varphi(D/2^a) f (x) 
&= 2^{an} \int_{\R^n} 
\check\varphi\big(2^a(x-y)\big) f(y) \,dy
\\ &\lesssim
2^{an} \int_{\R^n} (1+2^a|x-y|)^{-N} 
|f(y)-f_{Q(x,2^{-a})}| \,dy + |f_{Q(x,2^{-a})}|,
\end{align*}
where $N$ is a constant such that $N>n$ and
$f_{\Omega}=|\Omega|^{-1} \int_{\Omega} f(y) \, dy$ 
for the set $\Omega \subset \R^n$.
The first term is estimated by a constant times $\| f \|_{BMO}$
(see \cite[Proposition 3.1.5 (ii)]{grafakos 2014 modern}).
For the second term,
\begin{align*}
|f_{Q(x,2^{-a})}| 
\leq 
\sum_{0 \leq \ell \leq [a]} |f_{Q(x,2^{-a+\ell})}-f_{Q(x,2^{-a+\ell+1})}| +|f_{Q(x,2^{-a+([a]+1)})}|.
\end{align*}
Here, we have 
$|f_{Q(x,R)}-f_{Q(x,2R)}| \leq 2^n \| f \|_{BMO}$ for $R > 0$
(see \cite[Proposition 3.1.5 (i)]{grafakos 2014 modern})
and 
$|f_{Q(x,2^{-a+([a]+1)}}| \leq \| f \|_{bmo}$, 
since $-a+([a]+1) > 0$.
Hence, we obtain 
\[
|f_{Q(x,2^{-a})}| 
\lesssim (1+[a])\| f \|_{BMO} + \| f \|_{bmo}
\lesssim (1+a) \| f \|_{bmo},
\]
which completes the proof.

For (2),
since the assumption $\psi (0)=0$ means that
$\int_{\R^n} \check\psi(x)\,dx = 0$,
we have
\begin{align*}
\psi(D/2^a) f (x) 
&= 2^{an} \int_{\R^n} \check\psi\big(2^a(x-y)\big) \big( f(y) - f_{Q(x,2^{-a})} \big) \,dy.
\end{align*}
Since this is estimated by a constant times $\| f \|_{BMO}$,
we complete the proof.
\end{proof}



\section{Main theorems} \label{secmaintheorem}

\subsection{Key theorems}

In this subsection, we display two boundedness results
which immediately derive Theorem \ref{main-thm}.
These will be proved later in the next sections,
Sections \ref{secchangeform}, 
\ref{secl2l2h1}, and \ref{secl2bmol2}, by 
using lemmas which will be 
stated in Section \ref{seclemma}.

To state the results,
we give the definition of the Besov type symbol class.

\begin{definition} \label{sigmanorm12}
Let $0 \leq \rho <1$.
Let $\{ \psi_{k} \}_{k \in \N_0}$ and $\{ \Psi_{j} \}_{j\in \N_0}$
be Littlewood--Paley partition of unity 
on $\R^n$ and $(\R^n)^2$, respectively.
For $j \in \N_0$, 
$\boldsymbol{k} = (k_0,k_1,k_2) \in (\N_0)^3$,
$\boldsymbol{s} = (s_0,s_1,s_2) \in [0,\infty)^3$, and
$\sigma = \sigma(x,\xi,\eta) \in L^\infty((\R^n)^{3})$,
we write 
$\boldsymbol{k} \cdot \boldsymbol{s} = k_0s_0 + k_1s_1 + k_2s_2$,
\begin{equation*}
\Delta_{\boldsymbol{k}} \sigma(x,\xi,\eta) =
\psi_{k_0}(D_x)
\psi_{k_1}(D_\xi)
\psi_{k_2}(D_\eta)
\sigma (x,\xi,\eta),
\end{equation*}
and
\begin{equation*}
\sigma_j^\rho(x,\xi,\eta)
=
\sigma (2^{-j\rho}x,2^{j\rho}\xi, 2^{j\rho}\eta)
\Psi_j (2^{j\rho}\xi, 2^{j\rho}\eta).
\end{equation*}
Then, for $m \in \R$, we denote by
$BS^{m,\ast}_{\rho,\rho}(\boldsymbol{s}; \R^{n})$
the set of all $\sigma \in L^\infty ((\R^n)^{3})$ such that
\begin{equation*}
\| \sigma \|_{BS^{m,\ast}_{\rho,\rho}(\boldsymbol{s}; \R^{n})}
=
\sum_{k_0 \in \N_0}
\sup_{j\in\N_0}
\Bigg\{ \sum_{k_1,k_2 \in \N_0} 2^{-jm+\boldsymbol{k} \cdot \boldsymbol{s}}
\| \Delta_{\boldsymbol{k}} [ \sigma_j^\rho ] \|_{ L^{2}_{ul} ((\R^{n})^3) }
\Bigg\}
< \infty.
\end{equation*}
\end{definition}

Note that
$BS^{m}_{\rho,\rho}(\R^{n}) 
\subset BS^{m,\ast}_{\rho,\rho}(\boldsymbol{s}; \R^{n}) 
\subset BS^{m}_{\rho,\rho}(\boldsymbol{s}; \R^{n})$ 
for $\boldsymbol{s} \in [0,\infty)^3$.
See Lemma \ref{derivationlemma} below
for the first inclusion relation.
Now, we have the following theorems.

\begin{theorem}\label{l2l2h1bdd}
Let $0 \leq \rho < 1$, $m=-(1-\rho)n/2$, 
and $\boldsymbol{s}=(s_0,s_1,s_2) \in [0,\infty)^3$ satisfy
$s_0 > n/2$ and $s_1,s_2 \geq n/2$.
Then, if $\sigma \in BS_{\rho,\rho}^{m,\ast} (\boldsymbol{s};\R^{n})$,
the bilinear pseudo-differential operator $T_\sigma$ is bounded 
from $L^2 (\R^n) \times L^2 (\R^n) $ to $h^1 (\R^n) $.
\end{theorem}

\begin{theorem}\label{l2bmol2bdd}
Let $0 \leq \rho < 1$, $m=-(1-\rho)n/2$, and 
$\boldsymbol{s}=(s_0,s_1,s_2) \in [0,\infty)^3$ satisfy
$s_0,s_1,s_2 \geq n/2$.
Then, if $\sigma \in BS_{\rho,\rho}^{m,\ast} (\boldsymbol{s};\R^{n})$,
the bilinear pseudo-differential operator $T_\sigma$ is bounded 
from $L^2 (\R^n) \times bmo (\R^n) $ to $L^2 (\R^n) $.
\end{theorem}

\begin{remark} \label{rems0critical}
For the critical case $s_0 = n/2$
in Theorem \ref{l2l2h1bdd},
we can prove that 
$\op \big( BS_{\rho,\rho}^{m} (\boldsymbol{s};\R^{n}) \big)
\subset B(L^2 \times L^2 \to L^1)$
for $m=-(1-\rho)n/2$ and $s_0,s_1,s_2 \geq n/2$.
(See Appendix \ref{appl2l2l1bdd} below 
for the proof of this boundedness.)
In order to improve the target space $L^1$ to $h^1$,
we will use the small loss with respect to $s_0$.
\end{remark}



\subsection{Proofs of Theorem \ref{main-thm} and Corollary \ref{main-cor}}

Theorem \ref{main-thm} follows from Theorems \ref{l2l2h1bdd} and \ref{l2bmol2bdd} 
with the following inclusion relation among the symbol classes:
$BS^{m}_{\rho,\rho}( (s_0+\varepsilon, s_1, s_2) ; \R^{n})
\subset BS^{m,\ast}_{\rho,\rho}( (s_0, s_1, s_2) ; \R^{n})$
for any $\varepsilon > 0$.
Moreover, Corollary \ref{main-cor} is obtained
by using the following lemma,
Theorem \ref{main-thm},
and interpolation.
The idea contained in the argument
goes back to \cite[Proposition 4.7]{KMT-arxiv}.

\begin{lemma} \label{derivationlemma}
Let $0 \leq \rho < 1$, $m\in\R$, and 
$\boldsymbol{s}=(s_0,s_1,s_2) \in [0,\infty)^3$. 
Suppose that $\sigma \in C^{([s_0]+1, [s_1]+1, [s_2]+1)} ( (\R^{n})^3 )$ satisfies that
\begin{equation*}
|\partial^{\alpha}_x \partial^{\beta}_{\xi} \partial^{\gamma}_{\eta}\sigma(x,\xi,\eta)|
\leq C_{\alpha,\beta,\gamma} (1+|\xi|+|\eta|)^{m + \rho ( |\alpha|-|\beta|-|\gamma| )}
\end{equation*}
for all $\alpha, \beta, \gamma \in \N_{0}^n$ with 
$|\alpha|\leq [s_0]+1$, 
$|\beta| \leq [s_1]+1$, and 
$|\gamma| \leq [s_2]+1$.
Then, $\sigma \in BS^{m,*}_{\rho,\rho}(\boldsymbol{s}; \R^{n})$.
\end{lemma}

\begin{proof}
We let $N_i=[s_i]+1$, $i=0,1,2$.
We first consider $\Delta_{\boldsymbol{k}} [\sigma_j^\rho ]$ 
for $j,k_0,k_1,k_2 \geq 1$.
By using the Taylor expansion to the $\eta$ variable, 
the $\xi$ variable, and the $x$ variable (in this order),
together with the moment condition 
$\partial^\alpha \psi (0)=\int x^\alpha \check \psi = 0$,
we have
\begin{align} \label{taylor}
\begin{split}
&
\Delta_{\boldsymbol{k}} [ \sigma_j^\rho ](x,\xi,\eta)
=
2^{(k_0+k_1+k_2)n} 
\sum_{|\alpha| = N_0} \frac{1}{\alpha!} 
\sum_{|\beta| = N_1} \frac{1}{\beta!} 
\sum_{|\gamma| = N_2} \frac{1}{\gamma!} 
\\
&\times
\int_{(\R^{n})^3}
\int_{[0,1]^3} 
\check \psi (2^{k_0}x') (-x')^{\alpha} \,
\check \psi (2^{k_1}\xi') (-\xi')^{\beta} \,
\check \psi (2^{k_2}\eta') (-\eta')^{\gamma} 
\\
&\times \bigg( \prod_{i=0,1,2} N_i (1-t_i)^{N_i-1} \bigg)
\big( \partial_{x}^{\alpha} \partial_{\xi}^{\beta} \partial_{\eta}^{\gamma} (\sigma_j^\rho) \big)
(x-t_0x', \xi-t_1 \xi', \eta-t_2 \eta' ) 
\, dTdX',
\end{split}
\end{align}
where $dT = dt_0 dt_1 dt_2$ and $dX' = dx' d\xi' d\eta'$.
Here, we observe that
\begin{align*}
\left| \left( \partial_x^\alpha \partial_\xi^{\beta} \partial_\eta^{\gamma} \sigma \right) (2^{-j\rho}x, 2^{j\rho}\xi, 2^{j\rho}\eta) \right|
\lesssim
2^{jm + j\rho  ( |\alpha| - |\beta| - |\gamma| ) }
\end{align*}
on the support of $\Psi_j (2^{j\rho}\cdot , 2^{j\rho}\cdot)$.
Then, we have
\begin{equation}\label{partialsigmaest}
\left| \big( \partial_{x}^{\alpha} \partial_{\xi}^{\beta} \partial_{\eta}^{\gamma} (\sigma_j^\rho) \big) (x,\xi,\eta) \right|
\lesssim
\sum_{\beta'\leq\beta,\gamma'\leq\gamma}
2^{- j\rho (|\alpha| - |\beta'| - |\gamma'| ) } \, 2^{jm+ j\rho ( |\alpha| - |\beta'| - |\gamma'| ) }
\approx
2^{jm}.
\end{equation}
Collecting \eqref{taylor} and \eqref{partialsigmaest}, we have
\begin{equation*}
\left| \Delta_{\boldsymbol{k}} [ \sigma_j^\rho ](x,\xi,\eta) \right|
\lesssim
2^{jm} \, 2^{-k_0 N_0-k_1 N_1-k_2 N_2} ,
\end{equation*}
and hence
\begin{equation}
\label{sigmarhojkL2ul}
\| \Delta_{\boldsymbol{k}} [ \sigma_j^\rho ] \|_{L^{2}_{ul}}
\lesssim 2^{jm} \, 2^{-k_0 N_0-k_1 N_1-k_2 N_2}
\end{equation}
for $j,k_0,k_1,k_2 \geq 1$.
For the case $j\geq 1$ and at least one of $k_0,k_1,k_2$ is zero,
by avoiding the usage of the Taylor expansion for the corresponding variables,
we obtain the same conclusion as above.
Also, the case $j=0$ is similarly obtained.
Therefore, the estimate in 
\eqref{sigmarhojkL2ul} holds
for $j, k_0,k_1,k_2 \in \N_0$.
This means $\| \sigma \|_{BS^{m,\ast}_{\rho,\rho}(\boldsymbol{s}; \R^{n})} \lesssim 1$.
\end{proof}



\section{Lemmas}\label{seclemma}

\subsection{Elemental lemmas}
In this subsection, we denote by $S$ the operator
\begin{equation*}
S (f) (x) 
= \int_{\R^n} \frac{ f(y) }{ (1+|x-y|)^{n+1} } dy.
\end{equation*}
Note that $S$ is bounded on $L^p(\R^n)$, $1\leq p \leq \infty$.
We display basic properties of the operator $S$.
The proof can be found in \cite[Lemmas 4.1 and 4.3]{KMT-arxiv-2}.

\begin{lemma} \label{propertiesofS}
Let $1 \leq p \leq \infty$.
Then, 
the following assertions $(1)$-$(3)$ 
hold for all nonnegative functions 
$f, g$ on $\R^n$.
\begin{enumerate}
\setlength{\itemindent}{0pt} 
\setlength{\itemsep}{3pt} 
\item 
$S(f \ast g)(x) = \big( S(f) \ast g \big)(x) =  
\big( f \ast S(g) \big)(x)$.
\item 
$S(f)(x) \approx S(f)(y)$
for $x,y\in \R^n$ such that $|x-y|\lesssim1$.
\item 
$\| S(f)(x) \|_{L^p_x(\R^n)}
\approx \| S(f)(\nu) \|_{\ell^p_\nu(\Z^n)}$.
\item 
Let $\varphi$ be a function in $\calS(\R^n)$
with compact support.
Then, $| \varphi(D-\nu) f(x) |^{2}
\lesssim S( | \varphi(D-\nu) f |^{2} )(x)$
for any $f \in \calS(\R^{n})$, $\nu\in\Z^n$, and $x\in\R^n$.
\end{enumerate}
\end{lemma}

In the following lemma,
the first assertion for $R=1$ was 
proved by Miyachi--Tomita 
\cite[Lemma 3.2]{miyachi tomita 2018 AIFG}.
Moreover, the second assertion for $R=1$ 
was proved in the unpublished memo by Miyachi--Tomita.
We generalize them to the cases $R \geq 1$.

\begin{lemma} \label{estSR}
Let $2 \leq p \leq \infty$, $R \geq 1$, 
and $\varphi \in \calS (\R^n)$. 
Then, the following hold for any $f \in \calS (\R^n)$.
\begin{enumerate}
\setlength{\itemindent}{0pt} 
\setlength{\itemsep}{3pt} 
\item 
$\displaystyle
\left\| \left( \sum_{ \nu \in  \Z^n } \left|
\varphi \left( \frac{ D - \nu }{R}\right) f
\right|^2 \right)^{1/2} \right\|_{L^p(\R^n)}
\lesssim
R^{n/2} \| f \|_{L^p(\R^n)}
$. 

\item 
$\displaystyle
\left\| \left( 
\sum_{ \nu \in \Z^n : \varphi (-\nu/R)=0 } 
\left|
\varphi \left( \frac{ D - \nu }{R}\right) f
\right|^2 
\right)^{1/2} \right\|_{L^\infty(\R^n)}
\lesssim
R^{n/2} \| f \|_{BMO(\R^n)}
$.
\end{enumerate}
Here, the implicit constants above are
independent of $R \geq 1$.
\end{lemma}

\begin{proof}
We consider the assertion (1).
Since
$\R^n = \bigcup_{\nu' \in \Z^n } 2\pi \nu' + [-\pi, \pi)^n$,
we have
\begin{align*}
&
\varphi \left( \frac{D - \nu}{R} \right) f (x)
=
R^n \int_{\R^n} 
e^{ i (x-y ) \cdot \nu } 
\check \varphi \left( R ( x - y ) \right) f(y) 
\, dy
\\
&=
R^n e^{ i x \cdot \nu } 
\sum_{\nu' \in \Z^n} 
\int_{ 2\pi \nu' + [-\pi, \pi]^n}
e^{ - i y \cdot \nu } 
\check \varphi \left( R ( x - y ) \right) f(y) 
\, dy
\\&=
R^n e^{ i x \cdot \nu } 
\int_{ [-\pi, \pi]^n} e^{ - i y \cdot \nu } 
\bigg\{ 
\sum_{\nu' \in \Z^n}
\check \varphi \left( R ( x - y - 2\pi\nu' ) \right) f( y + 2\pi\nu' ) 
\bigg\} 
\, dy .
\end{align*}
Here, we realize that the function 
$\sum_{\nu' \in \Z^n} \check \varphi ( R ( x - y - 2\pi\nu' ) ) f ( y + 2\pi\nu' ) $
is $2\pi\Z^n$-periodic with respect to the $y$-variable.
Hence, we have by the Parseval identity
\begin{align*}
&
\left\|
\varphi \left( \frac{D - \nu}{R} \right) f (x)
\right\|_{\ell^2_\nu}^2
\\&= (2\pi)^{n} R^{2n} 
\int_{[-\pi, \pi]^n}
\left| 
\sum_{\nu' \in \Z^n} 
\check \varphi \left( R ( x - y - 2\pi\nu' ) \right) 
f( y + 2\pi\nu' ) 
\right|^2
\,dy.
\end{align*}
By applying the Cauchy--Schwarz inequality
to the sum over $\nu'$, we have
\begin{align*}
\left\|
\varphi \left( \frac{D - \nu}{R} \right) f (x)
\right\|_{\ell^2_\nu}^2
&\lesssim
R^{2n} \int_{ [-\pi, \pi]^n} \sum_{\nu' \in \Z^n}
\left| \check \varphi \left( R ( x - y - 2\pi\nu' ) \right) \right| | f( y+2\pi\nu' ) |^2 \, dy
\\&=
R^{2n} \left(| \check \varphi(R\cdot)| \ast |f|^2 \right) (x),
\end{align*}
where we used that
$\sum_{\nu' \in \Z^n} \left| \check \varphi \left( R ( z - 2\pi\nu' ) \right) \right|
\lesssim 1$
for $z\in\R^n$ and $R \geq 1$.
Taking the $L^{p/2}$ norm of the above,
since $2 \leq p \leq \infty$,
we have by the Young inequality
\begin{align*}
\left\| \left\|
\varphi \left( \frac{D - \nu}{R} \right) f (x)
\right\|_{\ell^2_\nu} \right\|_{L^p_x}^2
\lesssim
R^{n} \| f \|_{L^p}^2,
\end{align*}
which completes the proof of the assertion (1).

For the assertion (2), 
we have
for $\nu \in \Z^n$ such that $\varphi (-\nu/R)=0$
\begin{equation*}
\varphi \left( \frac{ D - \nu }{R}\right) f(x) 
=
R^n e^{ i x \cdot \nu }
\int_{\R^n} e^{ - i y \cdot \nu } 
\check \varphi \left( R ( x - y ) \right) \big( f(y) -c \big) dy 
\end{equation*}
for any constant $c \in \C$.
Repeating the same lines 
as above with $f(y) -c$, 
we obtain
\begin{align*}
\sum_{ \nu \in \Z^n : \varphi (-\nu/R)=0 }
\left| \varphi \left( \frac{ D - \nu }{R}\right) f(x) \right|^2 
\lesssim
R^{2n} \left(| \check \varphi(R\cdot)| \ast |f-c|^2 \right) (x).
\end{align*}
Choose 
$c = f_{Q(x,R^{-1})}=|Q(x,R^{-1})|^{-1} \int_{Q(x,R^{-1})} f(y) \, dy$
and observe that
\begin{align*}
R^{n} 
\int_{ \R^n} 
\left| \check \varphi \left(R ( x - y ) \right) \right|
\left| f( y ) -f_{Q(x,R^{-1})} \right|^2 
\, dy
\lesssim \| f \|_{BMO}^2
\end{align*}
(with the aid of 
\cite[Proposition 3.1.5 (i) and Corollary 3.1.8]{grafakos 2014 modern}).
Then, we have
\begin{equation*}
\sum_{ \nu \in \Z^n : \varphi (-\nu/R)=0 }
\left| \varphi \left( \frac{ D - \nu }{R}\right) f(x) \right|^2 
\lesssim R^{n} \| f \|_{BMO}^2,
\end{equation*}
which completes the proof of the assertion (2).
\end{proof}

\begin{corollary}\label{estSRcor}
Let $r > 0$, $R \geq 1$, and $\varphi , \phi \in \calS (\R^n)$. 
Then,
\begin{equation*}
\left\| \left( 
\sum_{ \nu \in \Z^n : \varphi (-\nu/R)=0 } 
\left|
\varphi \left( \frac{ D - \nu }{R}\right) 
\phi \left( \frac{D}{r} \right) f
\right|^2 
\right)^{1/2} \right\|_{L^\infty(\R^n)}
\lesssim
R^{n/2} \| f \|_{BMO(\R^n)},
\end{equation*}
where the implicit constant above is
independent of $r>0$ and $R \geq 1$.
\end{corollary}

\begin{proof}
The summand can be written by
\begin{equation*}
\varphi \left( \frac{ D - \nu }{R}\right) 
\phi \left( \frac{D}{r} \right) f (x)
=
r^n \int_{\R^n}
\varphi \left( \frac{ D - \nu }{R}\right) f (x-y) 
\, \check \phi (ry)
\, dy.
\end{equation*}
We take the $\ell^2_{\nu}$ norm restricted to
the set $\{ \nu \in \Z^n : \varphi (-\nu/R)=0 \}$ of the above 
and use the Minkowski inequality for the integral.
Then, by Lemma \ref{estSR} (2),
we have the desired results, since
$r^{n} \| \check \phi (r\cdot) \|_{L^1} \approx 1$
for $r>0$.
\end{proof}


\subsection{Lemmas for Theorems \ref{l2l2h1bdd} and \ref{l2bmol2bdd}}

In this subsection,
we show some lemmas for the dual form of
bilinear pseudo-differential operators.
The basic idea for the argument used here 
goes back to Boulkhemair \cite{boulkhemair 1995}.
See also \cite[Proposition 5.1]{KMT-arxiv-2}.

Throughout this subsection, 
we will denote the Fourier multiplier operator 
$\kappa (D-\mu)$ by $\square_{\mu}$
for $\mu \in \Z^n$, where $\kappa \in \calS(\R^n)$.
Moreover,
the norm 
$\| \sigma_{\boldsymbol{\nu}} (x,\xi,\eta) 
\|_{L^2_{\xi,\eta}L^2_{ul,x} \ell^\infty_{\boldsymbol{\nu}}}$
will be abbreviated to $\| \sigma_{\boldsymbol{\nu}} \|$
for the sake of simplicity.


\begin{lemma} \label{nu12lemma}
Let $R_0,R_1,R_2 \geq 1$ and 
$2\leq p,q,r \leq \infty$ satisfy $1/p+1/q+1/r=1$.
Let $\Lambda, \Lambda_{1}, \Lambda_{2}$
be subsets of $\Z^n$.
Suppose that 
$\kappa \in \calS (\R^n)$ satisfies that
$\supp \kappa \subset [-1,1]^n$ and 
that $\{ \sigma_{\boldsymbol{\nu}} \} $, 
${\boldsymbol{\nu}} = (\nu_1, \nu_2 ) \in (\Z^{n})^2$, 
is a sequence of bounded functions on $(\R^n)^3$
satisfying that
$\supp \calF [\sigma_{\boldsymbol{\nu}}]
\subset B_{R_0} \times B_{R_1} \times B_{R_2}$ 
for ${\boldsymbol{\nu}} \in (\Z^{n})^2$.
Then, the following $(1)$ and $(2)$ hold
for any $f,g,h \in \calS(\R^n)$,
and $(3)$ holds for any
$f,g \in \calS(\R^n)$ and 
$\{ h_{\tau} \}_{\tau \in \Z^n} \subset \calS(\R^n)$.
\begin{align}
\tag{1}
&
\sum_{ {\boldsymbol{\nu}} \in\Lambda_1 \times \Lambda_2 }
\left|
\int_{\R^n} T_{ \sigma_{\boldsymbol{\nu}} }( \square_{\nu_1} f, \square_{\nu_2} g )(x)
\, h(x) \, dx
\right|
\\\nonumber
&\quad\lesssim
\min_{i=1,2} \left( |\Lambda_i|^{1/2} \right)
( R_0 R_1 R_2 )^{n/2} 
\| \sigma_{\boldsymbol{\nu}} \| \,
\| f \|_{L^{p}(\R^n)} 
\| \square_{\nu_2} g \|_{\ell^2_{\nu_2}(\Lambda_2) L^{q}(\R^n)} 
\| h \|_{L^{r}(\R^n)}.
\\
\tag{2}&
\sum_{ {\boldsymbol{\nu}} \in \Lambda_1 \times \Lambda_2 }
\left| 
\int_{\R^n}
T_{ \sigma_{\boldsymbol{\nu}} } (\square_{\nu_1} f, \square_{\nu_2} g)(x)
\, h(x) \, dx
\right|
\\
\nonumber
&\quad\lesssim
\left( |\Lambda_1| \, |\Lambda_2| \right)^{1/2} 
( R_1 R_2 )^{n/2}
\| \sigma_{\boldsymbol{\nu}} \| \,
\| f \|_{L^{p}(\R^n)} 
\| \square_{\nu_2} g \|_{\ell^2_{\nu_2}(\Lambda_2) L^{q}(\R^n)} 
\| h \|_{L^{r}(\R^n)}.
\\
\tag{3}&
\sum_{ \tau \in \Lambda } \sum_{ \nu_1 \in \Z^n : \, \nu_1 + \nu_2 = \tau }
\left|
\int_{\R^n}
T_{ \sigma_{\boldsymbol{\nu}} } (\square_{\nu_1} f, \square_{\nu_2} g)(x)
\, h_{\tau} (x) \, dx
\right|
\\
\nonumber
&\quad\lesssim
|\Lambda| ^{1/2} 
( R_1 R_2 )^{n/2} 
\| \sigma_{\boldsymbol{\nu}} \| \,
\| f \|_{L^{p}(\R^n)} 
\| g \|_{L^{q}(\R^n)} 
\| h_{\tau} \|_{ \ell^2_{\tau}(\Lambda) L^{r}(\R^n) }.
\end{align}
Here, the absolute values 
of $\Lambda$, $\Lambda_{1}$, and $\Lambda_{2}$
are the cardinality of these sets,
and the implicit constants are independent of $R_0, R_1, R_2$.
In particular,
$\| \square_{\nu_2} g \|_{\ell^2_{\nu_2}(\Lambda_2) L^{q}(\R^n)}$
in $(1)$ and $(2)$
can be replaced 
by $\| g \|_{L^q(\R^n)}$.
\end{lemma}

To prove this, we use the following lemma.

\begin{lemma} \label{kyoutuulem}
Let $R_1,R_2 \geq 1$.
Suppose that 
$\kappa \in \calS (\R^n)$ satisfies that $\supp \kappa \subset [-1,1]^n$ and 
that $\{ \sigma_{\boldsymbol{\nu}} \} $, 
${\boldsymbol{\nu}} = (\nu_1, \nu_2 ) \in (\Z^{n})^2$, 
is a sequence of bounded functions on $(\R^n)^3$
such that
$\supp \calF_{1,2} [\sigma_{\boldsymbol{\nu}}] (x, \cdot, \cdot)
\subset B_{R_1} \times B_{R_2}$
for ${\boldsymbol{\nu}} \in (\Z^{n})^2$ and 
$x \in \R^n$.
Then, 
\begin{align} \label{kyoutuulemest}
\begin{split}
&
\left|
\int_{\R^n} 
T_{ \sigma_{\boldsymbol{\nu}} }( \square_{\nu_1} f, \square_{\nu_2} g )(x) 
\, h_{\nu_1 + \nu_2}(x) 
\, dx
\right|
\\&\lesssim
\| \sigma_{\boldsymbol{\nu}} \|
\sum_{ \nu_0\in\Z^n }
\left\{ S \left( \mathbf{1}_{B_{R_1}} \ast \big| \square_{\nu_1} f \big|^2 \right)(\nu_0) \right\}^{1/2}
\left\{ S \left( \mathbf{1}_{B_{R_2}} \ast \big| \square_{\nu_2} g \big|^2 \right)(\nu_0) \right\}^{1/2} 
\\
&\quad\times \left\| h_{\nu_1 + \nu_2} (x+\nu_0) \right\|_{L^2_x(Q)}
\end{split}
\end{align}
for any ${\boldsymbol{\nu}} \in (\Z^{n})^2$,
$x \in \R^n$,
$f,g\in \calS(\R^n)$, 
and $\{ h_{\tau} \}_{ \tau \in \Z^n } \subset \calS(\R^n)$. 
\end{lemma}

\begin{proof}
We simply denote by $I$ 
the left hand side of \eqref{kyoutuulemest}.
Observe
that
\begin{align*}
&
(2\pi)^{2n} \, T_{ \sigma_{\boldsymbol{\nu}} } ( \square_{\nu_1} f, \square_{\nu_2} g ) (x)
\\
&=
\int_{(\R^n)^2}
	\calF_{1,2} [\sigma_{\boldsymbol{\nu}}] (x, y-x, z-x)
	\, \mathbf{1}_{B_{R_1}}(x-y) \square_{\nu_1} f (y) 
	\, \mathbf{1}_{B_{R_2}}(x-z) \square_{\nu_2} g (z) \,dydz.
\end{align*}
Then, the Cauchy--Schwarz inequalities 
and the Plancherel theorem give
\begin{align*}
\left| 
T_{ \sigma_{\boldsymbol{\nu}} } \left( \square_{\nu_1} f, \square_{\nu_2} g \right) (x)
\right|^2
\lesssim
\big\| \sigma_{\boldsymbol{\nu}} (x,\xi,\eta) \big\|_{L^2_{\xi,\eta}}^2
\Big( \mathbf{1}_{B_{R_1}} \ast \big|\square_{\nu_1} f\big|^2 \Big) (x)
\Big( \mathbf{1}_{B_{R_2}} \ast \big|\square_{\nu_2} g\big|^2 \Big) (x)
\end{align*}
for any ${\boldsymbol{\nu}} \in (\Z^{n})^2$ and $x \in \R^n$.
From this, it holds that
\begin{align}\label{beforedecomposed}
\begin{split}
I &\lesssim
\int_{\R^n} 
\big\| \sigma_{\boldsymbol{\nu}} (x,\xi,\eta) \big\|_{L^2_{\xi,\eta}}
\left| h_{\nu_1 + \nu_2}(x) \right|
\\&\quad\times
\Big\{ \Big( \mathbf{1}_{B_{R_1}} \ast \big| \square_{\nu_1} f \big|^2 \Big) (x) \Big\}^{1/2}
\Big\{ \Big( \mathbf{1}_{B_{R_2}} \ast \big| \square_{\nu_2} g \big|^2 \Big) (x) \Big\}^{1/2} \, dx.
\end{split}
\end{align}
We separate the integral 
by using \eqref{cubicdiscretization}.
Then, the inequality \eqref{beforedecomposed}
coincides with
\begin{align*}
I &\lesssim
\sum_{ \nu_0\in\Z^n }
\int_{Q} 
\big\| \sigma_{\boldsymbol{\nu}} (x+\nu_0,\xi,\eta) \big\|_{L^2_{\xi,\eta}}
\left| h_{\nu_1 + \nu_2}(x+\nu_0) \right|
\\&\quad\times
\Big\{ \Big( \mathbf{1}_{B_{R_1}} \ast \big| \square_{\nu_1} f \big|^2 \Big) (x+\nu_0) \Big\}^{1/2}
\Big\{ \Big( \mathbf{1}_{B_{R_2}} \ast \big| \square_{\nu_2} g \big|^2 \Big) (x+\nu_0) \Big\}^{1/2} \, dx.
\end{align*}
Since we have
by Lemma \ref{propertiesofS} (4), (1), and (2)
\begin{align*}
\Big( \mathbf{1}_{B_{R_1}} \ast \big|\square_{\nu_1} f\big|^2 \Big) (x+\nu_0)
\lesssim
S \Big( \mathbf{1}_{B_{R_1}} \ast \big|\square_{\nu_1} f \big|^2 \Big) (\nu_0)
\end{align*}
for $x \in Q$, we obtain
\begin{align*}
I &\lesssim
\sum_{ \nu_0\in\Z^n }
\Big\{ S \Big( \mathbf{1}_{B_{R_1}} \ast \big| \square_{\nu_1} f \big|^2 \Big) (\nu_0) \Big\}^{1/2}
\Big\{ S \Big( \mathbf{1}_{B_{R_2}} \ast \big| \square_{\nu_2} g \big|^2 \Big) (\nu_0) \Big\}^{1/2}
\\&\quad\times 
\int_{Q} 
\big\| \sigma_{\boldsymbol{\nu}} (x+\nu_0,\xi,\eta) \big\|_{L^2_{\xi,\eta}}
\left| h_{\nu_1 + \nu_2}(x+\nu_0) \right| \, dx
\\&\leq
\sum_{ \nu_0\in\Z^n }
\Big\{ S \Big( \mathbf{1}_{B_{R_1}} \ast \big| \square_{\nu_1} f \big|^2 \Big) (\nu_0) \Big\}^{1/2}
\Big\{ S \Big( \mathbf{1}_{B_{R_2}} \ast \big| \square_{\nu_2} g \big|^2 \Big) (\nu_0) \Big\}^{1/2}
\\
&\quad\times
\big\| \sigma_{\boldsymbol{\nu}} (x+\nu_0,\xi,\eta) \big\|_{L^2_{\xi,\eta}L^2_{x}(Q)}
\left\| h_{\nu_1 + \nu_2}(x+\nu_0) \right\|_{L^2_x(Q)}.
\end{align*}
Taking the supremum over $\nu_0$ and $\boldsymbol{\nu}$ 
of the factor of $\sigma_{\boldsymbol{\nu}}$, 
we complete the proof.
\end{proof}

Now, we shall prove Lemma \ref{nu12lemma}.

\begin{proof}[Proof of Lemma \ref{nu12lemma} (1)]
We may assume that 
$|\Lambda_{1}| < \infty$ or $|\Lambda_{2}| < \infty$.
We simply denote by $I$ 
the left hand side of the assertion (1).
We first observe that
\begin{align*}
&
\calF \left[ T_{ \sigma_{\boldsymbol{\nu}} } ( \square_{\nu_1} f, \square_{\nu_2} g ) \right](\zeta)
\\&=
\frac{1}{(2\pi)^{2n}}
\int_{ (\R^{n})^2 }
\calF_0 [\sigma_{\boldsymbol{\nu}}] \big( \zeta - (\xi + \eta ) , \xi , \eta \big) 
\kappa (\xi-{\nu_1}) \widehat{f}(\xi)
\kappa (\eta-{\nu_2}) \widehat{g}(\eta)
\, d\xi d \eta.
\end{align*}
Here, since
$\supp \calF_{0} [\sigma_{\boldsymbol{\nu}}] (\cdot,\xi,\eta) 
\subset B_{R_0}$ 
and $\supp \kappa(\cdot-\nu_i) \subset \nu_i+[-1,1]^n$,
we see that
\begin{equation*}
\supp \calF
\left[ T_{ \sigma_{\boldsymbol{\nu}} } ( \square_{\nu_1} f, \square_{\nu_2} g ) \right]
\subset 
\big\{ \zeta \in \R^n : | \zeta-(\nu_1 + \nu_2) | \lesssim R_0 \big\}.
\end{equation*}
We take a function $\varphi \in \calS(\R^n)$ satisfying that 
$\varphi=1$ on $\{ \zeta\in\R^n: |\zeta|\lesssim 1\}$.
Then,
\begin{equation*}
I=
\sum_{ {\boldsymbol{\nu}} \in\Lambda_1 \times \Lambda_2 }
\left|
\int_{\R^n} T_{ \sigma_{\boldsymbol{\nu}} }( \square_{\nu_1} f, \square_{\nu_2} g )(x)
\, \varphi \left( \frac{D+\nu_1 + \nu_2}{R_0} \right) h (x) \, dx
\right|.
\end{equation*}
By the use of Lemma \ref{kyoutuulem},
\begin{align}\label{h1lemma1afterlem}
\begin{split}
I&\lesssim
\| \sigma_{\boldsymbol{\nu}} \|
\sum_{ \nu_0\in\Z^n }
\sum_{ {\boldsymbol{\nu}} \in\Lambda_1 \times \Lambda_2 }
\left\|\varphi \left( \frac{D+\nu_1 + \nu_2}{R_0} \right) h(x+\nu_0) \right\|_{L^2_x(Q)}
\\
&\quad\times 
\Big\{ S \Big( \mathbf{1}_{B_{R_1}} \ast \big| \square_{\nu_1} f \big|^2 \Big)(\nu_0) \Big\}^{1/2}
\Big\{ S \Big( \mathbf{1}_{B_{R_2}} \ast \big| \square_{\nu_2} g \big|^2 \Big)(\nu_0) \Big\}^{1/2}.
\end{split}
\end{align}
In what follows, we simply write each summand by
\begin{align*}
F(\nu_1,\nu_0)
&=\left\{ S \left( \mathbf{1}_{B_{R_1}} \ast \big| \square_{\nu_1} f \big|^2 \right)(\nu_0) \right\}^{1/2},
\\
G(\nu_2,\nu_0)
&=\left\{ S \left( \mathbf{1}_{B_{R_2}} \ast \big| \square_{\nu_2} g \big|^2 \right)(\nu_0) \right\}^{1/2},
\\
H(\nu_3,\nu_0)
&= \left\|\varphi \left( \frac{D+\nu_3}{R_0} \right) h(x+\nu_0) \right\|_{L^2_x(Q)}.
\end{align*}
Then the inequality \eqref{h1lemma1afterlem} is rewritten as
\begin{equation} \label{h1lemma1I}
I\lesssim
\| \sigma_{\boldsymbol{\nu}} \| \,
\II
\end{equation}
with
\begin{equation*}
\II=
\sum_{ \nu_0\in\Z^n } \sum_{\nu_1 \in \Lambda_1 } \sum_{ \nu_2 \in \Lambda_2 }
F(\nu_1,\nu_0)G(\nu_2,\nu_0)H(\nu_1+\nu_2,\nu_0).
\end{equation*}

Let us estimate $\II$.
We first consider the case $|\Lambda_1| \leq |\Lambda_2|$.
We apply the Cauchy--Schwarz inequality 
firstly to the sum over $\nu_2$
and secondly to the sum over $\nu_1$, and thirdly
apply the H\"older inequality with $1/p+1/q+1/r=1$
to the sum over $\nu_0$.
Then, 
\begin{align*}
\II 
&\leq
\sum_{ \nu_0\in\Z^n }
\sum_{ \nu_1 \in\Lambda_1 }
F(\nu_1,\nu_0) 
\| G(\nu_2,\nu_0) \|_{\ell^2_{\nu_2}(\Lambda_2)} 
\| H(\nu_3,\nu_0) \|_{\ell^2_{\nu_3}(\Z^n)}
\\
&\leq
| \Lambda_1|^{1/2}
\sum_{ \nu_0\in\Z^n }
\|F(\nu_1,\nu_0)\|_{\ell^2_{\nu_1}(\Lambda_1)} 
\|G(\nu_2,\nu_0)\|_{\ell^2_{\nu_2}(\Lambda_2)} 
\|H(\nu_3,\nu_0)\|_{\ell^2_{\nu_3}}
\\
&\leq
| \Lambda_1|^{1/2}
\|F(\nu_1,\nu_0)\|_{\ell^2_{\nu_1}\ell^p_{\nu_0}} 
\|G(\nu_2,\nu_0)\|_{\ell^2_{\nu_2}(\Lambda_2)\ell^q_{\nu_0}}
\|H(\nu_3,\nu_0)\|_{\ell^2_{\nu_3} \ell^r_{\nu_0}}.
\end{align*}
For the opposite case $| \Lambda_2 | \leq | \Lambda_1 |$,
we switch the order to use the Cauchy--Schwarz inequalities 
with respect to $\nu_1$ and $\nu_2$. 
Thus, we obtain
\begin{equation} \label{h1lemma1II}
\II \leq
\min_{i=1,2} \left( |\Lambda_i|^{1/2} \right)
\|F(\nu_1,\nu_0)\|_{\ell^2_{\nu_1}\ell^p_{\nu_0}} 
\|G(\nu_2,\nu_0)\|_{\ell^2_{\nu_2}(\Lambda_2)\ell^q_{\nu_0}} 
\|H(\nu_3,\nu_0)\|_{\ell^2_{\nu_3} \ell^r_{\nu_0}}.
\end{equation}

We shall estimate the norms of $F$, $G$, and $H$.
Firstly, for the norm of $F$,
we have by Lemma \ref{propertiesofS} (3)
\begin{equation}\label{normofA}
\big\|F(\nu_1,\nu_0)\big\|_{\ell^2_{\nu_1}\ell^p_{\nu_0}}
\approx
\left\|
S \Big( \mathbf{1}_{B_{R_1}} \ast \big\| \square_{\nu_1} f \big\|_{\ell^2_{\nu_1}}^2 \Big)(x) 
\right\|_{L^{p/2}_{x}}^{1/2}.
\end{equation}
By the boundedness of the operator $S$
on $L^{p/2}(\R^n)$, $p\geq 2$,
and the Young inequality,
the right hand side of \eqref{normofA}
is estimated by a constant times
\begin{equation}\label{normofAtochuu}
\left\| 
\mathbf{1}_{B_{R_1}} \ast \big\| \square_{\nu_1} f \big\|_{\ell^2_{\nu_1}}^2 
\right\|_{L^{p/2}}^{1/2}
\lesssim
R_1^{n/2} 
\left\| 
\big\| \square_{\nu_1} f \big\|_{\ell^2_{\nu_1}}
\right\|_{L^{p}}.
\end{equation}
Since
$\square_{\nu_1} = \kappa (D-\nu_1)$ 
with $\kappa\in \calS(\R^n)$,
we have
$\| \square_{\nu_1} f \|_{\ell^2_{\nu_1} L^{p}}
\lesssim \| f \|_{L^p}$
by Lemma \ref{estSR} (1) with $R=1$.
Therefore, we obtain
\begin{align} \label{h1lemma1A}
\big\|F(\nu_1,\nu_0) \big\|_{\ell^2_{\nu_1}\ell^p_{\nu_0}}
\lesssim
R_1^{n/2} \| f \|_{L^p}.
\end{align}
For the norm of $G$,
repeating the same line as for $F$,
we have by \eqref{normofA} and \eqref{normofAtochuu}
\begin{equation} \label{h1lemma1B}
\| G(\nu_2,\nu_0) \|_{\ell^2_{\nu_2} (\Lambda_2) \ell^q_{\nu_0}}
\lesssim R_2^{n/2} 
\left\| 
\big\| \square_{\nu_2} g \big\|_{\ell^2_{\nu_2}(\Lambda_2)}
\right\|_{L^{q}}.
\end{equation}
Lastly, we consider the norm of $H$.
Since
$L^{r}(Q) \hookrightarrow L^{2}(Q)$ 
for $2 \leq r \leq \infty$,
we have by Lemma \ref{estSR} (1)
\begin{align} \label{h1lemma1C}
\begin{split}
\|H(\nu_3,\nu_0)\|_{\ell^2_{\nu_3} \ell^r_{\nu_0}}
&=
\left\|\varphi \left( \frac{D+\nu_3}{R_0} \right) h(x+\nu_0)
 \right\|_{ \ell^2_{\nu_3} L^{2}_x(Q) \ell^r_{\nu_0}}
 \\&\leq
\left\|\varphi \left( \frac{D+\nu_3}{R_0} \right) h(x)
 \right\|_{ \ell^2_{\nu_3} L^{r}_x(\R^n)}
\lesssim
R_0^{n/2} \| h \|_{L^{r}}.
\end{split}
\end{align}
Thus, collecting 
\eqref{h1lemma1I}, \eqref{h1lemma1II}, 
\eqref{h1lemma1A}, \eqref{h1lemma1B},
and \eqref{h1lemma1C}, 
we obtain the desired estimate.
Moreover, by virtue of Lemma \ref{estSR} (1),
$\left\| \square_{\nu_2} g 
\right\|_{\ell^2_{\nu_2}(\Lambda_2) L^{q}}$
can be replaced by $\| g \|_{L^q}$.
\end{proof}

\begin{proof}[Proof of Lemma \ref{nu12lemma} (2)]
We may assume $|\Lambda_{1}|, | \Lambda_{2} | < \infty$.
We simply denote by $I$ 
the left hand side of the inequality of the assertion (2).
Since $h$ is independent of $\nu_1+\nu_2$, 
we have by Lemma \ref{kyoutuulem}
\begin{align*}
I&\lesssim
\| \sigma_{\boldsymbol{\nu}} \| \,
\sum_{ \nu_0\in\Z^n }
\| h (x+\nu_0) \|_{L^2_x(Q)}
\\
&\quad\times 
\sum_{ {\boldsymbol{\nu}} \in\Lambda_1 \times \Lambda_2 }
\left\{ S \left( \mathbf{1}_{B_{R_1}} \ast \big| \square_{\nu_1} f \big|^2 \right)(\nu_0) \right\}^{1/2}
\left\{ S \left( \mathbf{1}_{B_{R_2}} \ast \big| \square_{\nu_2} g \big|^2 \right)(\nu_0) \right\}^{1/2} .
\end{align*}
We use the Cauchy--Schwarz inequalities 
to the sums over $\nu_1$ and $\nu_2$
and then use the H\"older inequality 
to the sum over $\nu_0$
with $1/p+1/q+1/r=1$.
Then,
\begin{align*}
&
I\lesssim
\left( |\Lambda_1| \, |\Lambda_2| \right)^{1/2} \| \sigma_{\boldsymbol{\nu}} \|\,
\| h (x+\nu_0) \|_{L^2_x(Q)\ell^r_{\nu_0}}
\\
&\times 
\left\| 
	\left\{ S \left( \mathbf{1}_{B_{R_1}} \ast \big| \square_{\nu_1} f \big|^2 \right)(\nu_0) \right\}^{1/2}
\right\|_{\ell^2_{\nu_1}(\Lambda_1)\ell^p_{\nu_0}}
\left\| 
	\left\{ S \left( \mathbf{1}_{B_{R_2}} \ast \big| \square_{\nu_2} g \big|^2 \right)(\nu_0) \right\}^{1/2}
\right\|_{\ell^2_{\nu_2}(\Lambda_2)\ell^q_{\nu_0}}.
\end{align*}
Here, from \eqref{h1lemma1A} and \eqref{h1lemma1B},
the factors of $f$ and $g$ are bounded by 
a constant times $R_1^{n/2} \| f \|_{L^{p}}$ and 
$R_2^{n/2} \| \square_{\nu_2} g \|_{\ell^2_{\nu_2}(\Lambda_2) L^{q}}$,
respectively.
From the embedding
$L^{r}(Q) \hookrightarrow L^{2}(Q)$ for $2 \leq r \leq \infty$,
$\| h (x+\nu_0) \|_{L^2_x(Q)\ell^r_{\nu_0}} \leq \|h\|_{L^{r}}$.
Hence, we obtain the desired estimate.
Also, by Lemma \ref{estSR} (1),
$\left\| \square_{\nu_2} g 
\right\|_{\ell^2_{\nu_2}(\Lambda_2) L^{q}}$
can be replaced by $\| g \|_{L^q}$.
\end{proof}


\begin{proof}[Proof of Lemma \ref{nu12lemma} (3)]
We may assume $|\Lambda| < \infty$.
We denote by $I$ 
the left hand side of the assertion (3).
By Lemma \ref{kyoutuulem},
\begin{align*}
I&\lesssim
\| \sigma_{\boldsymbol{\nu}} \|
\sum_{\nu_0\in\Z^n}
\sum_{ \tau \in \Lambda } 
\sum_{ \nu_1 :\, \nu_1 + \nu_2 = \tau }
\| h_{\tau} (x+\nu_0) \|_{L^2_x(Q)}
\\
&\quad\times 
\Big\{ S \Big( \mathbf{1}_{B_{R_1}} \ast \big| \square_{\nu_1} f \big|^2 \Big)(\nu_0) \Big\}^{1/2}
\Big\{ S \Big( \mathbf{1}_{B_{R_2}} \ast \big| \square_{\nu_2} g \big|^2 \Big)(\nu_0) \Big\}^{1/2}.
\end{align*}
Firstly, we apply the Cauchy--Schwarz inequality to the sum over $\nu_1$,
since $\nu_2=\tau-\nu_1$.
Secondly, we use the Cauchy--Schwarz inequality to the sum over $\tau$
and thirdly use the H\"older inequality to the sum over $\nu_0$ 
with $1/p+1/q+1/r = 1$.
Then,
\begin{align*}
&
I\lesssim
|\Lambda| ^{1/2}
\| \sigma_{\boldsymbol{\nu}} \|
\| h_{\tau} (x+\nu_0) \|_{L^2_x(Q) \ell_{\tau}^{2} (\Lambda) \ell^{r}_{\nu_0} }
\\
&\times 
\left\| 
	\Big\{ S \Big( \mathbf{1}_{B_{R_1}} \ast \big| \square_{\nu_1} f \big|^2 \Big)(\nu_0) \Big\}^{1/2} 
\right\|_{ \ell^2_{\nu_1} \ell^{p}_{\nu_0} }
\left\| 
	\Big\{ S \Big( \mathbf{1}_{B_{R_2}} \ast \big| \square_{\nu_2} g \big|^2 \Big)(\nu_0) \Big\}^{1/2} 
\right\|_{ \ell^2_{\nu_2} \ell^{q}_{\nu_0} }.
\end{align*}
By \eqref{h1lemma1A},
the factors of $f$ and $g$ are bounded by 
a constant times $R_1^{n/2} \| f \|_{L^{p}}$
and $R_2^{n/2} \| g \|_{L^{q}}$, respectively.
For the factor of $h_{\tau}$, since
$L^{r}(Q) \hookrightarrow L^{2}(Q)$ for $2 \leq r \leq \infty$,
$\| h_{\tau} (x+\nu_0) \|_{L^2_x(Q) \ell_{\tau}^{2} (\Lambda) \ell^{r}_{\nu_0} }
\leq \| h_{\tau} \|_{ \ell^2_{\tau}(\Lambda) L^{r} (\R^n)}$.
Hence, we obtain the desired estimate.
\end{proof}



\section{Decomposition of symbols}
\label{secchangeform}

In this section, we decompose symbols
of the bilinear operator $T_{\sigma}$
by Littlewood--Paley partitions and the following lemma
given by Sugimoto \cite[Lemma 2.2.1]{sugimoto 1988 JMSJ}.
An explicit proof can be found in \cite[Lemma 4.4]{KMT-arxiv-2}.

\begin{lemma} \label{unifdecom}
There exist functions 
$\kappa \in \calS(\R^n)$ and $\chi \in \calS(\R^n)$ 
satisfying that
$\supp \kappa \subset [-1,1]^n$, 
$\supp \widehat \chi \subset B_1$, 
$| \chi | \geq c > 0$ on $[-1,1]^n$,
and
\begin{equation*}
\sum_{\nu\in\Z^n} \kappa (\xi - \nu) \chi (\xi - \nu ) = 1, \quad \xi \in \R^n.
\end{equation*}
\end{lemma}

As a first step of this section, 
we decompose symbols by
a Littlewood--Paley partition of unity 
$\{ \Psi_{j} \}_{ j \in \N_0 }$ on $(\R^{n})^2$
as follows:
\begin{equation*}
\sigma ( x, \xi, \eta )
=
\sum_{j \in \N_0} \sigma ( x, \xi, \eta ) \Psi_j ( \xi, \eta )
=\sum_{j\in\N_0} \sigma_j ( x, \xi, \eta )
\end{equation*}
with
\begin{equation} \label{symbolj}
\sigma_j (x,\xi,\eta) =\sigma ( x, \xi, \eta ) \Psi_j ( \xi, \eta ).
\end{equation}

Here, we observe the following two identities for $\Psi_{j}$.
Firstly, we let a function $\phi \in \calS(\R^n)$ satisfy that
$\phi = 1$ on $\{ \zeta \in \R^n : |\zeta| \leq 2 \}$ and
$\supp \phi \subset \{ \zeta \in \R^n : |\zeta| \leq 4 \}$,
and write $\phi_{j} = \phi (\cdot/2^{j})$.
Then, we have
for $j \geq 0$
\begin{equation*}
\Psi_j (\xi,\eta) = \Psi_j (\xi,\eta) \, \phi_j(\xi) \phi_j(\eta).
\end{equation*}
Secondly, there exist functions
$\phi' , \psi' , \psi'' \in \calS (\R^n) $
satisfying that
for $j \geq 1$
\begin{equation*}
\Psi_j ( \xi, \eta )
=
\Psi_j ( \xi, \eta )
\phi_{j}' (\xi) \psi_{j}' (\eta)
+
\Psi_j ( \xi, \eta ) 
\psi_{j}' (\xi) \phi_{j}' (\eta)
+
\Psi_j ( \xi, \eta ) 
\psi_{j}'' (\xi) \psi_{j}'' (\eta),
\end{equation*}
where $\phi_{j}' = \phi' (\cdot/2^{j})$, $\psi_{j}' = \psi' (\cdot/2^{j})$,
and $\psi_{j}'' = \psi'' (\cdot/2^{j})$,
and that
\begin{align} \label{suppdecomphi'}
\supp \phi' 
&\subset \left\{ \zeta \in\R^n : |\zeta| \leq 2^{-5}\right\},
\quad
\\ \label{suppdecompsi'}
\supp \psi'
&\subset \left\{ \zeta \in\R^n : 2^{-4} \leq |\zeta| \leq 2^{2}\right\},
\\ \label{suppdecompsi''}
\supp \psi''
&\subset \left\{ \zeta \in\R^n : 2^{-6} \leq |\zeta| \leq 2^{2}\right\}.
\end{align}
See Appendix \ref{appexistdecom} 
for the existence of such functions.
(By this, the annulus of $\supp \Psi_{j}$ is decomposed into the three parts:
$|\xi| \ll |\eta| $, $|\xi| \gg |\eta|$, $|\xi| \approx |\eta|$.)
Then, we have
\begin{align*}
\sigma ( x, \xi, \eta )
&=
\sum_{ j\lesssim1 }
\sigma_{j} ( x, \xi, \eta ) 
\, \phi_j(\xi) \phi_j(\eta)
+
\sum_{ j\gg1 }
\sigma_{j} ( x, \xi, \eta ) 
\, \phi_{j}' (\xi) \psi_{j}' (\eta)
\\
&\quad+
\sum_{ j\gg1 }
\sigma_{j} ( x, \xi, \eta ) 
\, \psi_{j}' (\xi) \phi_{j}' (\eta)
+
\sum_{ j\gg1 }
\sigma_{j} ( x, \xi, \eta ) 
\, \psi_{j}'' (\xi) \psi_{j}'' (\eta),
\end{align*}
where the implicit constant of the sum over $j$
depends only on $\rho$ and dimensions.


As a second step,
we rewrite the dual form 
of $T_\sigma (f,g)$.
By using the decomposition of symbols just above, 
the dual form of $T_\sigma (f,g)$ can be written as
\begin{align} \label{formIi}
\begin{split}
\int_{\R^n} T_\sigma (f,g) (x) \, h(x) \, dx
 &=
\sum_{ j\lesssim1 } \int_{\R^n} 
T_{\sigma_j} \big( \phi_{j}(D) f , \phi_{j}(D) g \big) (x) \, h(x) \, dx
\\&+
\sum_{ j\gg1} \int_{\R^n} 
T_{\sigma_j} \big( \phi_{j}' (D) f, \psi_{j}' (D) g \big) (x) \, h(x) \, dx
\\&+
\sum_{ j\gg1} \int_{\R^n} 
T_{\sigma_j} \big( \psi_{j}' (D) f, \phi_{j}' (D) g \big) (x) \, h(x) \, dx
\\&+
\sum_{ j\gg1} \int_{\R^n} 
T_{\sigma_j} \big( \psi_{j}'' (D) f, \psi_{j}'' (D) g \big) (x) \, h(x) \, dx
\\&
=: I_{0} + I_{1} + I_{2} + I_{3}.
\end{split}
\end{align}
In what follows, we rewrite this forms $I_{i}$, $i=0,1,2,3$.

To do this,
we shall consider the following:
\begin{equation*}
I:=\int_{\R^n} T_{\sigma_j} (F,G) (x) \, h(x) \, dx,
\end{equation*}
where $\sigma_j$ is in \eqref{symbolj}. 
By changes of variables, we have
\begin{align*}
I&=
\frac{2^{-j\rho n}}{(2\pi)^{2n}} 
\int_{(\R^{n})^3} e^{i x \cdot (\xi + \eta) }
\sigma_j (2^{-j\rho} x, 2^{j\rho} \xi, 2^{j\rho} \eta) 
\widehat{F (2^{-j\rho} \cdot )} (\xi) 
\widehat{G(2^{-j\rho} \cdot )} (\eta) \, 
{h(2^{-j\rho} x) }
\, dX
\\
&=
2^{-j\rho n} \int_{\R^n} 
T_{\sigma_j^\rho}	(F_j, G_j )(x) \, {h_j(x)}
\, dx,
\end{align*}
where we simply wrote $dX = dx d\xi d\eta$,
\begin{align} \label{sigmajrho}
\sigma_j^\rho 
&= \sigma_j ( 2^{-j\rho} \cdot, 2^{j\rho} \cdot, 2^{j\rho} \cdot),
\\\nonumber
F_j 
&= F (2^{-j\rho} \cdot ),\quad G_j = G (2^{-j\rho} \cdot ),
\andd h_j=h(2^{-j\rho} \cdot),
\end{align}
with $\sigma_j$ in \eqref{symbolj}. 
We next decompose the symbol $\sigma_j^\rho$.
We use the product type operator $\Delta_{\boldsymbol{k}}$
defined in Definition \ref{sigmanorm12}
and then use Lemma \ref{unifdecom}
to have
\begin{align*}
&\sigma_j^\rho (x,\xi,\eta)
=
\sum_{ \boldsymbol{k} = (k_0,k_1,k_2) \in (\N_0)^3}
\Delta_{\boldsymbol{k}}
[ \sigma_j^\rho ] (x,\xi, \eta)
\\&=
\sum_{ \boldsymbol{k} \in (\N_0)^3}
\sum_{ (\nu_1, \nu_2) \in (\Z^{n})^2 }
\Delta_{\boldsymbol{k}}
[ \sigma_j^\rho ] (x,\xi, \eta)
\chi (\xi-\nu_1) \chi (\eta-\nu_2)
\,\kappa (\xi-\nu_1) \kappa (\eta-\nu_2).
\end{align*}
Then, by writing as
\begin{equation}\label{sigmajknu}
\sigma_{ j, \boldsymbol{k}, \boldsymbol{\nu}}^\rho (x,\xi, \eta)  
=
\Delta_{\boldsymbol{k}} [ \sigma_j^\rho ] (x,\xi, \eta)
\chi(\xi - \nu_1) \chi(\eta - \nu_2)
\end{equation}
for $\boldsymbol{k} = (k_0,k_1,k_2) \in (\N_0)^3$ and 
$\boldsymbol{\nu} =(\nu_1, \nu_2)\in(\Z^n)^2$,
we can see that
\begin{equation*}
T_{\sigma_j^\rho} (F_j, G_j )(x)=
\sum_{ \boldsymbol{k} \in (\N_0)^3} \sum_{ \boldsymbol{\nu} \in (\Z^{n})^2}
T_{ \sigma_{ j, \boldsymbol{k}, \boldsymbol{\nu}}^\rho } 
\left( \kappa (D-\nu_1) F_j, \kappa (D-\nu_2) G_j \right) (x) .
\end{equation*}
By denoting the Fourier multiplier operator
$\kappa (D-\nu_i)$ 
by $\square_{\nu_i}$, $i=1,2$,
we have
\begin{equation}
\label{dualformdec}
I=
\sum_{ \boldsymbol{k} \in (\N_0)^3}
\sum_{ \boldsymbol{\nu} \in (\Z^{n})^2} 
2^{-j\rho n}
\int_{\R^n} T_{ \sigma_{ j, \boldsymbol{k}, \boldsymbol{\nu}}^\rho } (\square_{\nu_1} F_j, \square_{\nu_2} G_j ) (x) 
\, {h_j(x)} \, dx.
\end{equation}

Now, we actually rewrite the form $I_{i}$ defined in \eqref{formIi}.
For $I_{0}$, substituting $F=\phi_j (D) f$ and $G=\phi_j(D) g$ 
into \eqref{dualformdec}, then we have
\begin{equation} \label{i0}
I_{0}= 
\sum_{ j\lesssim1 }
\sum_{\boldsymbol{k} \in (\N_0)^3}
\sum_{\boldsymbol{\nu} \in (\Z^{n})^2}
2^{-j\rho n}
\int_{\R^n} 
T_{ \sigma_{ j, \boldsymbol{k}, \boldsymbol{\nu}}^\rho }
	(\square_{\nu_1} f_{j}, \square_{\nu_2} g_{j} ) (x) \, {h_j(x)} \, dx,
\end{equation}
where $\sigma_{ j, \boldsymbol{k}, \boldsymbol{\nu}}^\rho$
is in \eqref{sigmajknu} with \eqref{sigmajrho} and \eqref{symbolj},
\begin{equation*}
f_j = \phi_j (D) f (2^{-j\rho} \cdot ) , \quad
g_j = \phi_j (D) g (2^{-j\rho} \cdot ), \andd
h_j = h (2^{-j\rho} \cdot ).
\end{equation*}
Also, since
$\supp \phi \subset \{ \zeta \in \R^n : |\zeta| \leq 4 \}$,
$f_j$ and $g_j$ satisfy that
\begin{align}
\label{suppfgjrho0}
\begin{split}
\supp \widehat{f_j}
&\subset \left\{ \xi \in\R^n : |\xi| \leq 2^{j(1-\rho)+2}\right\} ,
\\
\supp \widehat{g_j}
&\subset \left\{ \eta \in\R^n : |\eta| \leq 2^{j(1-\rho)+2}\right\} .
\end{split}
\end{align}
Substituting
$(F,G) = (\phi_{j}'(D) f , \psi_{j}'(D) g )$, 
$(\psi_{j}'(D) f , \phi_{j}'(D) g )$, and 
$ (\psi_{j}''(D) f , \psi_{j}''(D) g )$ 
into \eqref{dualformdec},
the corresponding forms 
$I_{1}$, $I_{2}$, and $I_{3}$ 
are respectively given as follows.
For $I_{1}$, we have
\begin{equation} \label{i1}
I_{1}= 
\sum_{ j\gg1 }
\sum_{\boldsymbol{k} \in (\N_0)^3}
\sum_{\boldsymbol{\nu} \in (\Z^{n})^2} 2^{-j\rho n}
\int_{\R^n} 
T_{ \sigma_{ j, \boldsymbol{k}, \boldsymbol{\nu}}^\rho }
	(\square_{\nu_1} f_{j}^{(1)}, \square_{\nu_2} g_{j}^{(1)} ) (x) \, {h_j(x)} \, dx,
\end{equation}
where
\begin{equation} \label{fgjrho'}
f_{j}^{(1)} = \phi_{j}'(D) f (2^{-j\rho} \cdot ) , \quad
g_{j}^{(1)} = \psi_{j}'(D) g (2^{-j\rho} \cdot ), \andd 
h_j = h(2^{-j\rho} \cdot ).
\end{equation}
Also, $f_{j}^{(1)}$ and $g_{j}^{(1)}$ satisfy 
from \eqref{suppdecomphi'} and \eqref{suppdecompsi'}
that
\begin{align} \label{suppfgjrho'}
\begin{split}
\supp \widehat{f_{j}^{(1)}}
&\subset \left\{ \xi \in\R^n : |\xi| \leq 2^{j(1-\rho)-5}\right\} ,
\\
\supp \widehat{g_{j}^{(1)}}
&\subset \left\{ \eta \in\R^n : 2^{j(1-\rho)-4} \leq |\eta| \leq 2^{j(1-\rho)+2}\right\}.
\end{split}
\end{align}
Moreover, in this case, it should be remarkable that
\begin{equation} \label{xi+eta'}
2^{j(1-\rho)-5} \leq |\xi+\eta| \leq 2^{j(1-\rho)+3},
\quad \textrm{if} \quad 
(\xi,\eta) \in \supp \widehat{f_{j}^{(1)}} \times \supp \widehat{g_{j}^{(1)}}.
\end{equation}

For $I_{2}$, we have
\begin{equation} \label{i2}
I_{2}= 
\sum_{ j\gg1 }
\sum_{\boldsymbol{k} \in (\N_0)^3}
\sum_{\boldsymbol{\nu} \in (\Z^{n})^2} 2^{-j\rho n}
\int_{\R^n} 
T_{ \sigma_{ j, \boldsymbol{k}, \boldsymbol{\nu}}^\rho }
	(\square_{\nu_1} f_{j}^{(2)}, \square_{\nu_2} g_{j}^{(2)} ) (x) \, {h_j(x)} \, dx,
\end{equation}
where 
\begin{equation} \label{fgjrho''}
f_{j}^{(2)} = \psi_{j}'(D) f (2^{-j\rho} \cdot ), \quad
g_{j}^{(2)} = \phi_{j}'(D) g (2^{-j\rho} \cdot ), \andd
h_j = h(2^{-j\rho} \cdot ).
\end{equation}
Also, $f_{j}^{(2)}$ and $g_{j}^{(2)}$ satisfy 
from \eqref{suppdecompsi'} and \eqref{suppdecomphi'} that
\begin{align} \label{suppfgjrho''}
\begin{split}
\supp \widehat{f_{j}^{(2)}}
&\subset \left\{ \xi \in\R^n : 2^{j(1-\rho)-4} \leq |\xi| \leq 2^{j(1-\rho)+2}\right\},
\\
\supp \widehat{g_{j}^{(2)}}
&\subset \left\{ \eta \in\R^n : |\eta| \leq 2^{j(1-\rho)-5}\right\}.
\end{split}
\end{align}
Again, in this case, it should be remarkable that 
\begin{equation} \label{xi+eta''}
2^{j(1-\rho)-5} \leq |\xi+\eta| \leq 2^{j(1-\rho)+3},
\quad \textrm{if} \quad 
(\xi,\eta) \in \supp \widehat{f_{j}^{(2)}} \times \supp \widehat{g_{j}^{(2)}}.
\end{equation}

For $I_{3}$, we have
\begin{equation} \label{i3}
I_{3}= 
\sum_{ j\gg1 }
\sum_{\boldsymbol{k} \in (\N_0)^3}
\sum_{\boldsymbol{\nu} \in (\Z^{n})^2} 2^{-j\rho n}
\int_{\R^n} 
T_{ \sigma_{ j, \boldsymbol{k}, \boldsymbol{\nu}}^\rho }
	(\square_{\nu_1} f_{j}^{(3)}, \square_{\nu_2} g_{j}^{(3)} ) (x) \, {h_j(x)} \, dx,
\end{equation}
where 
\begin{equation} \label{fgjrhodag}
f_{j}^{(3)} = \psi_{j}''(D) f (2^{-j\rho} \cdot ),
\quad
g_{j}^{(3)} = \psi_{j}''(D) g (2^{-j\rho} \cdot ), 
\andd
h_j = h(2^{-j\rho} \cdot ).
\end{equation}
Also, $f_{j}^{(3)}$ and $g_{j}^{(3)}$ satisfy 
from \eqref{suppdecompsi''} that
\begin{align} \label{suppfgjrho34}
\begin{split}
\supp \widehat{f_{j}^{(3)}},\;
\supp \widehat{g_{j}^{(3)}}
&\subset \left\{ \zeta \in\R^n : 2^{j(1-\rho)-6} \leq |\zeta| \leq 2^{j(1-\rho)+2}\right\} .
\end{split}
\end{align}
We remark that a fact like \eqref{xi+eta'} and \eqref{xi+eta''}
does not hold in this case.

In the forthcoming sections, we will estimate the rewritten dual forms
in \eqref{i0}, \eqref{i1}, \eqref{i2}, and \eqref{i3}
to show the boundedness in Theorems \ref{l2l2h1bdd} and \ref{l2bmol2bdd}.

\begin{remark}
Through the argument above,
we have been able to examine precisely the annulus 
$\{ |(\xi,\eta)| \approx 2^{j(1-\rho)} \}$
by $\boldsymbol{\nu} + [-1,1]^{2n}$, 
$\boldsymbol{\nu} \in \Z^n \times \Z^n$.
On the other hands,
in \cite[Section 3]{miyachi tomita 2018 AIFG},
Miyachi--Tomita employed the method examining
the annulus $\{|(\xi,\eta)| \approx 2^{j} \}$
by $2^{j\rho}\boldsymbol{\nu} + [-2^{j\rho},2^{j\rho}]^{2n}$.
These two are essentially the same,
and so the idea used in this section goes back 
to their method in \cite{miyachi tomita 2018 AIFG}.
Moreover, the essential idea of
decomposing the annulus into the three parts, 
the sets $\{ |\xi| \ll |\eta| \}$ and $\{ |\xi| \gg |\eta| \}$ 
to have $|\xi+\eta| \gg 1$
and the set $\{ |\xi| \approx |\eta| \}$, also
goes back to \cite{miyachi tomita 2018 AIFG}.

The changes of variables concerned with $2^{\pm j\rho}$
were used to show the boundedness on $L^2$
of linear pseudo-differential operators
with symbols in the exotic class
$S^{0}_{\rho,\rho}$, $0 \leq \rho < 1$.
See, e.g., \cite[Chapter VII, Section 2.5]{stein 1993}.

The idea of decomposing symbols 
by using the functions $\kappa$ and $\chi$
in Lemma \ref{unifdecom}
comes from Sugimoto \cite{sugimoto 1988 JMSJ}.
\end{remark}



\section{Proof of Theorem \ref{l2l2h1bdd}} \label{secl2l2h1}

In this section, we prove Theorem \ref{l2l2h1bdd}.
To this end, we show that 
the absolute values of $I_{i}$, 
given in 
\eqref{i0}, \eqref{i1}, \eqref{i2}, and \eqref{i3},
are bounded by a constant times
\begin{equation*}
\| \sigma \|_{BS^{m,\ast}_{\rho,\rho}(\boldsymbol{s}; \R^{n})} 
\| f \|_{L^2} \| g \|_{L^2} \| h \|_{bmo} 
\end{equation*}
with $m=-(1-\rho)n/2$ and 
$\boldsymbol{s} = (n/2+\varepsilon,n/2,n/2)$ 
for any $\varepsilon > 0$.
Then, these complete the proof 
of Theorem \ref{l2l2h1bdd}
by recalling the expression in \eqref{formIi}.
In what follows, 
we will first consider $I_{1}$.
However, we omit the proof for $I_{2}$
because of symmetry between $I_{2}$ and $I_{1}$.
After that, we mention about $I_{3}$ 
and finally  about $I_{0}$.
The basic idea contained in the proof goes back to 
Miyachi--Tomita \cite{miyachi tomita 2018 AIFG}.

Before the proof, we shall give two remarks.
In order to obtain the desired boundedness, 
we will apply Lemma \ref{nu12lemma}
to the dual forms $I_{i}$.
This means that this lemma will be used
under the setting
$\sigma_{\boldsymbol{\nu}} = 
\sigma_{ j, \boldsymbol{k}, \boldsymbol{\nu}}^\rho$.
Recall from \eqref{sigmajknu} and Lemma \ref{unifdecom} that
\begin{equation*}
\sigma_{ j, \boldsymbol{k}, \boldsymbol{\nu}}^\rho (x,\xi, \eta)  
= \Delta_{\boldsymbol{k}} [ \sigma_j^\rho ] (x,\xi, \eta) \, \chi(\xi - \nu_1) \chi(\eta - \nu_2)
\end{equation*}
with $\chi \in \calS(\R^n)$ 
satisfying that $\supp \widehat \chi \subset B_1$ and 
$| \chi | \geq c > 0$ on $[-1,1]^n$.

Firstly, let us investigate 
the support of the Fourier transform of 
$\sigma_{ j, \boldsymbol{k}, \boldsymbol{\nu}}^\rho$.
Since $\supp \widehat \chi \subset B_1$,
\begin{equation} \label{fouriersupportofsigma}
\supp \calF [\sigma_{ j, \boldsymbol{k}, \boldsymbol{\nu}}^\rho]
\subset B_{2^{k_0+1}} \times B_{2^{k_1+2}} \times B_{2^{k_2+2}}
\end{equation}
holds.
Hence, we will use Lemma \ref{nu12lemma}
with $R_0=2^{k_0+1}$ and $R_i =2^{k_i+2}$, $i=1,2$.

Secondly, we have
\begin{equation}
\label{symbolequivalence}
\big\| \sigma_{ j, \boldsymbol{k}, \boldsymbol{\nu}}^\rho (x,\xi,\eta) 
\big\|_{L^2_{\xi,\eta} L^2_{ul,x} \ell^\infty_{\boldsymbol{\nu}}}
\approx \| \Delta_{\boldsymbol{k}} [ \sigma_j^\rho ] \|_{L^{2}_{ul}((\R^n)^3)}.
\end{equation}
To see this, we separate the integrals of $L^2_{\xi,\eta}$ 
by using \eqref{cubicdiscretization}.
Then, we have
\begin{align*}
&
\big\| \sigma_{ j, \boldsymbol{k}, \boldsymbol{\nu} }^\rho (x+\nu_0,\xi,\eta) 
\big\|_{L^2_{\xi,\eta}L^2_{x}(Q)}^2
\\&=
\sum_{\mu_1,\mu_2\in\Z^n} \int_{Q^3} 
\big|\Delta_{\boldsymbol{k}} [ \sigma_j^\rho ](x+\nu_0,\xi+\mu_1,\eta+\mu_2) 
	\chi(\xi+\mu_1-\nu_1) \chi(\eta+\mu_2-\nu_2) \big|^2 \, dX
\\&\lesssim
\sum_{\mu_1,\mu_2\in\Z^n} 
	\prod_{i=1,2}(1+|\mu_i-\nu_i|)^{-(n+1)}
\int_{Q^3} 
\big|\Delta_{\boldsymbol{k}} [ \sigma_j^\rho ](x+\nu_0,\xi+\mu_1,\eta+\mu_2) \big|^2 \, dX
\end{align*}
for $\nu_0\in\Z^n$ and $\boldsymbol{\nu}=(\nu_1,\nu_2) \in (\Z^n)^2$,
where $dX = dx d\xi d\eta$.
Thus, it holds that
\begin{equation*}
\big\| 
\sigma_{ j, \boldsymbol{k}, \boldsymbol{\nu} }^{\rho} (x+\nu_0,\xi,\eta) 
\big\|_{L^2_{\xi,\eta}L^2_{x}(Q)}
\lesssim 
\| \Delta_{\boldsymbol{k}} [ \sigma_j^\rho ] \|_{L^{2}_{ul}}
\end{equation*}
for $\nu_0, \nu_1, \nu_2\in\Z^n$,
which yields the inequality $\lesssim$ of \eqref{symbolequivalence}.
The opposite inequality $\gtrsim$ of \eqref{symbolequivalence} 
can be proved in a similar way by using the fact 
$|\chi| \geq c > 0$ on $[-1,1]^n$.
Therefore,
when we use Lemma \ref{nu12lemma}
under the situation
$\sigma_{\boldsymbol{\nu}} = 
\sigma_{ j, \boldsymbol{k}, \boldsymbol{\nu}}^\rho$,
the equivalence \eqref{symbolequivalence} allows us to replace
$\| \sigma_{ j, \boldsymbol{k}, \boldsymbol{\nu}}^\rho \| =
\| \sigma_{ j, \boldsymbol{k}, \boldsymbol{\nu}}^\rho (x,\xi,\eta) 
\|_{L^2_{\xi,\eta}L^2_{ul,x} \ell^\infty_{\boldsymbol{\nu}}}$ by 
$\| \Delta_{\boldsymbol{k}} [ \sigma_j^\rho ] \|_{L^{2}_{ul}}$.


\subsection{Estimate for $I_{1}$} \label{sech1T1}

In this subsection, we consider the dual form
$I_{1}$ in \eqref{i1}.
We decompose the factor of $f$
by a Littlewood--Paley partition $\{ \psi_{\ell} \}$ on $\R^n$
as
\[\square_{\nu_1} f_{j}^{(1)} 
= \sum_{\ell \in \N_0} \Delta_{\ell} \left[\square_{\nu_1} f_{j}^{(1)} \right] 
= \sum_{\ell \in \N_0} \square_{\nu_1} \Delta_{\ell} f_{j}^{(1)} .\]
Then, $I_{1}$ can be expressed by
\begin{equation*}
I_{1}
= 
\sum_{ j\gg1 }
\sum_{\boldsymbol{k} \in (\N_0)^3}
\sum_{\ell \in \N_0}
\sum_{\boldsymbol{\nu} \in (\Z^{n})^2}
I_{j,\boldsymbol{k},\ell,\boldsymbol{\nu}}^{(1)}
\end{equation*}
with
\begin{equation} \label{h1I1jklnu}
I_{j,\boldsymbol{k},\ell,\boldsymbol{\nu}}^{(1)}
= 
2^{-j\rho n}
\int_{\R^n} 
T_{ \sigma_{ j, \boldsymbol{k}, \boldsymbol{\nu}}^\rho }
	(\square_{\nu_1} \Delta_{\ell} f_{j}^{(1)}, \square_{\nu_2} g_{j}^{(1)} ) (x) \, {h_j(x)} \, dx.
\end{equation}
The sums over
$\ell$ and $\boldsymbol{\nu}$
are restricted
by the factors $\square_{\nu_1} \Delta_\ell f_{j}^{(1)}$ 
and $\square_{\nu_2} g_{j}^{(1)}$.
 Firstly, recalling the notation 
$\square_{\nu_i} = \kappa (D-\nu_i)$, $i=1,2$,
with $\kappa \in \calS(\R^n)$ such that
$\supp \kappa \subset [-1,1]^n$,
we have by \eqref{suppfgjrho'}
\begin{align}\label{h1T1nu12}
\nu_1 \in \Lambda_{1,\ell} 
= \left\{ \nu_1 \in \Z^n : |\nu_1| \lesssim 2^{\ell} \right\},
\quad
\nu_2 \in \Lambda_{2,j} 
= \left\{ \nu_2 \in \Z^n : |\nu_2| \lesssim 2^{j(1-\rho)} \right\}.
\end{align}
Secondly, from the factor $\Delta_\ell f_{j}^{(1)}$, we see that $2^{\ell-1} \leq 2^{j(1-\rho)-5}$,
which implies that 
\begin{equation} \label{h1T1ell}
\ell \leq j(1-\rho)-4 \leq j(1-\rho).
\end{equation}
Note that the set $\{\ell \in \N_0 : \ell \leq j(1-\rho)-4\}$ is not empty,
since $j \gg 1$.
Hence, 
\begin{equation*}
I_{1}
= 
\sum_{ j\gg1 }
\sum_{\boldsymbol{k} \in (\N_0)^3}
\sum_{\ell \in \N_0 :\, \ell \leq j(1-\rho) }
\sum_{\boldsymbol{\nu} \in \Lambda_{1,\ell} \times \Lambda_{2,j} }
I_{j,\boldsymbol{k},\ell,\boldsymbol{\nu}}^{(1)}.
\end{equation*}
For this expression,
we further separate the sum over $j$ as
\begin{align}\label{h1T1dec}
\begin{split}
I_{1}&=
\sum_{k_0 \in \N_0}
\Bigg\{ \sum_{ \substack{ j\gg 1: \\ j(1-\rho) \leq k_0 + L } } 
	+ \sum_{ \substack{ j\gg 1: \\ j(1-\rho) > k_0 + L } }  \Bigg\}
\sum_{k_1,k_2 \in \N_0} 
\sum_{ \substack{ \ell \in N_0 : \\ \ell \leq j(1-\rho) } }
\sum_{\boldsymbol{\nu} \in \Lambda_{1,\ell} \times \Lambda_{2,j} }
I_{j,\boldsymbol{k},\ell,\boldsymbol{\nu}}^{(1)}
\\&=:I^{(1,1)} + I^{(1,2)}
\end{split}
\end{align}
for some sufficiently large constant $L > 0$,
where $I_{j,\boldsymbol{k},\ell,\boldsymbol{\nu}}^{(1)}$
is in \eqref{h1I1jklnu}.


\subsubsection{Estimate of $I^{(1,1)}$ in \eqref{h1T1dec}} \label{sech1I11}

Firstly, we observe that
\begin{align*}
&
\calF \left[ T_{ \sigma_{ j, \boldsymbol{k}, \boldsymbol{\nu}}^\rho }
	(\square_{\nu_1} \Delta_{\ell} f_{j}^{(1)}, \square_{\nu_2} g_{j}^{(1)} ) \right](\zeta)
\\&=
\frac{1}{(2\pi)^{2n}}
\int_{ (\R^{n})^2 }
\calF_0 [\sigma_{ j, \boldsymbol{k}, \boldsymbol{\nu}}^\rho] \big( \zeta - (\xi + \eta ) , \xi , \eta \big) 
\widehat{ \square_{\nu_1} \Delta_{\ell} f_{j}^{(1)} }(\xi)
\widehat{ \square_{\nu_2} g_{j}^{(1)} }(\eta)
\, d\xi d \eta.
\end{align*}
Combining this with the fact 
$\supp \calF_{0} [\sigma_{j,\boldsymbol{k},\boldsymbol{\nu}}^{\rho}]
(\cdot,\xi,\eta) \subset B_{2^{k_0+1}}$
from \eqref{fouriersupportofsigma}, we have
\begin{align} \label{suppFouTsigma}
\begin{split}
&\supp 
\calF \left[ T_{ \sigma_{j,\boldsymbol{k},\boldsymbol{\nu}}^{\rho} } 
	(\square_{\nu_1} \Delta_{\ell} f_{j}^{(1)}, \square_{\nu_2}g_{j}^{(1)}) \right]
\\&\subset 
\left\{ \zeta \in \R^n : | \zeta-(\xi+\eta) | \leq2^{k_0+1}, \;
\xi \in \supp \widehat{ \square_{\nu_1} \Delta_{\ell} f_{j}^{(1)} }, \;
\eta \in \supp \widehat{ \square_{\nu_2}g_{j}^{(1)} }
\right\}.
\end{split}
\end{align}
Here, since we are considering the sum over 
$j$ such that $j(1-\rho) \leq k_0 + L$, 
we have $|\xi|, |\eta| \lesssim 2^{k_0}$ for
$\xi \in 
\supp \widehat{ \square_{\nu_1} \Delta_{\ell} f_{j}^{(1)} } $
and
$\eta \in 
\supp \widehat{ \square_{\nu_2}g_{j}^{(1)} }$.
Hence, by \eqref{suppFouTsigma}
\begin{equation*}
\supp 
\calF \left[ T_{ \sigma_{j,\boldsymbol{k},\boldsymbol{\nu}}^\rho } 
	(\square_{\nu_1} \Delta_{\ell} f_{j}^{(1)}, \square_{\nu_2}g_{j}^{(1)}) \right]
\subset 
\big\{ \zeta \in \R^n : | \zeta | \lesssim 2^{k_0} \big\}.
\end{equation*}
We take a function $\varphi \in \calS(\R^n)$ such that $\varphi = 1$ 
on $\{ \zeta\in\R^n: |\zeta|\lesssim 1\}$. 
Then, $I^{(1,1)}$
can be rewritten as
\begin{align*}
I^{(1,1)}
&=
\sum_{k_0 \in \N_0}
\sum_{ j : \, j(1-\rho) \leq k_0 + L }
\sum_{k_1,k_2 \in \N_0} 
\sum_{ \ell : \, \ell \leq j(1-\rho) }
\sum_{ \boldsymbol{\nu} \in \Lambda_{1,\ell} \times \Lambda_{2,j} }
\\&\quad\times
2^{-j\rho n} 
\int_{\R^n} 
T_{ \sigma_{ j, \boldsymbol{k}, \boldsymbol{\nu}}^\rho }
	(\square_{\nu_1} \Delta_{\ell} f_{j}^{(1)}, \square_{\nu_2} g_{j}^{(1)} ) (x) 
	\, \varphi (D/2^{k_0})[ h_{j} ] (x)  \, dx.
\end{align*}

In what follows, we shall estimate this rewritten $I^{(1,1)}$.
Again, recall that $\square_{\nu_i} = \kappa (D-\nu_i)$, $i=1,2$,
with $\kappa \in \calS(\R^n)$ such that
$\supp \kappa \subset [-1,1]^n$.
Then, since 
$\min( | \Lambda_{1,\ell} |, | \Lambda_{2,j} | ) \lesssim 2^{\ell n}$,
we have by Lemma \ref{nu12lemma} (1)
with $p=q=2$ and $r=\infty$
\begin{align*}&
\left| I^{(1,1)} \right|
\lesssim
\sum_{k_0 \in \N_0}
\sum_{ j : \, j(1-\rho) \leq k_0 + L }
\sum_{k_1,k_2 \in \N_0} 
\sum_{ \ell : \, \ell \leq j(1-\rho) }
\\&\times
2^{-j\rho n} \,
2^{ \ell n/2 }
\, 2^{ (k_0 + k_1 + k_2 ) n/2} 
\| \Delta_{\boldsymbol{k}} [ \sigma_j^\rho ] \|_{ L^{2}_{ul} }
\| \Delta_\ell f_{j}^{(1)} \|_{L^2}
\| g_{j}^{(1)} \|_{L^2}
\| \varphi (D/2^{k_0})[ h_{j} ] \|_{L^\infty}.
\end{align*}

Here, recall the notations of $f_{j}^{(1)}$, $g_{j}^{(1)}$, and $h_{j}$ 
from \eqref{fgjrho'}.
Then, we have
\begin{align} \label{h1I11fg}
\begin{split}
\| \Delta_\ell f_{j}^{(1)} \|_{L^2}
&=
\big\| \Delta_\ell \big[ \phi_{j}'(D) f (2^{-j\rho} \cdot ) \big] \big\|_{L^2}
\lesssim
2^{j\rho n/2 } \| \Delta_{\ell + j\rho} f \|_{L^2} ,
\\
\| g_{j}^{(1)} \|_{L^2}
&=
2^{j\rho n/2 } \| \psi_{j}'(D) g \|_{L^2},
\end{split}
\end{align}
where, $\Delta_{\ell + j\rho} = \psi_{\ell}(D/2^{j\rho})$ for $\ell \geq 0$.
Moreover,
since $j(1-\rho) \leq k_0 + L$,
by Lemma \ref{foumulopbmo} (1),
we have for any $\varepsilon>0$
\begin{align}\label{h1I11h}
\begin{split}
\| \varphi (D/2^{k_0})[ h_{j} ] \|_{L^\infty}
&=\| \varphi (D/2^{k_0+j\rho}) h (2^{-j\rho} \cdot) \|_{L^\infty}
\\&\lesssim (1+k_0+j\rho) \| h \|_{bmo}
\leq C_\varepsilon 2^{k_0 \varepsilon} \| h \|_{bmo}.
\end{split}
\end{align}
(The role of the condition $j(1-\rho) \leq k_0 + L$ is finished
by obtaining this estimate.)

Hence,
by denoting the Fourier multiplier operator 
$\psi_{j}'(D)$ by $\Delta_{j}'$,
we have
\begin{align}\label{h1I11beforeschur}
\begin{split}
&\left| I^{(1,1)} \right|
\lesssim
\sum_{k_0 \in \N_0}
\sum_{ j \gg 1 } 
\sum_{k_1,k_2 \in \N_0} 
\sum_{ \ell : \, \ell \leq j(1-\rho) }
\\&\quad\quad\quad\times
2^{ \ell n/2 }
\, 2^{ k_0 (n/2+\varepsilon)} 
\, 2^{ (k_1 + k_2 ) n/2} 
\| \Delta_{\boldsymbol{k}} [ \sigma_j^\rho ] \|_{ L^{2}_{ul} } 
\left\| \Delta_{\ell + j\rho} f \right\|_{L^2}
\| \Delta_{j}' g \|_{L^2}
\| h \|_{bmo}
\\&\leq
\| \sigma \|_{BS^{m,\ast}_{\rho,\rho}(\boldsymbol{s}; \R^{n})} 
\| h \|_{bmo}
\sum_{ j \gg 1 }
\sum_{ \ell:\, \ell \leq j(1-\rho) }
2^{-j(1-\rho)n/2}
\,2^{ \ell n/2 }
\| \Delta_{\ell + j\rho} f \|_{L^2}
\| \Delta_{j}' g \|_{L^2}
\\&=
\| \sigma \|_{BS^{m,\ast}_{\rho,\rho}(\boldsymbol{s}; \R^{n})} 
\| h \|_{bmo} \,
\II,
\end{split}
\end{align}
where $m=-(1-\rho)n/2$,
$\boldsymbol{s} = (n/2+\varepsilon,n/2,n/2)$
for any $\varepsilon>0$, and 
\begin{equation}\label{h1I11II}
\II:=
\sum_{ j \gg 1 }
\sum_{ 0 \leq \ell \leq j(1-\rho) }
2^{-j(1-\rho)n/2}
2^{ \ell n/2 }
\| \Delta_{\ell + j\rho} f \|_{L^2}
\| \Delta_{j}' g \|_{L^2}.
\end{equation}

In what follows, 
we shall show that
\begin{equation}\label{estofII}
\II \lesssim \| f \|_{L^2} \| g \|_{L^2}.
\end{equation}
To this end, we divide $\II$ 
into the two parts $\ell=0$ and $\ell \geq 1$.
We write
\begin{equation} \label{pfh1I1elldivide}
\II = \II_{\ell=0} + \II_{\ell \geq 1}
\end{equation}
with
\begin{align*}
\II_{\ell=0} &=
\sum_{ j \gg 1 } 2^{-j(1-\rho)n/2}
\| \psi_{0} (D/2^{j\rho}) f \|_{L^2} 
\| \Delta_{j}' g \|_{L^2},
\\
\II_{\ell \geq 1} &=
\sum_{ j \gg 1 }
\sum_{ 1 \leq \ell \leq j(1-\rho) }
2^{-j(1-\rho)n/2}
2^{ \ell n/2 }
\| \psi_{\ell}(D/2^{j\rho}) f \|_{L^2}
\| \Delta_{j}' g \|_{L^2}.
\end{align*}
The sum for $\ell=0$ is estimated as
\begin{equation}\label{h1I11IIell=0}
\II_{\ell=0}
\lesssim
\| f \|_{L^2} \| g \|_{L^2}
\sum_{ j \gg 1 } 2^{-j(1-\rho)n/2}
\approx
\| f \|_{L^2} \| g \|_{L^2},
\end{equation}
since $0 \leq \rho < 1$.
For the sum in the case $\ell\geq 1$,
we take a function 
$\psi^{\dag} \in \calS(\R^n)$ such that
$\psi^{\dag} = 1$ on 
$\{ \xi \in \R^n : 1/4 \leq |\xi| \leq 4 \}$ and 
$\supp \psi^{\dag} \subset
\{ \xi \in \R^n : 1/8 \leq |\xi| \leq 8 \}$.
Then, we realize that
$\psi^{\dag} ( \cdot / 2^{ \ell+[j\rho] } ) = 1$ on
$\supp \psi_{\ell} ( \cdot / 2^{ j\rho} )$,
since $[j\rho] \leq j \rho \leq [j\rho] + 1$.
Hence, by writing the operator 
$\Delta^{\dag}_{\ell+[j\rho]} = \psi^{\dag} ( D / 2^{ \ell+[j\rho] } )$,
we have
\begin{equation*}
\| \psi_{\ell}(D/2^{j\rho}) f \|_{L^2}
= \| \Delta^{\dag}_{\ell+[j\rho]} \psi_{\ell}(D/2^{j\rho}) f \|_{L^2}
\lesssim \| \Delta^{\dag}_{\ell+[j\rho]} f \|_{L^2}.
\end{equation*}
Therefore, we see that
\begin{equation*}
\II_{\ell \geq 1}
\lesssim
\sum_{ j \gg 1 }
\sum_{ 1 \leq \ell \leq j(1-\rho) }
2^{-j(1-\rho)n/2}
\, 2^{ \ell n/2 }
\| \Delta^{\dag}_{\ell+[j\rho]} f \|_{L^2}
\| \Delta_{j}' g \|_{L^2}.
\end{equation*}
By using the fact 
$2^{j\rho} \approx 2^{[j\rho]}$
and the translation as $\ell+[j\rho] \mapsto \ell'$, 
we have
\begin{equation*}
\II_{\ell \geq 1}
\lesssim
\sum_{ j \gg 1 }
\sum_{ 1 \leq \ell' \leq j }
2^{-j n/2}
\, 2^{ \ell' n/2 }
\| \Delta^{\dag}_{\ell'} f \|_{L^2}
\| \Delta_{j}' g \|_{L^2}.
\end{equation*}
Here, we verify that
\begin{equation*}
\sup_{j \in \N} \left\{
\sum_{ \ell' \in \N } \mathbf{1}_{\ell' \leq j} \, 2^{ - j n/2 } 2^{ \ell' n/2 }
\right\}
\approx 1
\andd
\sup_{\ell' \in \N} \left\{ 
\sum_{ j \in \N } \mathbf{1}_{\ell' \leq j} \, 2^{ \ell' n/2 } 2^{ - j n/2 } 
\right\}
\approx 1 .
\end{equation*}
Then, we obtain from Lemma \ref{schur}
\begin{equation}\label{h1I11IIellgeq0}
\II_{\ell \geq 1}
\lesssim
\big\| \Delta^{\dag}_{\ell'} f \big\|_{L^2\ell^2_{\ell'} (\N)}
\big\| \Delta_{j}' g \big\|_{L^2\ell^2_{j} (\N)}
\lesssim
\| f \|_{L^2} \| g \|_{L^2}.
\end{equation}
Therefore, from \eqref{pfh1I1elldivide},
\eqref{h1I11IIell=0}, and \eqref{h1I11IIellgeq0},
we obtain \eqref{estofII}.

Finally, collecting \eqref{h1I11beforeschur},
\eqref{h1I11II}, and \eqref{estofII},
we obtain
\begin{equation*}
\left| I^{(1,1)} \right|
\lesssim
\| \sigma \|_{BS^{m,\ast}_{\rho,\rho}(\boldsymbol{s}; \R^{n})} 
\| f \|_{L^2} \| g \|_{L^2} \| h \|_{bmo},
\end{equation*}
where $m=-(1-\rho)n/2$ and 
$\boldsymbol{s} = (n/2+\varepsilon,n/2,n/2)$
for any $\varepsilon>0$.


\subsubsection{Estimate for $I^{(1,2)}$ in \eqref{h1T1dec}} \label{sech1I12}

In this subsection,
since we are considering 
the sum over $j$ such that $j(1-\rho) > k_0 + L $
with some large $L > 0$,
we have by \eqref{suppFouTsigma}
\begin{align*}
&\supp 
\calF \left[ T_{ \sigma_{j,\boldsymbol{k},\boldsymbol{\nu}}^{\rho} } 
	(\square_{\nu_1} \Delta_{\ell} f_{j}^{(1)}, \square_{\nu_2}g_{j}^{(1)}) \right]
\\&\subset 
\left\{ \zeta \in \R^n : | \zeta-(\xi+\eta) | \leq2^{j(1-\rho)+1-L}, \;
\xi \in \supp \widehat{ \square_{\nu_1} \Delta_{\ell} f_{j}^{(1)} }, \;
\eta \in \supp \widehat{ \square_{\nu_2}g_{j}^{(1)} }
\right\}.
\end{align*}
Here, we recall from \eqref{xi+eta'} that 
if $\xi \in 
\supp \widehat{ f_{j}^{(1)} } $
and
$\eta \in 
\supp \widehat{ g_{j}^{(1)} }$, then
$2^{j(1-\rho)-5} \leq |\xi+\eta| \leq 2^{j(1-\rho)+3}$ holds.
Then, we see that
\begin{align} \label{suppfouT1}
\begin{split}
&
\supp 
\calF \left[ T_{ \sigma_{j,\boldsymbol{k},\boldsymbol{\nu}}^\rho } 
	(\square_{\nu_1} \Delta_{\ell} f_{j}^{(1)}, \square_{\nu_2}g_{j}^{(1)}) \right]
\\&\quad\subset 
\big\{ \zeta \in \R^n : 
2^{j(1-\rho)-6} \leq |\zeta| \leq 2^{j(1-\rho)+4} \big\}.
\end{split}
\end{align}
We take a function 
$\psi^{\ddag} \in \calS(\R^n)$ satisfying that 
$\psi^{\ddag} = 1$
on
$\{ \zeta \in \R^n : 2^{-6} \leq |\zeta | \leq 2^{4} \}$
and
$\supp \psi^{\ddag} \subset
\{ \zeta \in \R^n : 2^{-7} \leq | \zeta | \leq 2^5 \}$.
Then, $I^{(1,2)}$ can be rewritten as
\begin{align*}
I^{(1,2)} 
&=
\sum_{k_0 \in \N_0}
\sum_{ j : \,  j(1-\rho) > k_0 + L }
\sum_{k_1,k_2 \in \N_0} 
\sum_{ \ell : \, \ell \leq j(1-\rho) }
\sum_{ \boldsymbol{\nu} \in \Lambda_{1,\ell} \times \Lambda_{2,j} }
\\&\quad\times
2^{-j\rho n} 
\int_{\R^n} 
T_{ \sigma_{ j, \boldsymbol{k}, \boldsymbol{\nu}}^\rho }
	(\square_{\nu_1} \Delta_{\ell} f_{j}^{(1)}, \square_{\nu_2} g_{j}^{(1)} ) (x) \, 
	\psi^{\ddag} (D/2^{j(1-\rho)})[ h_{j} ] (x)  \, dx.
\end{align*}
(The role of the condition $j(1-\rho) > k_0 + L $
is finished by obtaining this expression.)

Now, we shall estimate the newly given $I^{(1,2)}$. Since
$\min( | \Lambda_{1,\ell} |, | \Lambda_{2,j} | ) \lesssim 2^{\ell n}$,
we have by using Lemma \ref{nu12lemma} (1)
with $p=q=2$ and $r=\infty$
\begin{align*}
&
\left| I^{(1,2)} \right|
\lesssim
\sum_{k_0 \in \N_0}
\sum_{ j \gg 1 } 
\sum_{k_1,k_2 \in \N_0} 
\sum_{ \ell : \, \ell \leq j(1-\rho) }
\\&\times
2^{-j\rho n} 
2^{ \ell n/2 }
\, 2^{ (k_0 + k_1 + k_2 ) n/2} 
\| \Delta_{\boldsymbol{k}} [ \sigma_j^\rho ] \|_{ L^{2}_{ul} }
\| \Delta_\ell f_{j}^{(1)} \|_{L^2}
\| g_{j}^{(1)} \|_{L^2}
\| \psi^{\ddag} (D/2^{j(1-\rho)})[ h_{j} ] \|_{L^\infty}.
\end{align*}
Here, for the factors of $f$ and $g$,
we have the estimate \eqref{h1I11fg}.
For the factor of $h$,
since $\psi^{\ddag}(0)=0$,
we have
by Lemma \ref{foumulopbmo} (2)
\begin{equation*}
\| \psi^{\ddag} (D/2^{j(1-\rho)})[ h_{j} ] \|_{L^\infty}
=\| \psi^{\ddag} (D/2^{j}) h (2^{-j\rho} \cdot) \|_{L^\infty}
\lesssim \| h \|_{BMO}.
\end{equation*}
Therefore, by writing as $\Delta_{j}' = \psi_{j}'(D)$, we have
\begin{align*}
&\left| I^{(1,2)} \right|
\lesssim
\sum_{k_0 \in \N_0}
\sum_{ j \gg 1 }
\sum_{k_1,k_2 \in \N_0} 
\sum_{ \ell : \, \ell \leq j(1-\rho) }
\\&\quad\quad\quad\quad\times
2^{ \ell n/2 }
\, 2^{ (k_0 + k_1 + k_2 ) n/2} 
\| \Delta_{\boldsymbol{k}} [ \sigma_j^\rho ] \|_{ L^{2}_{ul} } 
\| \Delta_{\ell + j\rho} f \|_{L^2}
\| \Delta_{j}' g \|_{L^2}
\| h \|_{BMO}
\\&\leq
\| \sigma \|_{BS^{m,\ast}_{\rho,\rho}(\boldsymbol{s}; \R^{n})} 
\| h \|_{BMO}
\sum_{ j \gg 1 }
\sum_{ \ell:\, \ell \leq j(1-\rho) }
2^{-j(1-\rho)n/2}
\,2^{ \ell n/2 }
\| \Delta_{\ell + j\rho} f \|_{L^2}
\| \Delta_{j}' g \|_{L^2},
\end{align*}
where $m=-(1-\rho)n/2$ and 
$\boldsymbol{s} = (n/2,n/2,n/2)$.
Since the sums over $j$ and $\ell$ are
exactly identical with $\II$ defined in \eqref{h1I11II},
we see 
from \eqref{estofII} that
\begin{equation*}
\left| I^{(1,2)} \right|
\lesssim
\| \sigma \|_{BS^{m,\ast}_{\rho,\rho}(\boldsymbol{s}; \R^{n})} 
\| f \|_{L^2} \| g \|_{L^2} \| h \|_{BMO},
\end{equation*}
where $m=-(1-\rho)n/2$ and 
$\boldsymbol{s} = (n/2,n/2,n/2)$.


\subsection{Estimate for $I_{2}$} \label{sech1T2}

In this subsection, 
we consider the dual form $I_{2}$ 
given in \eqref{i2}.
However, by comparing
\eqref{fgjrho''} with \eqref{fgjrho'},
\eqref{suppfgjrho''} with \eqref{suppfgjrho'}, and
\eqref{xi+eta''} with \eqref{xi+eta'},
we realize that $I_{2}$ and $I_{1}$
are in symmetrical positions.
Moreover, in this section, 
we are considering the boundedness 
on $L^2 \times L^2$.
Therefore, following the same lines as in Section \ref{sech1T1},
we obtain
\begin{equation*}
\left| I_{2} \right|
\lesssim
\| \sigma \|_{BS^{m,\ast}_{\rho,\rho}(\boldsymbol{s}; \R^{n})} 
\| f \|_{L^2} \| g \|_{L^2} \| h \|_{bmo}
\end{equation*}
with $m=-(1-\rho)n/2$ and 
$\boldsymbol{s} = (n/2+\varepsilon,n/2,n/2)$
for any $\varepsilon>0$.


\subsection{Estimate for $I_{3}$} \label{sech1T3}

In this subsection, we consider the dual form 
$I_{3}$ given in \eqref{i3}.
As in the previous subsections, we put
\begin{equation*}
I_{j,\boldsymbol{k},\boldsymbol{\nu}}^{(3)}
= 
2^{-j\rho n}
\int_{\R^n} 
T_{ \sigma_{ j, \boldsymbol{k}, \boldsymbol{\nu}}^\rho }
	(\square_{\nu_1} f_{j}^{(3)}, \square_{\nu_2} g_{j}^{(3)} ) (x) \, {h_j(x)} \, dx.
\end{equation*}
We separate the sum of $I_{3}$ into three parts
with slight changes of the way to sum over $\boldsymbol{\nu}$
as follows.
For some sufficiently large constants $N>0$ and $C>0$,
\begin{align}\label{h1T3}
\begin{split}
I_{3}
&= 
\sum_{k_0 \in \N_0}
\sum_{ \substack{ j \gg 1: \\ j(1-\rho) \leq k_0 + N } } 
\sum_{k_1,k_2 \in \N_0} 
\sum_{\boldsymbol{\nu} \in (\Z^{n})^2}
I_{j,\boldsymbol{k},\boldsymbol{\nu}}^{(3)}
\\&\quad+
\sum_{k_0 \in \N_0}
\sum_{ \substack{ j \gg 1: \\ j(1-\rho) > k_0 + N } }
\sum_{k_1,k_2 \in \N_0} 
\sum_{ \substack{ \tau \in \Z^{n}: \\ |\tau| \leq C 2^{k_0}} }
\sum_{ \substack{ \nu_1 \in \Z^n: \\ \nu_1 + \nu_2 = \tau} }
I_{j,\boldsymbol{k},\boldsymbol{\nu}}^{(3)}
\\
&\quad+
\sum_{k_0 \in \N_0}
\sum_{ \substack{ j \gg 1 : \\ j(1-\rho) > k_0 + N } }
\sum_{k_1,k_2 \in \N_0} 
\sum_{ \substack{ \tau \in \Z^{n}: \\ |\tau| > C 2^{k_0}} }
\sum_{ \substack{ \nu_1 \in \Z^n: \\ \nu_1 + \nu_2 = \tau} }
I_{j,\boldsymbol{k},\boldsymbol{\nu}}^{(3)}
\\&=:
I^{(3,1)} + I^{(3,2)} + I^{(3,3)}.
\end{split}
\end{align}


\subsubsection{Estimate for $I^{(3,1)}$ in \eqref{h1T3}}

By the factors $\square_{\nu_1} f_{j}^{(3)}$ 
and $\square_{\nu_2} g_{j}^{(3)}$
with \eqref{suppfgjrho34},
we have the restriction of 
the sum over $\boldsymbol{\nu}$:
\begin{equation*}
\nu_1, \nu_2 \in \Lambda_{j} = 
\left\{ \nu \in \Z^n : |\nu| \lesssim 2^{j(1-\rho)} \right\}.
\end{equation*}
We next observe that 
$|\xi|, |\eta|\lesssim2^{k_0}$ hold
for 
$\xi \in \widehat{ \square_{\nu_1} f_{j}^{(3)} }$ and 
$\eta \in \widehat{ \square_{\nu_2}g_{j}^{(3)} }$,
since $ j(1-\rho) \leq k_0 + N $.
Then, by referring to the argument around \eqref{suppFouTsigma},
we see that
\begin{equation*}
\supp 
\calF \left[ T_{ \sigma_{j,\boldsymbol{k},\boldsymbol{\nu}}^\rho } 
	(\square_{\nu_1} f_{j}^{(3)}, \square_{\nu_2}g_{j}^{(3)}) \right]
\subset \big\{ \zeta \in \R^n : | \zeta | \lesssim 2^{k_0} \big\}.
\end{equation*} 
We take a function $\varphi \in \calS(\R^n)$ such that 
$\varphi = 1$ on $\{ \zeta\in\R^n: |\zeta|\lesssim 1\}$.
Then, $I^{(3,1)}$ can be expressed as
\begin{align*}
I^{(3,1)}
&=
\sum_{k_0 \in \N_0}
\sum_{ j :\, j(1-\rho) \leq k_0 + N } 
\sum_{k_1,k_2 \in \N_0} 
\sum_{\boldsymbol{\nu} \in \Lambda_{j} \times \Lambda_{j}}
\\&\quad\times
2^{-j\rho n}
\int_{\R^n} 
T_{ \sigma_{ j, \boldsymbol{k}, \boldsymbol{\nu}}^\rho }
	(\square_{\nu_1} f_{j}^{(3)}, \square_{\nu_2} g_{j}^{(3)} ) (x) \, 
	\varphi (D/2^{k_0})[ h_{j} ] (x)  \, dx.
\end{align*}

By Lemma \ref{nu12lemma} (1) with 
$p=q=2$, $r=\infty$, and the fact
$| \Lambda_j | \lesssim 2^{j(1-\rho)n}$,
we have
\begin{align*}&
\left| I^{(3,1)} \right|
\lesssim
\sum_{k_0 \in \N_0}
\sum_{ j :\, j(1-\rho) \leq k_0 + N } 
\sum_{k_1,k_2 \in \N_0} 
\\&\times 
2^{-j \rho n} \, 2^{j(1-\rho)n/2} \, 2^{ (k_0 + k_1 + k_2 ) n/2} 
\| \Delta_{\boldsymbol{k}} [ \sigma_j^\rho ] \|_{ L^{2}_{ul} }
\| f_{j}^{(3)} \|_{L^2} 
\| g_{j}^{(3)} \|_{L^2} 
\| \varphi (D/2^{k_0})[ h_{j} ] \|_{L^\infty}.
\end{align*}
Here,
$\| f_{j}^{(3)} \|_{L^2}
= 2^{j\rho n/2 } \| \psi_{j}''(D) f \|_{L^2}$ and
$\| g_{j}^{(3)} \|_{L^2}
= 2^{j\rho n/2 } \| \psi_{j}''(D) g \|_{L^2}$
hold from \eqref{fgjrhodag}.
Also, since $j(1-\rho) \leq k_0 + N$,
$\| \varphi (D/2^{k_0})[ h_{j} ] \|_{L^\infty}
\leq C_\varepsilon 2^{k_0 \varepsilon} \| h \|_{bmo}$ holds
for $\varepsilon>0$
from \eqref{h1I11h}.
(Again, the role of the condition $j(1-\rho) \leq k_0 + N$ is finished
by obtaining this estimate.)
Hence, we obtain from the Cauchy--Schwarz inequality
\begin{align*}
&
\left| I^{(3,1)} \right|
\lesssim
\sum_{k_0 \in \N_0}
\sum_{ j \gg 1 } 
\sum_{k_1,k_2 \in \N_0} 
\\&\quad\times 
2^{j(1-\rho)n/2} \, 2^{ k_0 (n/2+\varepsilon)} \, 2^{ (k_1 + k_2 ) n/2} 
\| \Delta_{\boldsymbol{k}} [ \sigma_j^\rho ] \|_{ L^{2}_{ul} }
\| \psi_{j}''(D) f \|_{L^2} 
\| \psi_{j}''(D) g \|_{L^2}
\| h \|_{bmo}
\\&\leq
\| \sigma \|_{BS^{m,\ast}_{\rho,\rho}(\boldsymbol{s}; \R^{n})} 
\| h \|_{bmo}
\sum_{ j \gg 1} 
\| \psi_{j}''(D) f \|_{L^2} 
\| \psi_{j}''(D) g \|_{L^2}
\\ 
&\lesssim
\| \sigma \|_{BS^{m,\ast}_{\rho,\rho}(\boldsymbol{s}; \R^{n})} 
\| f \|_{L^2} 
\| g \|_{L^2} 
\| h \|_{bmo},
\end{align*}
where $m=-(1-\rho)n/2$ and 
$\boldsymbol{s} = (n/2+\varepsilon,n/2,n/2)$
for any $\varepsilon > 0$.


\subsubsection{Estimate for $I^{(3,2)}$ in \eqref{h1T3}} \label{sech1I32}

We here write $\Lambda_{k_0} = \{ \tau \in \Z^n: |\tau| \leq C 2^{k_0} \}$.

We shall consider the support of the Fourier transform of 
$T_{ \sigma_{j,\boldsymbol{\nu},\boldsymbol{k}}^\rho } (\cdots)$.
We observe that
$|\xi + \eta| \lesssim 2^{k_0}$ holds for
$(\xi,\eta) \in 
\supp \widehat{ \square_{\nu_1} f_{j}^{(3)} } \times
\supp \widehat{ \square_{\nu_2}g_{j}^{(3)} }$,
since $|\xi-\nu_1| \lesssim 1$ and $|\eta-\nu_2|\lesssim 1$
for this $(\xi,\eta)$ and
$\nu_1+\nu_2=\tau \in \Lambda_{k_0}$
in the sum of $I^{(3,2)}$.
Then, by referring to the argument around \eqref{suppFouTsigma}, 
\begin{equation*}
\supp 
\calF \left[ T_{ \sigma_{j,\boldsymbol{k},\boldsymbol{\nu}}^{\rho} } 
	(\square_{\nu_1} f_{j}^{(3)}, \square_{\nu_2}g_{j}^{(3)}) \right]
\subset 
\big\{ \zeta \in \R^n : | \zeta | \lesssim 2^{k_0} \big\}.
\end{equation*} 
Take a function $\varphi \in \calS(\R^n)$ such that 
$\varphi = 1$ on $\{ \zeta \in \R^n : |\zeta| \lesssim 1\}$.
Then, we have
\begin{align*}
I^{(3,2)}
&=
\sum_{k_0 \in \N_0}
\sum_{ j :\, j(1-\rho) > k_0 + N }
\sum_{k_1,k_2 \in \N_0} 
\sum_{ \tau \in \Lambda_{k_0} }
\sum_{ \nu_1 :\, \nu_1 + \nu_2 = \tau }
\\&\quad\times2^{-j \rho n}
\int_{\R^n} 
T_{ \sigma_{j,\boldsymbol{k},\boldsymbol{\nu}}^{\rho} } 
	(\square_{\nu_1} f_{j}^{(3)}, \square_{\nu_2} g_{j}^{(3)} ) (x)
	\, \varphi (D/2^{k_0})[ h_{j} ](x) \, dz.
\end{align*}

Now, we shall estimate this $I^{(3,2)}$.
By Lemma \ref{nu12lemma} (3) with 
$p=q=2$ and $r=\infty$, 
\begin{align*}
&
\left| I^{(3,2)} \right|
\lesssim
\sum_{k_0 \in \N_0}
\sum_{ j :\, j(1-\rho) > k_0 + N }
\sum_{k_1,k_2 \in \N_0} 
\\&\times
2^{-j \rho n}
\,2^{k_0 n/2}
\,2^{(k_1+k_2)n/2}
\| \Delta_{\boldsymbol{k}} [ \sigma_j^\rho ] \|_{ L^{2}_{ul} }
\| f_{j}^{(3)} \|_{L^2} 
\| g_{j}^{(3)} \|_{L^2}
\left\| \varphi (D/2^{k_0})[ h_{j} ] \right\|_{ \ell^2_{\tau}(\Lambda_{k_0}) L^{\infty} }.
\end{align*}
Here, recalling \eqref{fgjrhodag}, we have
$\| f_{j}^{(3)} \|_{L^2} \lesssim 2^{j\rho n/2 } \| f \|_{L^2}$ and 
$\| g_{j}^{(3)} \|_{L^2} \lesssim 2^{j\rho n/2 } \| g \|_{L^2}$.
Moreover, since $\varphi (D/2^{k_0})[ h_{j} ]$ 
is independent of $\tau$,
we have by Lemma \ref{foumulopbmo} (1)
\begin{align*}
&
\left\| \varphi (D/2^{k_0})[ h_{j} ] 
\right\|_{ \ell^2_{\tau}(\Lambda_{k_0}) L^{\infty} }
= | \Lambda_{k_0} |^{1/2} \| \varphi (D/2^{k_0})[ h_{j} ]  \|_{ L^\infty }
\\&
\lesssim 2^{k_0 n/2} (1+k_0+j\rho) \| h \|_{bmo}
\leq 2^{k_0 n/2} (1+j) \| h \|_{bmo},
\end{align*}
where we used the condition $ j(1-\rho) > k_0 + N$
in the last inequality.
From these, 
\begin{align*}
\left| I^{(3,2)} \right|
&\lesssim
\sum_{k_0 \in \N_0} \sum_{ j :\, j(1-\rho) > k_0 + N } \sum_{k_1,k_2 \in \N_0} 
\\&\quad\times 
(1+j) \,
2^{k_0 n/2}
\,2^{(k_0+k_1+k_2)n/2}
\| \Delta_{\boldsymbol{k}} [ \sigma_j^\rho ] \|_{ L^{2}_{ul} }
\| f \|_{L^2} 
\| g \|_{L^2} \| h \|_{bmo}
\\
&\leq
\sum_{k_0 \in \N_0}
\sup_{j\in\N_0} \Bigg\{
\sum_{k_1,k_2 \in \N_0} 2^{j(1-\rho)n/2} \, 2^{ (k_0 + k_1 + k_2 ) n/2} 
\| \Delta_{\boldsymbol{k}} [ \sigma_j^\rho ] \|_{ L^{2}_{ul} }
\Bigg\}
\\
&\quad\times 2^{ k_0 n/2 } 
\sum_{ j :\, j(1-\rho) > k_0 + N }
(1+j) 2^{-j(1-\rho)n/2} \,
\| f \|_{L^2} \| g \|_{L^2} \| h \|_{bmo}.
\end{align*}
Here, we have for $0< \varepsilon < n/2$
\begin{equation*}
\sum_{ j :\, j(1-\rho) > k_0 + N } (1+j) 2^{-j(1-\rho)n/2}
\leq 
C_{\varepsilon,\rho} 2^{ - k_0 (n/2 - \varepsilon) },\quad
k_0 \in \N_0,
\end{equation*}
since $0 \leq \rho < 1$.
Therefore, we obtain
\begin{equation*}
\left| I^{(3,2)} \right|
\lesssim
\| \sigma \|_{BS^{m,\ast}_{\rho,\rho}(\boldsymbol{s}; \R^{n})} 
\| f \|_{L^2} \| g \|_{L^2} \| h \|_{bmo},
\end{equation*}
where $m=-(1-\rho)n/2$ and 
$\boldsymbol{s} = (n/2+\varepsilon,n/2,n/2)$
for any $\varepsilon > 0$.


\subsubsection{Estimate for $I^{(3,3)}$ in \eqref{h1T3}}

The sum over $\tau$ is further restricted to
\begin{equation*}
\tau \in \Lambda_{j,k_0} = \{ \tau \in \Z^n : C 2^{k_0} < |\tau| \lesssim 2^{j(1-\rho)} \},
\end{equation*}
by the factors $\square_{\nu_1} f_{j}^{(3)}$ 
and $\square_{\nu_2} g_{j}^{(3)}$.
Here, note that the set $\Lambda_{j,k_0}$ is not empty, since $j(1-\rho) > k_0 + N$
with a sufficiently large $N$.
Moreover, referring to the argument around \eqref{suppFouTsigma}, we have
\begin{equation*}
\supp 
\calF \left[ T_{ \sigma_{j,\boldsymbol{k},\boldsymbol{\nu}}^\rho } 
	(\square_{\nu_1} f_{j}^{(3)}, \square_{\nu_2}g_{j}^{(3)}) \right]
\subset 
\big\{ \zeta \in \R^n : | \zeta - \tau | \lesssim 2^{k_0} \big\} ,
\end{equation*} 
since $| \xi - \nu_1 | \lesssim 1$ and $| \eta - \nu_2| \lesssim 1$ hold
for $\xi \in 
\supp \widehat{ \square_{\nu_1} f_{j}^{(3)} } $
and
$\eta \in 
\supp \widehat{ \square_{\nu_2}g_{j}^{(3)} }$,
where $\nu_1 + \nu_2 = \tau$.
We take a function $\varphi \in \calS(\R^n)$ such that 
$\varphi = 1$ on $\{ \zeta \in \R^n : |\zeta| \lesssim 1\}$.
Then, $I^{(3,3)}$ can be rewritten by
\begin{align*}
I^{(3,3)}
&=
\sum_{k_0 \in \N_0}
\sum_{ j :\, j(1-\rho) > k_0 + N }
\sum_{k_1,k_2 \in \N_0} 
\sum_{ \tau \in \Lambda_{j,k_0} }
\sum_{ \nu_1 :\, \nu_1 + \nu_2 = \tau }
\\&\quad\times2^{-j \rho n}
\int_{\R^n} 
T_{ \sigma_{j,\boldsymbol{k},\boldsymbol{\nu}}^\rho } 
	(\square_{\nu_1} f_{j}^{(3)}, \square_{\nu_2} g_{j}^{(3)} ) (x)
	\, \varphi \left(\frac{D+\tau}{2^{k_0}}\right)[h_{j}](x)
	\, dx.
\end{align*}

Now, we shall estimate this $I^{(3,3)}$.
By Lemma \ref{nu12lemma} (3) with 
$p=q=2$ and $r=\infty$,
\begin{align*}
&
\left| I^{(3,3)} \right|
\lesssim
\sum_{k_0 \in \N_0}
\sum_{ j :\, j(1-\rho) > k_0 + N }
\sum_{k_1,k_2 \in \N_0} 
\\&\times
2^{-j \rho n}
2^{j(1-\rho)n/2}
2^{(k_1+k_2)n/2}
\| \Delta_{\boldsymbol{k}} [ \sigma_j^\rho ] \|_{ L^{2}_{ul} }
\| f_{j}^{(3)} \|_{L^2} \| g_{j}^{(3)} \|_{L^2}
\left\| \varphi \left(\frac{D+\tau}{2^{k_0}}\right)
[h_{j}] \right\|_{\ell^2_{\tau}(\Lambda_{j,k_0}) L^{\infty}} .
\end{align*}
Here, from \eqref{fgjrhodag}, 
$\| f_{j}^{(3)} \|_{L^2} = 2^{j \rho n/2} \| \psi_{j}''(D) f \|_{L^2}$ and 
$\| g_{j}^{(3)} \|_{L^2} = 2^{j \rho n/2} \| \psi_{j}''(D) g \|_{L^2}$.
Moreover, since 
$\varphi (\tau/2^{k_0}) = 0$ for $\tau \in \Lambda_{j,k_0}$,
we see from Lemma \ref{estSR} (2) 
that
\begin{equation*}
\left\| \varphi \left(\frac{D+\tau}{2^{k_0}}\right)[h_{j}] \right\|_{\ell^2_{\tau}(\Lambda_{j,k_0}) L^{\infty}}
\lesssim
2^{k_0 n/2} \| h( 2^{-j\rho} \cdot) \|_{BMO}
=
2^{k_0 n/2} \| h \|_{BMO}.
\end{equation*}
Here, it is remarkable that 
$BMO$ is scaling invariant,
that is, $\| f(\lambda \cdot) \|_{BMO} = \| f \|_{BMO}$
for $\lambda > 0$,
although the space $bmo$ is not so in general.
See, e.g., \cite[Proposition 3.1.2 (6)]{grafakos 2014 modern}.
Hence, we obtain from the Cauchy--Schwarz inequality
\begin{align*}
&\left| I^{(3,3)} \right|
\lesssim
\sum_{k_0 \in \N_0}
\sum_{ j \gg 1 } 
\sum_{k_1,k_2 \in \N_0} 
\\&\quad\times 
2^{j(1-\rho)n/2}
\, 2^{(k_0+k_1+k_2)n/2}
\| \Delta_{\boldsymbol{k}} [ \sigma_j^\rho ] \|_{ L^{2}_{ul} }
\| \psi_{j}''(D) f \|_{L^2} \| \psi_{j}''(D) g \|_{L^2}
\| h \|_{BMO}
\\&\leq
\| \sigma \|_{BS^{m,\ast}_{\rho,\rho}(\boldsymbol{s}; \R^{n})} 
\| h \|_{BMO}
\sum_{ j \gg 1 }
\| \psi_{j}''(D) f \|_{L^2} \| \psi_{j}''(D) g \|_{L^2}
\\&\lesssim
\| \sigma \|_{BS^{m,\ast}_{\rho,\rho}(\boldsymbol{s}; \R^{n})} 
\| f \|_{L^2} \| g \|_{L^2}
\| h \|_{BMO},
\end{align*}
where $m=-(1-\rho)n/2$ and $\boldsymbol{s} = (n/2,n/2,n/2)$.


\subsection{Estimate for $I_{0}$} \label{boundl1i0}

In this subsection, we consider $I_{0}$ in \eqref{i0}.
Considering $\square_{\nu_1} f_j$ and $\square_{\nu_1} g_j$, 
we see from \eqref{suppfgjrho0} that
$\nu_1, \nu_2 \in \Z^n$ satisfy that $| \nu_1 | \lesssim 1$ and $|\nu_2| \lesssim1$,
since $ j \lesssim 1$. 
Moreover, we see that 
$\supp 
\calF [ T_{ \sigma_{j,\boldsymbol{k},\boldsymbol{\nu}}^\rho } 
	(\square_{\nu_1} f_{j}, \square_{\nu_2}g_{j}) ]
\subset 
\{ \zeta \in \R^n : | \zeta | \lesssim 2^{k_0} \} $.
We take a function $\varphi \in \calS(\R^n)$ such that
$\varphi = 1$ on $\{ \zeta \in \R^n : |\zeta| \lesssim 1\}$.
Then, by Lemma \ref{nu12lemma} (2)
and Lemma \ref{foumulopbmo} (1),
we have
\begin{align*}
\left| I_{0} \right|
&\leq
\sum_{ j\lesssim1 }
\sum_{\boldsymbol{k} \in (\N_0)^3} 
\sum_{ |\nu_1| , |\nu_2| \lesssim 1 } 2^{-j\rho n}
\left| 
\int_{\R^n} 
T_{ \sigma_{j,\boldsymbol{k},\boldsymbol{\nu}}^\rho } 
	(\square_{\nu_1} f_j, \square_{\nu_2} g_j ) ( x )
	\, \varphi(D/2^{k_0}) [h_{j} ] (x)
	\, dx \right|
\\
&\lesssim
\sum_{ j\lesssim1 }
\sum_{\boldsymbol{k} \in (\N_0)^3}
2^{(k_1+k_2)n/2}
\| \Delta_{\boldsymbol{k}} [ \sigma_j^\rho ] \|_{ L^{2}_{ul} }
\| f_j \|_{L^2} \| g_j \|_{L^2} \| \varphi(D/2^{k_0}) [h_{j}] \|_{L^\infty}
\\
&\lesssim
\sum_{ j\lesssim1 } 
\sum_{\boldsymbol{k} \in (\N_0)^3}
(1+k_0)2^{(k_1+k_2)n/2}
\| \Delta_{\boldsymbol{k}} [ \sigma_j^\rho ] \|_{ L^{2}_{ul} }
\| f \|_{L^2} \| g \|_{L^2} \| h \|_{bmo}
\\
&\lesssim
\| \sigma \|_{BS^{m,\ast}_{\rho,\rho}(\boldsymbol{s}; \R^{n})} 
\| f \|_{L^2} \| g \|_{L^2} \| h \|_{bmo},
\end{align*}
where $m=-(1-\rho)n/2$ 
and $\boldsymbol{s} = (\varepsilon,n/2,n/2)$
for any $\varepsilon > 0$.


\begin{remark}
The conclusion requiring
the class $BS^{m,\ast}_{\rho,\rho}(\boldsymbol{s}; \R^{n})$
is only in Section \ref{sech1I32} .
The conclusions in the other subsections hold for
the wider class $BS^{m}_{\rho,\rho}(\boldsymbol{s}; \R^{n})$.
\end{remark}


\section{Proof of Theorem \ref{l2bmol2bdd}} \label{secl2bmol2}

In this section we prove Theorem \ref{l2bmol2bdd}.
To this end, we will estimate the dual forms $I_{i}$, $i=0,1,2,3$,
given in \eqref{i0}, \eqref{i1}, \eqref{i2}, and \eqref{i3},
by a constant times
\begin{equation*}
\| \sigma \|_{BS^{m,\ast}_{\rho,\rho}(\boldsymbol{s}; \R^{n})} 
\| f \|_{L^2} \| g \|_{bmo} \| h \|_{L^2} 
\end{equation*}
for $m=-(1-\rho)n/2$ and $\boldsymbol{s} = (n/2,n/2,n/2)$.
Then, these complete the proof 
of Theorem \ref{l2bmol2bdd}.
The basic idea of the proof here again 
goes back to 
\cite{miyachi tomita 2018 AIFG}.

We will again use Lemma \ref{nu12lemma}
under the setting
$\sigma_{\boldsymbol{\nu}} = \sigma_{ j, \boldsymbol{k}, \boldsymbol{\nu} }^{\rho}$.
Hence, from the same reasons as stated in Section \ref{secl2l2h1},
we will use the lemma with
$R_0=2^{k_0+1}$, $R_i =2^{k_i+2}$, $i=1,2$,
and the equivalence \eqref{symbolequivalence}.
Now, we shall start the proof of Theorem \ref{l2bmol2bdd}.


\subsection{Estimate for $I_{1}$}

In this subsection, 
we consider the dual form
$I_{1}$ given in \eqref{i1}.
As was done in the previous subsections,
we divide the sum over $j$ as follows.
For some sufficiently large constant $L > 0$, 
\begin{align}\label{bmoT1dec}
\begin{split}
I_{1}
&=
\sum_{k_0 \in \N_0}
\Bigg\{ \sum_{ \substack{ j \gg 1: \\ j(1-\rho) \leq k_0 + L } } 
	+ \sum_{ \substack{ j \gg 1: \\ j(1-\rho) > k_0 + L } } \Bigg\}
\sum_{k_1,k_2 \in \N_0} 
\sum_{\boldsymbol{\nu} \in (\Z^n)^2 }
I_{j,\boldsymbol{k},\boldsymbol{\nu}}^{(1)}
\\&=:I^{(1,1)} + I^{(1,2)}
\end{split}
\end{align}
with
\begin{equation*}
I_{j,\boldsymbol{k},\boldsymbol{\nu}}^{(1)}
= 
2^{-j\rho n}
\int_{\R^n} 
T_{ \sigma_{ j, \boldsymbol{k}, \boldsymbol{\nu}}^\rho }
	(\square_{\nu_1} f_{j}^{(1)}, \square_{\nu_2} g_{j}^{(1)} ) (x) \, {h_j(x)} \, dx.
\end{equation*}


\subsubsection{Estimate for $I^{(1,1)}$ in \eqref{bmoT1dec}}
\label{secbmoI11}

The sum over $\boldsymbol{\nu}$ of $I^{(1,1)}$
is restricted to
\begin{equation*}
\nu_1, \nu_2 \in \Lambda_{j} 
= \left\{ \nu \in \Z^n : |\nu| \lesssim 2^{j(1-\rho)} \right\},
\end{equation*}
because of the factors $\square_{\nu_1} f_{j}^{(1)}$ 
and $\square_{\nu_2} g_{j}^{(1)}$ with \eqref{suppfgjrho'}.
Hence, $I^{(1,1)}$ is rewritten by
\begin{equation*}
I^{(1,1)}=
\sum_{k_0 \in \N_0}
\sum_{ j : \, j(1-\rho) \leq k_0 + L } 
\sum_{k_1,k_2 \in \N_0} 
\sum_{\boldsymbol{\nu} \in \Lambda_{j} \times \Lambda_{j} }
I_{j,\boldsymbol{k},\boldsymbol{\nu}}^{(1)}.
\end{equation*}
By Lemma \ref{nu12lemma} (2) with $p=r=2$ and $q=\infty$, we have
\begin{align*}&
\left| I^{(1,1)} \right|
\lesssim
\sum_{k_0 \in \N_0}
\sum_{ j : \, j(1-\rho) \leq k_0 + L } 
\sum_{k_1,k_2 \in \N_0} 
\\&\times
2^{ -j\rho n } \, 2^{j(1-\rho)n/2} \, 2^{j(1-\rho)n/2} \, 2^{(k_1+k_2)n/2}
\| \Delta_{\boldsymbol{k}} [ \sigma_j^\rho ] \|_{ L^{2}_{ul} }
\| f_{j}^{(1)} \|_{L^2} \| g_{j}^{(1)} \|_{L^\infty} \| h_{j} \|_{L^2}.
\end{align*}
Here,
$\| f_{j}^{(1)} \|_{L^{2}}
\lesssim 2^{j \rho n/2} \| f \|_{L^{2}}$
and $\| h_{j} \|_{L^{2}} = 2^{j \rho n/2} \| h \|_{L^{2}}$ hold.
Moreover, it holds from Lemma \ref{foumulopbmo} (2) that
\begin{equation}\label{bmoI11g}
\| g_{j}^{(1)} \|_{L^\infty}
= \| \psi_{j}'(D) g (2^{-j\rho} \cdot ) \|_{L^\infty}
\lesssim \| g \|_{BMO}.
\end{equation}
Hence, since
$\sum_{ j:\, j(1-\rho) \leq k_0 + L } 2^{j(1-\rho)n/2} 
\lesssim 2^{k_0 n/2}$,
we have
\begin{align*}
&
\left| I^{(1,1)} \right|
\lesssim
\| f \|_{L^2} \| g \|_{BMO} \| h \|_{L^2}
\\
&\times
\sum_{k_0 \in \N_0}
\sum_{ j :\, j(1-\rho) \leq k_0 + L }
2^{j(1-\rho)n/2}
\sup_{ j \in \N_0 }
\Bigg\{
\sum_{k_1,k_2 \in \N_0} 
2^{j(1-\rho)n/2} \,
2^{(k_1+k_2)n/2}
\| \Delta_{\boldsymbol{k}} [ \sigma_j^\rho ] \|_{ L^{2}_{ul} }
\Bigg\}
\\&\lesssim
\| \sigma \|_{BS^{m,\ast}_{\rho,\rho}(\boldsymbol{s}; \R^{n})} 
\| f \|_{L^2} \| g \|_{BMO} \| h \|_{L^2},
\end{align*}
where $m=-(1-\rho)n/2$ and $\boldsymbol{s}=(n/2,n/2,n/2)$.


\subsubsection{Estimate for $I^{(1,2)}$ in \eqref{bmoT1dec}}

We follow the same lines 
as in Section \ref{sech1T1}.

We use a Littlewood--Paley partition on $\R^n$
to decompose $\square_{\nu_1} f_{j}^{(1)} $ as
$\sum_{\ell} \square_{\nu_1} \Delta_{\ell} f_{j}^{(1)}$.
Then, the sum over $\boldsymbol{\nu}$
is restricted to
$\nu_1 \in \Lambda_{1,\ell}
= \{ \nu_1 \in \Z^n : |\nu_1| \lesssim 2^\ell \}$
and 
$\nu_2 \in \Lambda_{2,j} 
= \{ \nu_2 \in \Z^n : |\nu_2| \lesssim 2^{ j(1-\rho) } \}$
(see \eqref{h1T1nu12}), and the sum over $\ell$ is restricted to
$\ell \leq j(1-\rho)$
(see \eqref{h1T1ell}).
Moreover, recall \eqref{suppfouT1}
and take a function 
$\psi^{\ddag} \in \calS(\R^n)$ such that 
$\psi^{\ddag} = 1$ on
$\{ \zeta \in \R^n : 2^{-6} \leq |\zeta | \leq 2^{4} \}$ and
$\supp \psi^{\ddag} \subset
\{ \zeta \in \R^n : 2^{-7} \leq | \zeta | \leq 2^{5} \}$.
Then, $I^{(1,2)}$ can be expressed as
\begin{align*}
I^{(1,2)}
&= 
\sum_{k_0 \in \N_0}
\sum_{ j : \,  j(1-\rho) > k_0 + L }
\sum_{k_1,k_2 \in \N_0} 
\sum_{ \ell : \, \ell \leq j(1-\rho) }
\sum_{ \boldsymbol{\nu} \in \Lambda_{1,\ell} \times \Lambda_{2,j} }
\\&\quad\times
2^{-j\rho n}
\int_{\R^n} 
T_{ \sigma_{ j, \boldsymbol{k}, \boldsymbol{\nu}}^\rho }
	(\square_{\nu_1} \Delta_{\ell} f_{j}^{(1)}, \square_{\nu_2} g_{j}^{(1)} ) (x) \, 
	\psi^{\ddag} (D/2^{j(1-\rho)})[ h_{j} ] (x)  \, dx.
\end{align*}

We shall estimate this new $I^{(1,2)}$.
We use Lemma \ref{nu12lemma} (1) with 
$p=r=2$ and $q=\infty$ to the sum over $\boldsymbol{\nu}$, 
and use the fact
$\min( | \Lambda_{1,\ell} |, | \Lambda_{2,j} | ) \lesssim 2^{\ell n}$
to have
\begin{align*}
&
\left| I^{(1,2)} \right|
\lesssim
\sum_{k_0 \in \N_0}
\sum_{ j \gg 1 } 
\sum_{k_1,k_2 \in \N_0} 
\sum_{ \ell : \, \ell \leq j(1-\rho) }
\\&\times
2^{-j\rho n} 
2^{ \ell n/2 }
\, 2^{ (k_0 + k_1 + k_2 ) n/2} 
\| \Delta_{\boldsymbol{k}} [ \sigma_j^\rho ] \|_{ L^{2}_{ul} }
\| \Delta_\ell f_{j}^{(1)} \|_{L^2}
\| g_{j}^{(1)} \|_{L^{\infty}}
\| \psi^{\ddag} (D/2^{j(1-\rho)})[ h_{j} ] \|_{L^{2}}.
\end{align*}
Here, 
$\| \Delta_\ell f_{j}^{(1)} \|_{L^2} 
\lesssim 2^{j \rho n/2} \| \Delta_{\ell + j \rho} f \|_{L^{2}}$
and
$\| g_{j}^{(1)} \|_{L^{\infty}}
\lesssim \| g \|_{BMO}$ hold
(see \eqref{h1I11fg} and \eqref{bmoI11g},
respectively).
Also, 
$\| \psi^{\ddag} (D/2^{j(1-\rho)})[ h_{j} ] \|_{L^{2}} 
= 2^{j \rho n/2} \| \psi^{\ddag} (D/2^{j})h \|_{L^{2}}$.
Hence, by denoting $\psi^{\ddag} (D/2^{j})$ by $\Delta_{j}^{\ddag}$,
we have
\begin{align*}
\left| I^{(1,2)} \right|
&\lesssim
\| \sigma \|_{BS^{m,\ast}_{\rho,\rho}(\boldsymbol{s}; \R^{n})} 
\| g \|_{BMO}
\sum_{ j \gg 1 }
\sum_{ \ell:\, \ell \leq j(1-\rho) }
2^{-j(1-\rho)n/2}
\,2^{ \ell n/2 }
\| \Delta_{\ell + j\rho} f \|_{L^2}
\| \Delta_{j}^{\ddag} h \|_{L^{2}},
\\
&\lesssim
\| \sigma \|_{BS^{m,\ast}_{\rho,\rho}(\boldsymbol{s}; \R^{n})} 
\| f \|_{L^2} \| g \|_{BMO} \| h \|_{L^2},
\end{align*}
where $m=-(1-\rho)n/2$ and $\boldsymbol{s}=(n/2,n/2,n/2)$,
since the sums over $j$ and $\ell$ are identical with
$\II$ in \eqref{h1I11II}
whose factor $\| \Delta_{j}^{(1)} g \|_{L^2}$ is replaced 
by $\| \Delta_{j}^{\ddag} h \|_{L^{2}}$
and the inequality
$\| \Delta_{j}^{\ddag} h \|_{L^{2} \ell^{2}_{j}(\N)} \lesssim \| h \|_{L^2}$ holds.


\subsection{Estimate for $I_{2}$}

As in the previous sections, 
we divide the sum of $I_{2}$ defined in \eqref{i2} as follows.
For some sufficiently large constants $M > 0$ and $C>0$, 
\begin{align}\label{bmoT2dec}
\begin{split}
I_{2}
&=
\sum_{ j\gg1 } 
\sum_{\boldsymbol{k} \in (\N_0)^{3}} 
\sum_{ \substack{ \nu_1 \in \Z^n, \\ \nu_{2} \in \Z^n :\, |\nu_2| \leq C } }
I_{j,\boldsymbol{k},\boldsymbol{\nu}}^{(2)}
\\&\quad+ 
\sum_{k_0 \in \N_0} 
\sum_{ \substack{ j\gg1: \\ j(1-\rho) \leq k_0 + M } } 
\sum_{k_1,k_2 \in \N_0} 
\sum_{ \substack{ \nu_1 \in \Z^n, \\ \nu_{2} \in \Z^n :\, |\nu_2| > C } }
I_{j,\boldsymbol{k},\boldsymbol{\nu}}^{(2)}
\\&\quad+ 
\sum_{k_0 \in \N_0} 
\sum_{ \substack{ j\gg1 \\ j(1-\rho) > k_0 + M } }
\sum_{k_1,k_2 \in \N_0} 
\sum_{ \substack{ \nu_1 \in \Z^n, \\ \nu_{2} \in \Z^n :\, |\nu_2| > C } }
I_{j,\boldsymbol{k},\boldsymbol{\nu}}^{(2)}
\\&
=: I^{(2,1)} + I^{(2,2)} + I^{(2,3)}
\end{split}
\end{align}
with
\begin{equation*}
I_{j,\boldsymbol{k},\boldsymbol{\nu}}^{(2)}
= 2^{-j\rho n}
\int_{\R^n} 
T_{ \sigma_{ j, \boldsymbol{k}, \boldsymbol{\nu}}^\rho }
	(\square_{\nu_1} f_{j}^{(2)}, \square_{\nu_2} g_{j}^{(2)} ) (x) \, {h_j(x)} \, dx.
\end{equation*}


\subsubsection{Estimate for $I^{(2,1)}$ in \eqref{bmoT2dec}}

We write $\Lambda_{2} = \{ \nu_2 \in \Z^n : |\nu_2| \leq C \}$
and observe that $\min (|\Z^n|,|\Lambda_{2}|)\lesssim1$.
Then, by Lemma \ref{nu12lemma} (1)
with $p=r=2$ and $q=\infty$
\begin{equation*}
\left| I^{(2,1)} \right|
\lesssim
\sum_{ j\gg1 } 
\sum_{\boldsymbol{k} \in (\N_0)^{3}} 
2^{-j\rho n} \, 2^{ (k_0 + k_1 + k_2 ) n/2 }
\| \Delta_{\boldsymbol{k}} [ \sigma_j^\rho ] \|_{ L^{2}_{ul} }
\| f_{j}^{(2)} \|_{L^2} \| g_{j}^{(2)} \|_{L^\infty} \| h_{j} \|_{L^2}.
\end{equation*}
Here, we have
$\| f_{j}^{(2)} \|_{L^2} \lesssim 2^{j\rho n/2 } \| f \|_{L^2} $
and
$\| h_{j} \|_{L^{2}} = 2^{j \rho n/2} \| h \|_{L^{2}}$
by \eqref{fgjrho''},
and 
\begin{equation*}
\| g_{j}^{(2)} \|_{L^\infty}
= \| \phi_{j}^{(2)}(D) g (2^{-j\rho} \cdot) \|_{L^\infty}
\lesssim (1+j) \| g \|_{bmo} 
\end{equation*}
by Lemma \ref{foumulopbmo} (1).
Hence, recalling the assumption $0 \leq \rho < 1$, we obtain
\begin{align*}
\left| I^{(2,1)} \right|
&\lesssim
\| \sigma \|_{BS^{m,\ast}_{\rho,\rho}(\boldsymbol{s}; \R^{n})}
\| f \|_{L^2} \| g \|_{bmo} \| h \|_{L^2}
\sum_{ j \gg 1 } (1+j) \, 2^{-j(1-\rho) n/2 }
\\&\lesssim
\| \sigma \|_{BS^{m,\ast}_{\rho,\rho}(\boldsymbol{s}; \R^{n})}
 \| f \|_{L^2} \| g \|_{bmo} \| h \|_{L^2},
\end{align*} 
where $m=-(1-\rho)n/2$ and $\boldsymbol{s}=(n/2,n/2,n/2)$.


\subsubsection{Estimate for $I^{(2,2)}$ in \eqref{bmoT2dec}} \label{secbmoI22}

The sum over $\boldsymbol{\nu}$ is restricted to
\begin{align*}
\nu_1 \in 
\left\{ \nu_1 \in \Z^n : |\nu_1| \lesssim 2^{j(1-\rho)} \right\}
\andd
\nu_2 \in 
\left\{ \nu_2 \in \Z^n : C < |\nu_2| \lesssim 2^{j(1-\rho)} \right\}.
\end{align*}
Denote by $\Lambda_{2,j}$ 
the set to which $\nu_2$ belongs.
Then, we have by Lemma \ref{nu12lemma} (2)
with $p=r=2$ and $q=\infty$
\begin{align*}
\left| I^{(2,2)} \right|
\lesssim&
\sum_{k_0 \in \N_0} 
\sum_{ j: \, j(1-\rho) \leq k_0 + M }
\sum_{k_1,k_2 \in \N_0} 
2^{-j\rho n} \, 2^{j(1-\rho)n/2} \, 2^{j(1-\rho)n/2} \, 2^{ ( k_1 + k_2 ) n/2 }
\\
&\times
\| \Delta_{\boldsymbol{k}} [ \sigma_j^\rho ] \|_{ L^{2}_{ul} }
\| f_{j}^{(2)} \|_{L^2}
\| \square_{\nu_2} g_{j}^{(2)} \|_{ \ell^2_{\nu_2} (\Lambda_{2,j}) L^{\infty} }
\| h_{j} \|_{L^2}.
\end{align*}
Here, we have
$\| f_{j}^{(2)} \|_{L^2} \lesssim 2^{j\rho n/2 } \| f \|_{L^2} $ and
$\| h_{j} \|_{L^{2}} = 2^{j \rho n/2} \| h \|_{L^{2}}$.
Moreover, we have
\begin{equation}\label{fl2gbmocalcu}
\| \square_{\nu_2} g_{j}^{(2)} \|_{ \ell^2_{\nu_2} (\Lambda_{2,j}) L^{\infty} }
\lesssim \| g \|_{BMO}.
\end{equation}
To see this, we observe that
\begin{align*}
\square_{\nu_2} g_{j}^{(2)} (x) 
&= \kappa(D-\nu_2) \left[ \phi_{j}' (D) g (2^{-j\rho} \cdot) \right](x) 
\\&= \kappa(D-\nu_2) \phi' (D/2^{j(1-\rho)}) [ g (2^{-j\rho} \cdot) ](x)
\end{align*}
and also that
$\kappa (-\nu_2) = 0$ 
for $\nu_2 \in \Lambda_{2,j}$,
since $\supp \kappa \subset [-1,1]^n$.
Then, by Corollary \ref{estSRcor}
and the scaling invariance of the space $BMO$, we have
\[
\| \square_{\nu_2} g_{j}^{(2)} \|_{ \ell^2_{\nu_2} (\Lambda_{2,j}) L^{\infty} }
\lesssim
\| g (2^{-j\rho} \cdot) \|_{BMO}
= \| g \|_{BMO}.
\]
Therefore, since $\sum_{ j :\, j(1-\rho) \leq k_0 + M } 
2^{j(1-\rho)n/2} \lesssim 2^{k_0 n/2}$,
we obtain
\begin{align*}
&\left| I^{(2,2)} \right|
\lesssim
\| f \|_{L^2} \| g \|_{BMO} \| h \|_{L^2}
\\
&\times
\sum_{k_0 \in \N_0} \sum_{ j:\, j(1-\rho) \leq k_0 + M } 2^{j(1-\rho)n/2}
\sup_{j\in\N_0} \Bigg\{ \sum_{k_1,k_2 \in \N_0} 
2^{j(1-\rho)n/2} \, 2^{ ( k_1 + k_2 ) n/2 }
\| \Delta_{\boldsymbol{k}} [ \sigma_j^\rho ] \|_{ L^{2}_{ul} }
\Bigg\}
\\
&\lesssim
\| \sigma \|_{BS^{m,\ast}_{\rho,\rho}(\boldsymbol{s}; \R^{n})} 
\| f \|_{L^2} \| g \|_{BMO} \| h \|_{L^2},
\end{align*} 
where $m=-(1-\rho)n/2$ and $\boldsymbol{s}=(n/2,n/2,n/2)$.


\subsubsection{Estimate for $I^{(2,3)}$ in \eqref{bmoT2dec}}

The sum over $\boldsymbol{\nu}$ in $I^{(2,3)}$
is restricted to
\begin{align*}
\nu_1 \in \Lambda_{1,j} = 
\{ \nu_1 \in \Z^n : |\nu_1| \lesssim 2^{j(1-\rho)} \}, \quad
\nu_2 \in \Lambda_{2,j} = 
\{ \nu_2 \in \Z^n : C < |\nu_2| \lesssim 2^{j(1-\rho)} \}.
\end{align*}
Moreover,
recalling the argument around \eqref{suppfouT1},
in the sum over $j$ such that $j(1-\rho) > k_0 + M$,
we have
\begin{equation*}
\supp 
\calF \left[ T_{ \sigma_{j,\boldsymbol{k},\boldsymbol{\nu}}^\rho } 
	(\square_{\nu_1} f_{j}^{(2)}, \square_{\nu_2}g_{j}^{(2)}) \right]
\subset 
\big\{ \zeta \in \R^n : 
2^{j(1-\rho)-6} \leq |\zeta| \leq 2^{j(1-\rho)+4} \big\},
\end{equation*}
since
$2^{j(1-\rho)-5} \leq |\xi+\eta| \leq 2^{j(1-\rho)+3}$ holds
for 
$\xi \in 
\supp \widehat{ f_{j}^{(2)} } $
and
$\eta \in 
\supp \widehat{ g_{j}^{(2)} }$
(see \eqref{xi+eta''}).
Hence, by taking a function $\psi^{\ddag}  \in \calS(\R^n)$ such that
$\psi^{\ddag} = 1$ on
$\{ \zeta \in \R^n : 2^{-6} \leq | \zeta| \leq 2^{4} \}$ and
$\supp \psi^{\ddag} \subset 
\{ \zeta \in \R^n : 2^{-7} \leq | \zeta| \leq 2^{5} \}$,
we have
\begin{align*}
I^{(2,3)} &=
\sum_{k_0 \in \N_0} 
\sum_{ j :\, j(1-\rho) > k_0 + M }
\sum_{k_1,k_2 \in \N_0} 
\sum_{ \boldsymbol{\nu} \in \Lambda_{1,j} \times \Lambda_{2,j} }
\\&\quad\times
2^{-j\rho n}
\int_{\R^n} 
T_{ \sigma_{ j, \boldsymbol{k}, \boldsymbol{\nu}}^\rho }
	(\square_{\nu_1} f_{j}^{(2)}, \square_{\nu_2} g_{j}^{(2)} ) (x) \, 
	\psi^{\ddag} (D/2^{j(1-\rho)})[ h_{j} ] (x)  \, dx.
\end{align*}

Now, we shall estimate the rewritten $I^{(2,3)}$ above.
Using Lemma \ref{nu12lemma} with $p=r=2$ and $q=\infty$
and recalling the estimate \eqref{fl2gbmocalcu},
we have
\begin{align*}
\left| I^{(2,3)} \right|
&\lesssim
\sum_{k_0 \in \N_0} 
\sum_{ j \gg 1 } 
\sum_{k_1,k_2 \in \N_0} 
2^{- j \rho n} \, 2^{j(1-\rho)n/2} \, 2^{ ( k_0+k_1+k_2 ) n/2 }
\\&\quad\times
\| \Delta_{\boldsymbol{k}} [ \sigma_j^\rho ] \|_{ L^{2}_{ul} }
\| f_{j}^{(2)} \|_{L^2}
\| \square_{\nu_2} g_{j}^{(2)} \|_{ \ell^2_{\nu_2} (\Lambda_{2,j}) L^{\infty} }
\| \psi^{\ddag} (D/2^{j(1-\rho)})[ h_{j} ] \|_{L^2}
\\&\lesssim
\| \sigma \|_{BS^{m,\ast}_{\rho,\rho}(\boldsymbol{s}; \R^{n})} 
\| g \|_{BMO}
\sum_{ j \gg 1 }
\| \psi_{j}''(D) f \|_{L^2}
\| \psi^{\ddag} (D/2^{j}) h \|_{L^2}
\\&\lesssim
\| \sigma \|_{BS^{m,\ast}_{\rho,\rho}(\boldsymbol{s}; \R^{n})} 
\| f \|_{L^2} \| g \|_{BMO}\| h \|_{L^2},
\end{align*}
where $m=-(1-\rho)n/2$ and $\boldsymbol{s}=(n/2,n/2,n/2)$.


\subsection{Estimate for $I_{3}$} \label{secbmoT3}

As in Section \ref{sech1T3}, 
we split the sum into two parts
and slightly change
the way to sum over $\boldsymbol{\nu}$.
For some sufficiently large $N > 0$,
\begin{align}\label{bmoT3dec}
\begin{split}
I_{3}
&=
\sum_{k_0 \in \N_0}
\sum_{ \substack{ j \gg 1 : \\ j(1-\rho) \leq k_0 + N } } 
\sum_{k_1,k_2 \in \N_0} 
\sum_{\boldsymbol{\nu} \in (\Z^{n})^2}
I_{j,\boldsymbol{k},\boldsymbol{\nu}}^{(3)}
\\&\quad+
\sum_{k_0 \in \N_0}
\sum_{ \substack{ j \gg 1 : \\ j(1-\rho) > k_0 + N } }
\sum_{k_1,k_2 \in \N_0} 
\sum_{ \tau \in \Z^{n} }
\sum_{ \substack{ \nu_1 \in \Z^n: \\ \nu_1 + \nu_2 = \tau} }
I_{j,\boldsymbol{k},\boldsymbol{\nu}}^{(3)}
\\&
=:I^{(3,1)} + I^{(3,2)}
\end{split}
\end{align}
with
\begin{equation*}
I_{j,\boldsymbol{k},\boldsymbol{\nu}}^{(3)}
= 
2^{-j\rho n}
\int_{\R^n} 
T_{ \sigma_{ j, \boldsymbol{k}, \boldsymbol{\nu}}^\rho }
	(\square_{\nu_1} f_{j}^{(3)}, \square_{\nu_2} g_{j}^{(3)} ) (x) \, {h_j(x)} \, dx.
\end{equation*}


\subsubsection{Estimate of $I^{(3,1)}$ in \eqref{bmoT3dec}}

The sum over $\boldsymbol{\nu}$ is restricted to 
$\nu_1, \nu_2 \in\{ \nu \in \Z^n : |\nu| \lesssim 2^{j(1-\rho)} \}.$
Then, by Lemma \ref{nu12lemma} (2) with $p=r=2$ and $q=\infty$, 
we have
\begin{align*}&
\left| I^{(3,1)} \right|
\lesssim
\sum_{k_0 \in \N_0}
\sum_{ j  :\,  j(1-\rho) \leq k_0 + N } 
\sum_{k_1,k_2 \in \N_0} 
\\&\times
2^{-j\rho n} \, 2^{j(1-\rho) n/2 }  \, 2^{j(1-\rho) n/2 } \, 2^{(k_1+k_2)n/2}
\| \Delta_{\boldsymbol{k}} [ \sigma_j^\rho ] \|_{ L^{2}_{ul} }
\| f_{j}^{(3)} \|_{L^2} \| g_{j}^{(3)} \|_{L^\infty} \| h_{j} \|_{L^2}.
\end{align*}
Hence, since $\| g_{j}^{(3)} \|_{L^\infty}
\lesssim \| g \|_{BMO}$ from Lemma \ref{foumulopbmo} (2), 
we obtain
\begin{align*}
&\left| I^{(3,1)} \right| 
\lesssim
\| f \|_{L^2} \| g \|_{BMO} \| h \|_{L^2}
\\&\times
\sum_{k_0 \in \N_0} 
\sum_{ j  :\,  j(1-\rho) \leq k_0 + N } 
2^{j(1-\rho) n/2 }  
\sup_{j\in\N_0} \left\{ 
\sum_{k_1,k_2 \in \N_0} 
2^{j(1-\rho) n/2 } \, 2^{(k_1+k_2)n/2}
\| \Delta_{\boldsymbol{k}} [ \sigma_j^\rho ] \|_{ L^{2}_{ul} }
\right\}
\\&\lesssim
\| \sigma \|_{BS^{m,\ast}_{\rho,\rho}(\boldsymbol{s}; \R^{n})} 
\| f \|_{L^2} \| g \|_{BMO}\| h \|_{L^2} ,
\end{align*}
where $m=-(1-\rho)n/2$ and $\boldsymbol{s}=(n/2,n/2,n/2)$.


\subsubsection{Estimate of $I^{(3,2)}$ in \eqref{bmoT3dec}}

Let $\{ \psi_{\ell} \}_{\ell\in\N_0}$ 
be the Littlewood--Paley partition of unity on $\R^n$.
Then, it holds that
\begin{equation*}
\sum_{ \ell \in \N_{0} } \psi_\ell (\zeta / 2^{k_0}) = 1, \quad
\zeta \in \R^{n}, \quad
k_{0} \in \N_{0}.
\end{equation*}
By using this, we further decompose $I^{(3,2)}$ 
as follows:
\begin{equation*}
I^{(3,2)}
=
\sum_{k_0 \in \N_0}
\sum_{ j : \, j(1-\rho) > k_0 + N }
\sum_{k_1,k_2 \in \N_0} 
\sum_{ \ell \in \N_{0} } 
\sum_{ \tau \in \Z^{n} }
\sum_{ \nu_1: \, \nu_1 + \nu_2 = \tau}
I_{j,\boldsymbol{k},\ell,\boldsymbol{\nu}}^{(3)}
\end{equation*}
with
\begin{equation*}
I_{j,\boldsymbol{k},\ell,\boldsymbol{\nu}}^{(3)}
=
2^{- j \rho n} 
\int_{\R^n} 
T_{ \sigma_{ j,\boldsymbol{k},\boldsymbol{\nu} }^{\rho} } 
	(\square_{\nu_1} f_{j}^{(3)}, \square_{\nu_2} g_{j}^{(3)}) (x) 
	\, \psi_\ell (D/2^{k_0}) [ h_{j} ] (x) \, dx.
\end{equation*}
As in the previous sections, observe that
\begin{equation*}
\supp 
\calF \left[ T_{ \sigma_{j,\boldsymbol{k},\boldsymbol{\nu}}^\rho } 
	(\square_{\nu_1} f_{j}^{(3)}, \square_{\nu_2}g_{j}^{(3)}) \right]
\subset 
\left\{ \zeta \in \R^n : 
| \zeta - \tau | \lesssim 2^{k_0}
\right\}
\end{equation*}
and take a function $\varphi \in \calS(\R^n)$ such that 
$\varphi = 1$ on $\{ \zeta \in \R^n : |\zeta| \lesssim 1\}$.
Then,
$I_{j,\boldsymbol{k},\ell,\boldsymbol{\nu}}^{(3)}$
can be rewritten as
\begin{equation} \label{L2Ijkellnyu3}
I_{j,\boldsymbol{k},\ell,\boldsymbol{\nu}}^{(3)}
=
2^{- j \rho n} 
\int_{\R^n} 
T_{ \sigma_{ j,\boldsymbol{k},\boldsymbol{\nu} }^{\rho} } 
	(\square_{\nu_1} f_{j}^{(3)}, \square_{\nu_2} g_{j}^{(3)}) (x) 
	\, \varphi \left( \frac{D+\tau}{2^{k_0}} \right) \psi_\ell (D/2^{k_0}) [ h_{j} ] (x) \, dx.
\end{equation}

We next investigate restrictions of
$\ell$ and $\tau$
by considering the multipliers operated to $h_{j}$. 
Observe that $|\nu_1|, |\nu_2|\lesssim 2^{j(1-\rho)}$
with $\nu_2 = \tau - \nu_1$.
Then, since $j(1-\rho) > k_0 + N$,
\begin{align} \label{suppFTtauk0}
\supp 
\varphi \left( \frac{\cdot+\tau}{2^{k_0}} \right)
&\subset 
\left\{ \zeta \in \R^n : 
| \zeta+\tau | \lesssim 2^{k_0 }
\right\}
\\& \label{suppFTj}
\subset 
\left\{ \zeta \in \R^n : 
|\zeta| \lesssim 2^{j(1-\rho)}
\right\}.
\end{align}
Furthermore,
\begin{align} \label{suppDelhk0ell}
\begin{split}
\supp \psi_{0} (\cdot / 2^{k_0})
&\subset \left\{ \zeta \in \R^n : |\zeta| \leq 2^{k_0 + 1 } \right\}, 
\quad\textrm{if}\quad
\ell = 0,
\\
\supp \psi_\ell (\cdot / 2^{k_0})
&\subset \left\{ \zeta \in \R^n : 2^{k_0 + \ell - 1 } \leq |\zeta| \leq 2^{k_0 + \ell + 1 } \right\}, 
\quad\textrm{if}\quad
\ell \geq 1.
\end{split}
\end{align}
From \eqref{suppFTj} and \eqref{suppDelhk0ell},
the sum over $\ell$ is restricted to
\begin{equation*}
\ell \in \Omega_{j,k_0} =
\left\{ \ell \in \N_{0} : \ell \leq j(1-\rho) - k_0 + C \right\}
\end{equation*}
with a suitable constant $C>0$ depending on dimensions.
From \eqref{suppFTtauk0} and \eqref{suppDelhk0ell},
the sum over $\tau$ is restricted to
\begin{equation*}
\tau \in \Lambda_{k_0,\ell} = 
\left\{ \tau \in \Z^n : |\tau| \lesssim 2^{k_0 + \ell} \right\}.
\end{equation*}
Therefore, with $I_{j,\boldsymbol{k},\ell,\boldsymbol{\nu}}^{(3)}$ 
newly given in \eqref{L2Ijkellnyu3},
$I^{(3,2)}$ can be rewritten by
\begin{equation*}
I^{(3,2)}
=
\sum_{k_0 \in \N_0}
\sum_{ j : \, j(1-\rho) > k_0 + N }
\sum_{k_1,k_2 \in \N_0} 
\sum_{ \ell \in \Omega_{j,k_0} } 
\sum_{ \tau \in \Lambda_{k_0,\ell} }
\sum_{ \nu_1: \, \nu_1 + \nu_2 = \tau}
I_{j,\boldsymbol{k},\ell,\boldsymbol{\nu}}^{(3)}.
\end{equation*}

We shall estimate this new $I^{(3,2)}$.
By Lemma \ref{nu12lemma} (3) with 
$p=r=2$ and $q=\infty$
\begin{align*}&
\left| I^{(3,2)} \right|
\lesssim
\sum_{k_0 \in \N_0}
\sum_{ j : \, j(1-\rho) > k_0 + N }
\sum_{k_1,k_2 \in \N_0} 
\sum_{ \ell \in \Omega_{j,k_0} } 
2^{- j \rho n} \, 2^{(k_0 + \ell) n/2} \, 2^{ (k_1 + k_2) n/2 }
\\
&\quad\times
\| \Delta_{\boldsymbol{k}} [ \sigma_j^\rho ] \|_{ L^{2}_{ul} }
\| f_{j}^{(3)} \|_{L^2} \| g_{j}^{(3)} \|_{L^\infty} 
\left\| \varphi \left( \frac{D+\tau}{2^{k_0}} \right) 
	\psi_\ell (D/2^{k_0}) [ h_{j} ] \right\|_{ \ell^2_{\tau}(\Lambda_{k_0,\ell}) L^2 }.
\end{align*}
Here,
$\| g_{j}^{(3)} \|_{L^\infty} \lesssim \| g \|_{BMO}$ hold
from Lemma \ref{foumulopbmo} (2).
Moreover, we have by Lemma \ref{estSR} (1) and changes of variable
\begin{align*}
\left\| \varphi \left( \frac{D+\tau}{2^{k_0}} \right) 
	\psi_\ell (D/2^{k_0}) [ h_{j} ] \right\|_{ \ell^2_{\tau}(\Z^n) L^2 }
&\lesssim
2^{k_0 n/2}
\| \psi_\ell (D/2^{k_0}) [ h_{j} ] \|_{ L^2 }
\\&=
2^{k_0 n/2} \, 2^{j \rho n/2}
\| \psi_\ell (D/2^{k_0+j\rho}) h \|_{ L^2 }.
\end{align*}
Hence, by denoting the operators
$\psi_{j}''(D)$ by $\Delta_{j}''$
and $\psi_\ell (D/2^{k_0+j\rho})$ by $\Delta_{\ell+k_0+j\rho}$,
\begin{align}\label{bmoI32beforeschur}
\begin{split}
&
\left| I^{(3,2)} \right|
\lesssim
\sum_{k_0 \in \N_0}
\sum_{ j : \, j(1-\rho) > k_0 + N }
\sum_{k_1,k_2 \in \N_0} 
\sum_{ \ell \in \Omega_{j,k_0} } 
2^{(k_0 + \ell) n/2} \, 2^{ (k_0 + k_1 + k_2) n/2 }
\\
&\quad\times
\| \Delta_{\boldsymbol{k}} [ \sigma_j^\rho ] \|_{ L^{2}_{ul} }
\| \Delta_{j}'' f \|_{L^2} \| g \|_{BMO} 
\| \Delta_{\ell+k_0+j\rho} h \|_{ L^2 }
\\&\leq
 \| g \|_{BMO} 
\sum_{k_0 \in \N_0}
\sup_{j\in\N_0} \Bigg\{ \sum_{k_1,k_2 \in \N_0} 
2^{j(1-\rho) n/2} \, 2^{ (k_0+k_1+k_2) n/2 }
\| \Delta_{\boldsymbol{k}} [ \sigma_j^\rho ] \|_{ L^{2}_{ul} }
\Bigg\}
\\
&\quad\times 
\sum_{ j :\, j(1-\rho) > k_0 + N } 
\sum_{ \ell \in \Omega_{j,k_0} } 
2^{(k_0 + \ell) n/2} \, 2^{-j(1-\rho) n/2}
\| \Delta_{j}'' f \|_{L^2} 
\| \Delta_{\ell+k_0+j\rho} h \|_{ L^2 }
\\
&\leq
\| \sigma \|_{BS^{m,\ast}_{\rho,\rho}(\boldsymbol{s}; \R^{n})} 
\| g \|_{BMO}
\sup_{k_0 \in \N_{0} } \III_{k_0}
\end{split}
\end{align}
with
\begin{equation*}
\III_{k_0}=
\sum_{ j :\,  j(1-\rho) > k_0 + N } 
\sum_{ \ell \in \Omega_{j,k_0} } 
2^{(k_0 + \ell) n/2} \, 2^{-j(1-\rho) n/2}
\| \Delta_{j}'' f \|_{L^2} 
\| \Delta_{\ell+k_0+j\rho} h \|_{ L^2 }.
\end{equation*}

We shall estimate $\III_{k_0}$
and prove that this is bounded by 
a constant times $\| f \|_{L^2} \| h \|_{L^2}$
for all $k_0 \in \N_0$.
The way to estimate $\III_{k_0}$ is almost the same 
as was done to have \eqref{estofII}.
We divide $\III_{k_0}$ 
into the two parts $\ell=0$ and $\ell \geq 1$.
That is,
\begin{equation} \label{bmoI32ell}
\III_{k_0} = \III_{k_0}^{\ell=0} + \III_{k_0}^{\ell \geq 1}
\end{equation}
with
\begin{align*}
\III_{k_0}^{\ell=0} &=
\sum_{ j :\,  j(1-\rho) > k_0 + N } 
2^{k_0 n/2} \, 2^{-j(1-\rho) n/2}
\| \Delta_{j}'' f \|_{L^2} 
\| \psi_{0} (D/2^{k_0+j\rho}) h \|_{ L^2 },
\\
\III_{k_0}^{\ell\geq1} &=
\sum_{ j :\,  j(1-\rho) > k_0 + N }
\sum_{ \ell \in \N \cap \Omega_{j,k_0} }
2^{(k_0 + \ell) n/2} \, 2^{-j(1-\rho) n/2}
\| \Delta_{j}'' f \|_{L^2} 
\| \psi_\ell (D/2^{k_0+j\rho}) h \|_{ L^2 }.
\end{align*}
For the first sum $\III_{k_0}^{\ell=0}$, we have
\begin{equation}\label{bmoI32ell0}
\III_{k_0}^{\ell=0}
\lesssim
\| f \|_{L^2} 
\| h \|_{ L^2 }
\sum_{ j :\, j(1-\rho) > k_0 + N }
2^{k_0 n/2} \, 2^{-j(1-\rho) n/2}
\lesssim \| f \|_{L^2} \| h \|_{L^2}.
\end{equation}
We next consider the second sum $\III_{k_0}^{\ell\geq1}$.
We take a function $\psi^{\dag} \in \calS(\R^n)$ such that
$\psi^{\dag} = 1$ on
$\{ \xi \in \R^n : 1/4 \leq |\xi| \leq 4 \}$ and 
$\supp \psi^{\dag} \subset
\{ \xi \in \R^n : 1/8 \leq |\xi| \leq 8 \}$.
Then, by writing as
$\Delta^{\dag}_{\ell+k_0+[j\rho]} = \psi^{\dag} ( D / 2^{ \ell+k_0+[j\rho] } )$,
we have
\begin{equation*}
\| \psi_\ell (D/2^{k_0+j\rho}) h \|_{ L^2 } 
=
\| \Delta^{\dag}_{\ell+k_0+[j\rho]} \psi_\ell (D/2^{k_0+j\rho}) h \|_{ L^2 }
\lesssim
\big\| \Delta^{\dag}_{\ell + k_0 + [j\rho] } h \big\|_{L^2}.
\end{equation*}
From this and the fact $2^{j\rho} \approx 2^{[j\rho]}$,
it holds that
\begin{equation*}
\III_{k_0}^{\ell\geq1} \lesssim
\sum_{ j :\, j(1-\rho) > k_0 + N }
\sum_{ \ell : \, 1 \leq \ell \leq j(1-\rho) - k_0 + C }
2^{ (\ell + k_0 +[j\rho]) n/2 } \, 2^{-jn/2}
\| \Delta_{j}'' f \|_{L^2} \| \Delta^{\dag}_{\ell + k_0 + [j\rho] } h \|_{L^2}.
\end{equation*}
Changing the sum of $\ell$ as $\ell + k_0 + [j\rho] \mapsto \ell'$,
we have by Lemma \ref{schur}
\begin{equation} \label{bmoI32ellgeq1}
\III_{k_0}^{\ell\geq1} \lesssim
\sum_{ j \geq 1 } \sum_{ \ell' :\, 1 \leq \ell' \leq j+C }
2^{ \ell' n/2 } \, 2^{-j n/2}
\| \Delta_{j}'' f \|_{L^2} \| \Delta^{\dag}_{\ell'} h \|_{L^2}
\lesssim 
\| f \|_{L^2} \| h \|_{L^2}.
\end{equation}

Lastly, collecting 
\eqref{bmoI32beforeschur}, \eqref{bmoI32ell}, 
\eqref{bmoI32ell0}, and \eqref{bmoI32ellgeq1},
we obtain that
\begin{equation*}
\left| I^{(3,2)} \right|
\lesssim 
\| \sigma \|_{BS^{m,\ast}_{\rho,\rho}(\boldsymbol{s}; \R^{n})} 
\| f \|_{L^2} \| g \|_{BMO} \| h \|_{L^2},
\end{equation*}
where $m=-(1-\rho)n/2$ and $\boldsymbol{s}=(n/2,n/2,n/2)$.


\subsection{Estimate for $I_{0}$}

In this section, we consider $I_{0}$ in \eqref{i0}.
Considering $\square_{\nu_1} f_j$ and $\square_{\nu_1} g_j$, 
we see from \eqref{suppfgjrho0} that
$\nu_1, \nu_2 \in \Z^n$ satisfy that $|\nu_1|, |\nu_2| \lesssim1$,
since $ j \lesssim 1$. 
Hence, we have by Lemma \ref{nu12lemma} (2) with $p=r=2$, $r=\infty$
and Lemma \ref{foumulopbmo} (1),
\begin{align*}
\left| I_0 \right|
&\leq
\sum_{ j\lesssim1 }
\sum_{\boldsymbol{k} \in (\N_0)^3} 
\sum_{ |\nu_1| , |\nu_2| \lesssim 1 } 2^{-j\rho n}
\left| 
\int_{\R^n} 
T_{ \sigma_{j,\boldsymbol{k},\boldsymbol{\nu}}^\rho } 
	(\square_{\nu_1} f_j, \square_{\nu_2} g_j ) ( x )
	\, h_{j} (x)
	\, dx \right|
\\
&\lesssim
\sum_{ j\lesssim1 }
\sum_{ {\boldsymbol{k}} \in (\N_0)^3}
2^{-j\rho n} \, 2^{ (k_1 + k_2 ) n/2} 
\| \Delta_{\boldsymbol{k}} [ \sigma_j^\rho ] \|_{ L^{2}_{ul} }
\| f_j \|_{L^2} \| g_j \|_{L^\infty} \| h_{j} \|_{L^2}
\\
&\lesssim
\| \sigma \|_{BS^{m,\ast}_{\rho,\rho}(\boldsymbol{s}; \R^{n})} 
\| f \|_{L^2} \| g \|_{bmo} \| h \|_{L^2},
\end{align*}
where $m=-(1-\rho)n/2$ and $\boldsymbol{s} = (0,n/2,n/2)$.



\section{Sharpness of Theorem \ref{main-thm}} \label{secsharp}

In this section, we consider the sharpness of the conditions
of the order $m$ and the smoothness $\boldsymbol{s}=(s_0,s_1,s_2)$ 
stated in Theorem \ref{main-thm}.


\subsection{Sharpness of the order $m$}

In this subsection, we show the following.
\begin{proposition}
Let $0 \leq \rho <1$, $m \in \R$, 
$\boldsymbol{s}=(s_0,s_1,s_2) \in [0,\infty)^3$, 
and $1 \leq r \leq 2\leq p,q \leq \infty$ 
satisfy $1/p+1/q=1/r$.
If all bilinear pseudo-differential operators with symbols in 
$BS^m_{\rho,\rho}(\boldsymbol{s};\R^{n})$
are bounded from $L^p \times L^q$ to $L^r$,
then $m \leq -(1-\rho)n/2$.
\end{proposition}

This is immediately obtained by
the fact that $BS^m_{\rho,\rho}(\R^n) 
\subset BS^{m}_{\rho,\rho}(\boldsymbol{s}; \R^{n})$,
$\boldsymbol{s} \in [0,\infty)^3$,
from Lemma \ref{derivationlemma} and
the following theorem proved by Miyachi--Tomita 
\cite[Theorem A.2]{miyachi tomita 2013 IUMJ}.

\begin{theorem}
Let $0 \leq \rho <1$, $m \in \R$, and 
$1 \leq r \leq 2\leq p,q \leq \infty$ 
satisfy $1/p+1/q=1/r$.
If all bilinear pseudo-differential operators with symbols in $BS^m_{\rho,\rho}(\R^n)$
are bounded from $L^p \times L^q$ to $L^r$,
then $m \leq -(1-\rho)n/2$.
\end{theorem}


\subsection{Sharpness of the smoothness $s_1$ and $s_2$}

In this subsection, we show the following.
The idea of the proof comes from
Miyachi--Tomita \cite[Section 7]{miyachi tomita 2013 RMI}.

\begin{proposition}\label{rhosharps12}
Let $0\leq \rho<1$, $m \in \R$, $\boldsymbol{s}=(s_0,s_1,s_2)\in (0,\infty)^3$, 
and $1 \leq p,q,r \leq \infty$ satisfy $1/p+1/q=1/r$. 
Suppose that the inequality
\begin{equation}\label{rhoassumptionofs12}
\| T_\sigma (f,g) \|_{L^r}
\lesssim 
\| \sigma \|_{BS^{m}_{\rho,\rho}(\boldsymbol{s}; \R^{n})} \| f \|_{L^p} \| g \|_{L^q}
\end{equation}
holds for all $\sigma \in BS^{m}_{\rho,\rho}(\boldsymbol{s}; \R^{n})$ and $f,g\in\calS(\R^n)$.
Then, $s_1,s_2 \geq n/2$.
\end{proposition}

We will use the following fact.
See, e.g., \cite[Proposition 1.1 (i)]{sugimoto 1988 TMJ}.

\begin{lemma}\label{sugimoto}
Let $1 \leq p,q \leq \infty$ and $s>0$. Then, we have
\begin{equation*}
\| f(\lambda\cdot) \|_{B^s_{p,q}(\R^n)} 
\lesssim \lambda^{-n/p}\max(1,\lambda^s) \| f \|_{B^s_{p,q}(\R^n)}, \quad
\lambda > 0 .
\end{equation*}
\end{lemma}

\begin{proof}[Proof of Proposition \ref{rhosharps12}]
It suffices to prove $s_1 \geq n/2$.
Let $u,v \in \calS (\R^n)$ satisfy that
\begin{align*}
\supp \widehat u &\subset \{ \xi \in \R^n : | \xi | \leq 1\},
\\
\supp \widehat v &\subset \{ \eta \in \R^n : 9/10 \leq | \eta | \leq 11/10\}, \quad
\\
\widehat v &= 1 \quad \textrm{on} \quad
 \{ \eta \in \R^n : 19/20 \leq |\eta| \leq 21/20 \},
\end{align*}
and put for $\varepsilon > 0$
\begin{align}\label{Mfgtest}
\begin{split}
\sigma (x,\xi,\eta)
&= \sigma (\xi,\eta)
= \widehat u (\xi/\varepsilon) \, \widehat v (\eta),
\\
\widehat f (\xi) 
&= \varepsilon^{n/p-n} \, \widehat u (\xi/\varepsilon), \quad
\widehat g (\xi) 
= \varepsilon^{n/q-n} \, \widehat u \big( (\eta - e_1)/\varepsilon \big),
\end{split}
\end{align}
where $e_1=(1,0,\dots,0)\in\R^n$.
Then,
\begin{equation} \label{assumptionright23}
\|T_{\sigma} (f,g)\|_{L^r} \approx 1,\quad
\| f \|_{L^p} \approx 1,
\andd
\| g \|_{L^q} \approx 1
\end{equation}
for a sufficiently small $\varepsilon >0$.
In fact, since $\widehat{v}=1$ on the support of $\widehat{u} ((\cdot-e_1)/\varepsilon)$
by choosing a suitably small $\varepsilon>0$, we have
\begin{align*}
T_{\sigma} (f,g) (x)
=
\varepsilon^{n/p} \big( u \ast u \big) (\varepsilon x)
\, \varepsilon^{n/q} e^{ix_1} u (\varepsilon x),
\end{align*}
where $x=(x_1,\dots,x_n)\in\R^n$.
Hence, $\|T_{\sigma} (f,g)\|_{L^r} \approx 1$,
since $1/p+1/q=1/r$.
The second and third equivalences are obvious.

Next, we let a function 
$\Psi_0 \in \calS(\R^{n} \times \R^{n})$ satisfy that
for a sufficiently small $\delta > 0$
\begin{align*}
\supp \Psi_0
&\subset \left\{ (\xi,\eta) \in \R^n\times\R^n : |(\xi,\eta)| \leq 2^{1/2 + \delta} \right\} ,
\\
\Psi_0 &= 1 \quad\textrm{on}\quad 
\left\{ (\xi,\eta) \in \R^n\times\R^n : |(\xi,\eta)| \leq 2^{1/2 - \delta} \right\}.
\end{align*}
We put $\Psi = \Psi_{0} - \Psi_{0}(2\cdot,2\cdot)$ and $\Psi_{j} = \Psi (\cdot/2^j,\cdot/2^j)$, $j \in \N$.
Then, we have
\begin{align*}
\supp \Psi 
&\subset \left\{ (\xi,\eta) \in \R^n\times\R^n : 2^{-1/2 - \delta} \leq |(\xi,\eta)| \leq 2^{1/2 + \delta} \right\} ,
\\
\Psi &= 1 \quad\textrm{on}\quad 
\left\{ (\xi,\eta) \in \R^n\times\R^n : 2^{-1/2 + \delta} \leq |(\xi,\eta)| \leq 2^{1/2 - \delta} \right\},
\end{align*}
and $\sum_{j\in\N_{0}} \Psi_{j} = 1$ on $\R^n \times \R^n$.
This $\{ \Psi_{j} \}_{j\in\N_{0}}$ is 
a Littlewood--Paley partition of unity on $\R^n \times \R^n$.
Now, observe that
$\supp \sigma \subset
\{ (\xi,\eta) \in \R^n\times\R^n: |(\xi,\eta)| \leq 2^{1/2 - \delta} \}$,
choosing $\varepsilon > 0$ and $\delta > 0$ suitably.
Then, we see that
\begin{equation*}
\sigma_{j}^\rho (\xi,\eta)
=
\sigma (2^{j\rho}\xi, 2^{j\rho}\eta)
\Psi_j (2^{j\rho}\xi, 2^{j\rho}\eta)
=
\left\{
\begin{array}{ll}
\sigma,
& \textrm{if} \quad j=0,
\\
0,
& \textrm{if} \quad j\neq0.
\end{array}
\right.
\end{equation*}
Moreover, since the symbol $\sigma$ is independent of $x$,
we have
\begin{equation*}
\Delta_{\boldsymbol{k}} [\sigma_{j}^\rho] (x,\xi, \eta)
= \left\{
\begin{array}{ll}
\psi_{k_1} (D_\xi) \psi_{k_2} (D_\eta) [\sigma_{j}^\rho] (\xi, \eta), &\textrm{if}\quad k_0=0,
\\
0, & \textrm{if}\quad k_0\neq0.
\end{array}
\right.
\end{equation*}
These two facts mean that
\begin{align*}
\| \sigma \|_{BS^{m}_{\rho,\rho}(\boldsymbol{s}; \R^{n})}
= \sum_{k_1,k_2 \in \N_{0}} 2^{s_1k_1+s_2k_2} 
\| \psi_{k_1} (D_\xi) \psi_{k_2} (D_\eta) \sigma (\xi, \eta) \|_{ L^{2}_{ul} ((\R^{n})^2) }.
\end{align*}
Therefore, from
the embedding $L^{2} \hookrightarrow L^{2}_{ul}$
and Lemma \ref{sugimoto},
we obtain
\begin{align} \label{assumptionright1}
\begin{split}
\| \sigma \|_{BS^{m}_{\rho,\rho}(\boldsymbol{s}; \R^{n})}
&\leq \sum_{k_1,k_2 \in \N_{0}} 2^{s_1k_1+s_2k_2} 
\| \psi_{k_1} (D_\xi) \psi_{k_2} (D_\eta) \sigma (\xi, \eta) \|_{ L^2 ((\R^{n})^2) },
\\&= \| \widehat u (\cdot/\varepsilon) \|_{ B^{s_1}_{2,1} (\R^n) } \, 
\| \widehat v \|_{ B^{s_2}_{2,1} (\R^n) }
\lesssim \varepsilon^{n/2-s_1}
\end{split}
\end{align} 
for $\varepsilon>0$ sufficiently small.

Test \eqref{Mfgtest} to the assumption \eqref{rhoassumptionofs12}.
Then, by \eqref{assumptionright23} and \eqref{assumptionright1},
we have $\varepsilon^{n/2-s_1} \gtrsim 1$ for $\varepsilon > 0$ sufficiently small.
This means that $s_1 \geq n/2$.
\end{proof}


\subsection{Sharpness of the smoothness $s_0$}

In this subsection, we show the following.

\begin{proposition}\label{rhosharps0}
Let $0 \leq \rho < 1$, $m=-(1-\rho)n/2$, 
$\boldsymbol{s}=(s_0,s_1,s_2)\in (0,\infty)^3$, 
and $1 \leq p,q,r \leq \infty$ satisfy $1/p+1/q=1/r$. 
Suppose that the inequality
\begin{equation}\label{rhoassumptionofs0}
\| T_\sigma (f,g) \|_{L^r}
\lesssim 
\| \sigma \|_{BS^{m}_{\rho,\rho}(\boldsymbol{s}; \R^{n})} \| f \|_{L^p} \| g \|_{L^q}
\end{equation}
holds for all $\sigma \in BS^{m}_{\rho,\rho}(\boldsymbol{s}; \R^{n})$ and $f,g\in\calS(\R^n)$.
Then, $s_0 \geq n/2$.
\end{proposition}

To prove this, 
we employ a strategy by Miyachi--Tomita 
\cite[Appendix A]{miyachi tomita 2013 IUMJ}.
Define
\begin{equation*}
\| \sigma \|_{BS_{\rho,\rho}^{m,\dag} (\boldsymbol{s};\R^{n})}
=
\sup_{ j\in\N_0, \boldsymbol{k} \in (\N_0)^3 }
2^{-jm + \boldsymbol{k} \cdot \boldsymbol{s}}
\big\|
\Delta_{\boldsymbol{k}} [ \sigma_{j}^{\rho} ]
\big\|_{ L^\infty ((\R^{n})^3) },
\end{equation*}
with the same notations 
as in Definition \ref{sigmanorm12}.
Then, we have the following.

\begin{lemma}\label{lemsharps0}
Let $0<\rho<1$, $m\in\R$, $\boldsymbol{s}=(s_0,s_1,s_2)\in (0,\infty)^3$, 
and $1 \leq p,q,r \leq \infty$ satisfy $1/p+1/q=1/r$.
Suppose that the inequality
\begin{equation}\label{rhoto0assumption}
\| T_\varsigma (f,g) \|_{L^r}
\lesssim 
\| \varsigma \|_{BS^{m,\dag}_{\rho,\rho}(\boldsymbol{s}; \R^{n})} \| f \|_{L^p} \| g \|_{L^q}
\end{equation}
holds for all $\varsigma \in BS^{m,\dag}_{\rho,\rho}(\boldsymbol{s}; \R^{n})$
and $f,g\in\calS(\R^n)$.
Then, the inequality
\begin{equation*}
\| T_\sigma (f,g) \|_{L^r}
\lesssim 
\| \sigma \|_{BS^{m',\dag}_{0,0}(\boldsymbol{s}; \R^{n})} \| f \|_{L^p} \| g \|_{L^q}
\end{equation*} 
holds for all $\sigma \in BS^{m',\dag}_{0,0}(\boldsymbol{s}; \R^{n})$ 
with $m' < m/(1-\rho)$
and $f,g\in\calS(\R^n)$.
\end{lemma}

To prove this, we will use the following lemma given by 
Sugimoto \cite{sugimoto 1988 TMJ}.

\begin{lemma}\label{sugimoto2}
Let $\{ \psi_{k} \}_{k \in \N_0}$
be the Littlewood--Paley partition of unity on $\R^n$.
\begin{enumerate}
\setlength{\itemindent}{0pt} 
\setlength{\itemsep}{5pt} 
\item 
Let $(s_1,s_2) \in (0,\infty)^2$. Then,
\begin{align*}
&
\sup_{k_1,k_2 \in \N_0}
2^{s_1k_1+s_2k_2}
\left\| \psi_{k_1}(D_\xi) \psi_{k_2}(D_\eta) [f_1 f_2](\xi,\eta)
\right\|_{ L^\infty((\R^n)^2) } 
\\&\lesssim
\prod_{i=1,2}
\sup_{k_1,k_2 \in \N_0}
2^{s_1k_1+s_2k_2}
\left\| \psi_{k_1}(D_\xi) \psi_{k_2}(D_\eta) [f_i](\xi,\eta)
\right\|_{ L^\infty((\R^n)^2) }.
\end{align*}

\item 
Let $\boldsymbol{s}=(s_0,s_1,s_2)\in(0,\infty)^3$.
Then,
\begin{align*}
&
\sup_{\boldsymbol{k} \in (\N_0)^3}
2^{\boldsymbol{k} \cdot \boldsymbol{s}}
\left\| \Delta_{\boldsymbol{k}} [f(\lambda_0 \cdot, \lambda_1 \cdot, \lambda_2 \cdot )]
\right\|_{ L^\infty((\R^n)^3) } 
\\&\lesssim
\max(1,\lambda_0^{s_0})
\max(1,\lambda_1^{s_1})
\max(1,\lambda_2^{s_2})
\sup_{\boldsymbol{k} \in (\N_0)^3}
2^{\boldsymbol{k} \cdot \boldsymbol{s}}
\left\| \Delta_{\boldsymbol{k}} f
\right\|_{ L^\infty((\R^n)^3) }
\end{align*}
for $\lambda_0,\lambda_1,\lambda_2 \in (0,\infty)$.
\end{enumerate}
\end{lemma}

The assertion (1) was stated in
\cite[Theorem 1.4]{sugimoto 1988 TMJ}.
See also \cite[Theorem and Remark 1 in Section 2.8.2]{triebel 1983}.
The two variables version of the assertion (2) 
was mentioned in \cite[Proposition 1.1]{sugimoto 1988 TMJ}.
Following the same lines as there,
one can show the assertion (2)
for the three variables.
Therefore, we omit these proofs.

\begin{proof}[Proof of Lemma \ref{lemsharps0}]
Assume $\sigma \in BS^{m',\dag}_{0,0}(\boldsymbol{s}; \R^{n})$
with $m' < m/(1-\rho)$.
Let $\{ \Psi_\ell \}_{\ell\in\N_0}$ be 
a Littlewood--Paley partition of unity on $(\R^n)^2$.
Then,
\begin{equation}\label{sigmaell}
\sigma ( x, \xi, \eta )
=
\sum_{\ell \in \N_0} \sigma ( x, \xi, \eta ) \Psi_\ell ( \xi, \eta )
=\sum_{\ell\in\N_0} \sigma_\ell ( x, \xi, \eta )
\end{equation}
with $\sigma_\ell ( x, \xi, \eta ) = \sigma ( x, \xi, \eta ) \Psi_\ell ( \xi, \eta )$.
For simplicity, we write $\varrho = \frac{\rho}{1-\rho}$ for $0<\rho<1$.
We have by changes of variables 
\begin{align}\label{sigmaelltotauell}
\begin{split}
&
T_{\sigma_\ell} (f,g) (x)
\\&=
\frac{1}{(2\pi)^{2n}} 
\int_{(\R^{n})^2} e^{i 2^{-\ell\varrho} x \cdot (\xi + \eta) }
\sigma_\ell (x, 2^{-\ell\varrho} \xi, 2^{-\ell\varrho} \eta) 
	\widehat{f (2^{\ell\varrho} \cdot )} (\xi) \widehat{g(2^{\ell\varrho} \cdot )} (\eta)
\, d\xi d\eta
\\&=
T_{\varsigma_\ell}	(f_\ell, g_\ell)(2^{-\ell\varrho}x)
\end{split}
\end{align}
with
\begin{equation*}
\varsigma_\ell 
= \sigma_\ell (2^{\ell\varrho} \cdot, 2^{-\ell\varrho}\cdot, 2^{-\ell\varrho}\cdot),
\quad
f_\ell 
= f (2^{\ell\varrho} \cdot ), 
\andd
g_\ell =g(2^{\ell\varrho} \cdot ).
\end{equation*}
Then, $\varsigma_\ell \in BS^{m,\dag}_{\rho,\rho}(\boldsymbol{s}; \R^{n})$,
and more precisely, the inequality
\begin{equation}\label{tauell}
\| \varsigma_\ell \|_{ BS^{m,\dag}_{\rho,\rho}(\boldsymbol{s}; \R^{n}) }
\lesssim
2^{\ell (m' - m/(1-\rho))} 
\| \sigma \|_{ BS^{m',\dag}_{0,0}(\boldsymbol{s}; \R^{n}) }
\end{equation}
holds for any $\ell \in \N_0$.
We shall prove \eqref{tauell}.
Since 
\[
\supp \varsigma_\ell \subset 
\left\{ (\xi,\eta) \in \R^n \times \R^n :
2^{\ell/(1-\rho) -1} \leq | (\xi,\eta) | \leq 2^{\ell/(1-\rho) +1}
\right\},
\]
recalling the notation of $(\varsigma_\ell)_j^\rho$:
\[
(\varsigma_\ell)_j^\rho (x,\xi,\eta)
=
\varsigma_\ell ( 2^{-j\rho} x, 2^{j\rho} \xi, 2^{j\rho} \eta) \, \Psi_j ( 2^{j\rho} \xi, 2^{j\rho} \eta),
\]
we see that $j$ must be in the set 
\begin{equation*}
\Omega_\ell = 
\left\{ j \in \N_0: 
\max \left( 0, \frac{\ell}{1-\rho} -2 \right) \leq j \leq \frac{\ell}{1-\rho} +2 \right\}
\end{equation*}
(otherwise, $(\varsigma_\ell)_j^\rho$ vanishes).
Then, we have
\begin{equation} \label{tauell1}
\| \varsigma_\ell \|_{ BS^{m,\dag}_{\rho,\rho}(\boldsymbol{s}; \R^{n}) }
\approx
2^{-\ell m/(1-\rho)}
\sup_{j\in\Omega_\ell} \sup_{\boldsymbol{k} \in (\N_0)^3}
2^{\boldsymbol{k} \cdot \boldsymbol{s}}
\left\| \Delta_{\boldsymbol{k}} \left[ (\varsigma_\ell)_j^\rho \right] \right\|_{ L^\infty } 
\end{equation}
for $\ell \in \N_{0}$. 
Here, since
\begin{align*}
\Delta_{\boldsymbol{k}} \left[ (\varsigma_\ell)_j^\rho \right]
=
\psi_{k_1}(D_\xi) \psi_{k_2}(D_\eta) \Big[ \psi_{k_0}(D_x) 
\big[ \varsigma_\ell ( 2^{-j\rho} x, 2^{j\rho} \xi, 2^{j\rho} \eta) \big] \times
\Psi_j ( 2^{j\rho} \xi, 2^{j\rho} \eta) \Big] ,
\end{align*}
it holds from Lemma \ref{sugimoto2} (1) that
\begin{equation}\label{tauell2}
\sup_{\boldsymbol{k} \in (\N_0)^3}
2^{\boldsymbol{k} \cdot \boldsymbol{s}}
\left\| \Delta_{\boldsymbol{k}} \left[ (\varsigma_\ell)_j^\rho \right] \right\|_{ L^\infty } 
\lesssim
\sup_{\boldsymbol{k} \in (\N_0)^3}
2^{\boldsymbol{k} \cdot \boldsymbol{s}}
\left\| \Delta_{\boldsymbol{k}} 
\big[ \varsigma_\ell (2^{-j\rho} \cdot, 2^{j\rho}\cdot, 2^{j\rho}\cdot) \big]
\right\|_{ L^\infty }
\end{equation}
for $j \in N_0$, where we used the fact that
\[
\sup_{k_1,k_2\in\N_0} 
2^{s_1k_1 + s_2k_2} 
\left\| 
\psi_{k_1} (D_\xi) \psi_{k_2} (D_\eta) 
	\left[ \Psi_j ( 2^{j\rho} \cdot,2^{j\rho} \cdot ) \right] 
\right\|_{L^\infty}
\lesssim1
\] holds
with the implicit constant
independent of $j \in \N_0$.
Moreover, since 
\[
\varsigma_\ell (2^{-j\rho} \cdot, 2^{j\rho}\cdot, 2^{j\rho}\cdot)
=
\sigma_\ell (2^{\ell\varrho-j\rho} \cdot, 2^{-(\ell\varrho-j\rho)} \cdot, 2^{-(\ell\varrho-j\rho)} \cdot)
\]
and $2^{\ell\varrho-j\rho} \approx 1$ for $j \in \Omega_\ell$,
we have by Lemma \ref{sugimoto2} (2)
\begin{equation}\label{tauell3}
\sup_{\boldsymbol{k} \in (\N_0)^3}
2^{\boldsymbol{k} \cdot \boldsymbol{s}}
\left\| \Delta_{\boldsymbol{k}} 
\big[ \varsigma_\ell (2^{-j\rho} \cdot, 2^{j\rho}\cdot, 2^{j\rho}\cdot) \big]
\right\|_{ L^\infty }
\lesssim
\sup_{\boldsymbol{k} \in (\N_0)^3}
2^{ \boldsymbol{k} \cdot \boldsymbol{s} }
\left\| \Delta_{\boldsymbol{k}} 
\big[ \sigma_\ell \big]
\right\|_{ L^\infty } , \quad
j \in \Omega_\ell.
\end{equation}
Collecting \eqref{tauell1}, \eqref{tauell2}, and \eqref{tauell3}, we obtain
\begin{align*}
\| \varsigma_\ell \|_{ BS^{m,\dag}_{\rho,\rho}(\boldsymbol{s}; \R^{n}) }
&\lesssim
2^{-\ell m/(1-\rho)}
\sup_{\boldsymbol{k} \in (\N_0)^3}
2^{ \boldsymbol{k} \cdot \boldsymbol{s} }
\left\| \Delta_{\boldsymbol{k}} 
\big[ \sigma_\ell \big]
\right\|_{ L^\infty }
\\&\leq
2^{-\ell m/(1-\rho)}
\, 2^{\ell m'}
\left\| \sigma \right\|_{ BS^{m',\dag}_{0,0}(\boldsymbol{s}; \R^{n}) },
\end{align*}
since $\sigma_\ell = \sigma_{\ell}^{0}$.
This is the desired inequality \eqref{tauell}.

Therefore, we have by 
\eqref{sigmaell}, \eqref{sigmaelltotauell}, \eqref{rhoto0assumption}, and \eqref{tauell}
\begin{align*}
\| T_\sigma (f,g) \|_{L^r} 
&\leq 
\sum_{\ell\in\N_0} 
\| T_{\sigma_\ell} (f,g) \|_{L^r} 
=
\sum_{\ell\in\N_0} 2^{\ell\varrho n/r}
\| T_{\varsigma_\ell} (f_\ell, g_\ell) \|_{L^r}
\\
&\lesssim
\sum_{\ell\in\N_0} 2^{\ell\varrho n/r}
\| \varsigma_\ell \|_{BS^{m,\dag}_{\rho,\rho}(\boldsymbol{s}; \R^{n})} \| f_\ell \|_{L^p} \| g_\ell \|_{L^q} 
\\
&\lesssim
\sum_{\ell\in\N_0} 
2^{\ell (m' - m/(1-\rho))} 
\left\| \sigma \right\|_{ BS^{m',\dag}_{0,0}(\boldsymbol{s}; \R^{n}) } \| f \|_{L^p} \| g \|_{L^q}
\\
&\approx
\left\| \sigma \right\|_{ BS^{m',\dag}_{0,0}(\boldsymbol{s}; \R^{n}) } \| f \|_{L^p} \| g \|_{L^q},
\end{align*}
since $1/p+1/q=1/r$ and $m' < m/(1-\rho)$.
This completes the proof.
\end{proof}

We also use the following lemma proved 
by Wainger \cite[Theorem 10]{wainger 1965}
and Miyachi--Tomita \cite[Lemma 6.1]{miyachi tomita 2013 IUMJ}.

\begin{lemma}\label{lemmawainger}
Let $1\leq p \leq \infty$, $0<a<1$, $0<b<n$ and $\varphi\in\calS(\R^n)$.
For $t>0$, put
\begin{equation*}
f_{a,b,t}(x)
= \sum_{k\in\Z^n\setminus\{0\}} e^{-t|k|} |k|^{-b} e^{i |k|^a} e^{i k \cdot x} \varphi(x).
\end{equation*}
Then, if $b > n(1 - a/2-1/p+a/p)$,
we have $\sup_{t>0} \| f_{a,b,t} \|_{L^p(\R^n)} < \infty$.
\end{lemma}

By using this lemma, we shall show the following.
See also \cite[Proposition 7.3]{KMT-arxiv-2}.

\begin{lemma}\label{0sharps0}
Let $m \geq -n$, $\boldsymbol{s}=(s_0,s_1,s_2)\in[0,\infty)^3$, and $1 \leq p,q,r \leq \infty$. 
Suppose that
\begin{equation}\label{assumptionofs0}
\| T_\sigma (f,g) \|_{L^r}
\lesssim 
\| \sigma \|_{BS^{m,\dag}_{0,0}(\boldsymbol{s}; \R^{n})} \| f \|_{L^p} \| g \|_{L^q}
\end{equation}
holds for all $\sigma \in BS^{m,\dag}_{0,0}(\boldsymbol{s}; \R^{n})$ and $f,g\in\calS(\R^n)$.
Then, $s_0 \geq m+n$.
\end{lemma}

\begin{proof}
In this proof, we will use nonnegative functions 
$\varphi, \widetilde{\varphi} \in \calS(\R^n)$
satisfying that 
$\supp \varphi \subset [-1/4,1/4]^n$, 
$\widetilde\varphi = 1$ on $[-1/4,1/4]^n$,
and $\supp \widetilde\varphi \subset [-1/2,1/2]^n$.
Define
\begin{align*}
\sigma_{a_1,a_2} (x,\xi,\eta)
&= \varphi(x) e^{-ix\cdot(\xi+\eta)}
\\&\quad\times
\sum_{k,\ell\in\Z^n\setminus\{0\}}
(1+|k|+|\ell|)^{m-s_0} e^{-i|k|^{a_1}} e^{-i|\ell|^{a_2}} \varphi(\xi-k) \varphi(\eta-\ell),
\\
f_{a_1,b_1,t}(x)
&=
\sum_{\nu\in\Z^n\setminus\{0\}}
e^{-t|\nu|} |\nu|^{-b_1} e^{i|\nu|^{a_1}} e^{i \nu \cdot x} \calFi \widetilde \varphi (x),
\\
g_{a_2,b_2,t}(x)
&=
\sum_{\mu\in\Z^n\setminus\{0\}}
e^{-t|\mu|} |\mu|^{-b_2} e^{i|\mu|^{a_2}} e^{i \mu \cdot x} \calFi \widetilde \varphi (x),
\end{align*}
where $t>0$, $0 < a_1,a_2 < 1$,
$b_1 = n (1 - a_1/2 - 1/p + a_1/p ) + \varepsilon_1$, and
$b_2 = n (1 - a_2/2 - 1/q + a_2/q ) + \varepsilon_2$,
with $\varepsilon_1,\varepsilon_2 > 0$.
Note that $f_{a_1,b_1,t}, g_{a_2,b_2,t} \in \calS(\R^n)$, 
thanks to the exponential decay factors.

For these functions, the following hold:
\begin{align}\label{s0testsigmafg}
\begin{split}
&\| \sigma_{a_1,a_2} \|_{ BS^{m,\dag}_{0,0}(\boldsymbol{s}; \R^{n}) } \lesssim 1,
\\
&\sup_{t>0} \| f_{a_1,b_1,t} \|_{L^p(\R^n)} \lesssim 1, \quad
\sup_{t>0} \| g_{a_2,b_2,t} \|_{L^q(\R^n)} \lesssim 1.
\end{split}
\end{align}
The second and third inequalities
follow from Lemma \ref{lemmawainger}.
We shall consider the first inequality.
We write $N_i=[s_i]+1$, $i=0,1,2$
and recall the notation
$\sigma_j^{0}(x,\xi,\eta)
=\sigma (x,\xi,\eta)\Psi_j (\xi, \eta)$
for $\rho = 0$.
Then, observing that
\begin{align*}
\left| \partial_x^\alpha \partial_\xi^\beta \partial_\eta^\gamma \sigma_{a_1,a_2} (x,\xi,\eta) \right|
\lesssim
(1+|\xi|+|\eta|)^{|\alpha|+m-s_0},
\end{align*}
and that $1+|\xi|+|\eta| \approx 2^j$ on the support of $\Psi_j$,
we realize that
\begin{align*}
\left| \partial_x^\alpha \partial_\xi^\beta \partial_\eta^\gamma (\sigma_{a_1,a_2})_j^0 (x,\xi,\eta) \right|
\lesssim
2^{j(|\alpha|+m-s_0)}.
\end{align*}
Hence, as in the proof of Lemma \ref{derivationlemma}, 
applying the Taylor expansion (see the argument around \eqref{taylor})
to $ \Delta_{\boldsymbol{k}} [(\sigma_{a_1,a_2})_j^0]$, 
we have
\begin{equation}\label{withN0}
\left\| \Delta_{\boldsymbol{k}} \left[ (\sigma_{a_1,a_2})_j^0 \right] \right\|_{L^\infty}
\lesssim 
2^{-k_0N_0-k_1N_1-k_2N_2} \, 2^{j(N_0+m-s_0)}
\end{equation}
for $j\geq0$ and $\boldsymbol{k} \in (\N_{0})^3$.
On the other hand,
in the derivation of \eqref{taylor},
by avoiding the Taylor expansion with respect to the $x$-variable, 
we have
\begin{align*}
&
\Delta_{\boldsymbol{k}} \sigma(x,\xi,\eta)
\\
&=
2^{n(k_0+k_1+k_2)} 
\sum_{|\beta| = N_1} \frac{1}{\beta!} 
\sum_{|\gamma| = N_2} \frac{1}{\gamma!} 
\int_{(\R^{n})^3}
\check \psi (2^{k_0}x') \,
\check \psi (2^{k_1}\xi') (-\xi')^{\beta} \,
\check \psi (2^{k_2}\eta') (-\eta')^{\gamma} 
\\
&\quad\times
\int_{[0,1]^2} \bigg( \prod_{i=1,2} N_i (1-t_i)^{N_i-1} \bigg)
\Big( \partial_{\xi}^{\beta} \partial_{\eta}^{\gamma} \sigma \Big)
(x-x', \xi-t_1 \xi', \eta-t_2 \eta' ) 
\, dTdX',
\end{align*}
where $dT = dt_1 dt_2$ and $dX' = dx' d\xi' d\eta'$.
Then, we have by the same lines as above
\begin{equation}\label{withoutN0}
\left\| \Delta_{\boldsymbol{k}} \left[ (\sigma_{a_1,a_2})_j^0 \right] \right\|_{L^\infty}
\lesssim 
2^{-k_1N_1-k_2N_2} \, 2^{j(m-s_0)}
\end{equation}
for $j\geq0$ and $\boldsymbol{k} \in (\N_{0})^3$.
Take $0< \theta < 1$ such that $N_0 \theta = s_0$.
By \eqref{withN0} and \eqref{withoutN0}
\begin{align*}
\left\| \Delta_{\boldsymbol{k}} \left[ (\sigma_{a_1,a_2})_j^0 \right] \right\|_{L^\infty}
&=
\left\| \Delta_{\boldsymbol{k}} \left[ (\sigma_{a_1,a_2})_j^0 \right] \right\|_{L^\infty}^{\theta}
\left\| \Delta_{\boldsymbol{k}} \left[ (\sigma_{a_1,a_2})_j^0 \right] \right\|_{L^\infty}^{1-\theta}
\\
&\lesssim
2^{-k_0 s_0} \, 2^{-k_1 N_1 -k_2 N_2} \, 2^{jm}.
\end{align*}
Therefore, we obtain the first inequality in \eqref{s0testsigmafg}.

We next investigate the left hand side of \eqref{assumptionofs0}.
Observe that $\supp \varphi(\cdot-k) \, \cap \, \supp \widetilde{\varphi}(\cdot-k')=\emptyset$
if $k\neq k'$,
and $\widetilde\varphi = 1$ on $\supp \varphi$.
Then, 
\begin{align*}
T_{\sigma_{a_1,a_2}} (f_{a_1,b_1,t}, g_{a_2,b_2,t})(x)
&=
\frac{\| \varphi \|_{L^1}^2}{ (2\pi)^{2n} }
\, \varphi(x)
\sum_{k,\ell}
e^{-t(|k|+|\ell|)} 
(1+|k|+|\ell|)^{m-s_0} 
|k|^{-b_1}  |\ell|^{-b_2}.
\end{align*}
Taking the $L^r$ norm of both sides, we have
\begin{align}\label{s0testleft}
\begin{split}
\| T_{\sigma_{a_1,a_2}} (f_{a_1,b_1,t}, g_{a_2,b_2,t}) \|_{L^r}
&\gtrsim
\sum_{k\in\Z^n\setminus\{0\}} 
\sum_{\ell:\, 0<|\ell|\leq|k|} 
e^{-2t|k|}
(1+|k|)^{m-s_0} 
|k|^{-b_1}  |k|^{-b_2}
\\
&\approx
\sum_{k\in\Z^n\setminus\{0\}} 
e^{-2t|k|}
(1+|k|)^{m-s_0-b_1-b_2+n}.
\end{split}
\end{align}

Collecting \eqref{assumptionofs0}, \eqref{s0testsigmafg}, and \eqref{s0testleft},
we obtain
\begin{equation*}
\sum_{k\in\Z^n\setminus\{0\}} 
e^{-2t|k|}
(1+|k|)^{m-s_0-b_1-b_2+n}
\lesssim 1
\end{equation*}
with the implicit constant independent of $t$.
Thus, we have by the Fatou lemma
\begin{equation*}
\sum_{k\in\Z^n\setminus\{0\}} 
(1+|k|)^{m-s_0-b_1-b_2+n}
\lesssim 1.
\end{equation*}
This yields $m-s_0-b_1-b_2+n < -n$,
which is identical with
\begin{equation*}
s_0 > m+2n
-n \left(1 - \frac{a_1}{2}-\frac{1}{p}+\frac{a_1}{p} \right) 
-n \left(1 - \frac{a_2}{2}-\frac{1}{q}+\frac{a_2}{q} \right) 
-\varepsilon_1-\varepsilon_2.
\end{equation*}
Since $0< a_i <1 $ and $\varepsilon_i>0$, $i=1,2$, are arbitrary,
if we take the limits as $a_i\to 1$ and $\varepsilon_i \to 0$, 
we obtain the condition $s_0 \geq m+n$, which gives the desired result.
\end{proof}

\begin{proof}[Proof of Proposition \ref{rhosharps0}]
We first observe that
$BS^{m,\dag}_{\rho,\rho}( \boldsymbol{s}_{\varepsilon} ; \R^{n})
\subset BS^{m}_{\rho,\rho}( \boldsymbol{s} ; \R^{n})$ holds
with 
$\boldsymbol{s}_{\varepsilon}
=(s_0+\varepsilon,s_1+\varepsilon,s_2+\varepsilon)$
for any $\varepsilon > 0$.
Now, we first consider the case $0 < \rho < 1$.
In this case, we have by \eqref{rhoassumptionofs0} and the inclusion relation above
\[
\| T_\sigma (f,g) \|_{L^r}
\lesssim 
\| \sigma \|_{BS^{-(1-\rho)n/2,\dag}_{\rho,\rho}(\boldsymbol{s}_{\varepsilon}; \R^{n})} 
\| f \|_{L^p} \| g \|_{L^q},
\]
and then, by Lemma \ref{lemsharps0}
\begin{equation*}
\| T_\sigma (f,g) \|_{L^r}
\lesssim 
\| \sigma \|_{BS^{-n/2-\delta,\dag}_{0,0}(\boldsymbol{s}_{\varepsilon}; \R^{n})} 
\| f \|_{L^p} \| g \|_{L^q}
\end{equation*}
for any $\delta>0$.
Thus, we conclude by Lemma \ref{0sharps0} 
that $s_0+\varepsilon \geq n/2-\delta$.
Since $\varepsilon>0$ and $\delta>0$ are both arbitrary,
if we take the limits as $\varepsilon \to 0$ and $\delta \to 0$,
we obtain $s_0 \geq n/2$,
which is the desired result for the case $0<\rho<1$.
The case $\rho=0$ is similarly proved
by the embedding above, Lemma \ref{0sharps0}, 
and a limit argument.
\end{proof}



\appendix
\section{Existence of decomposition}
\label{appexistdecom}

In this appendix, 
we determine functions 
used to decompose symbols 
in Section \ref{secchangeform}.
Let $\phi \in \calS (\R^n)$ satisfy that
$\supp \phi \subset \{ \xi \in\R^n : | \xi | \leq 2 \}$ and
$\phi = 1$ for $| \xi | \leq 1$.
For $k \in \Z$, we write $\phi_k = \phi (\cdot / 2^k )$.
Then, we see that
$\supp \phi_k \subset \{ \xi \in\R^n : | \xi | \leq 2^{k+1} \}$
and $\supp (1-\phi_k) \subset \{ \xi \in \R^n : | \xi | \geq 2^{k} \}$
for $k \in \Z$,
since $\phi_k = 1$ for $| \xi | \leq 2^{k}$.

Let $\{ \Psi_{j} \}_{j\in \N_0}$
be a Littlewood--Paley partition on $(\R^n)^2$.
Then, for $j \geq 1$, 
$\Psi_j$ can be expressed into the following form:
\begin{align*}
\Psi_j ( \xi, \eta )
&=
\Psi_j ( \xi, \eta ) 
\phi_{j+1} (\xi) \phi_{j+1} (\eta)
\\
&\times
\left\{
	\phi_{j-6} (\xi) + (1- \phi_{j-6} (\xi) )
\right\}
\left\{
	\phi_{j-6} (\eta) + (1- \phi_{j-6} (\eta) )
\right\}.
\end{align*}
Since $\phi_{j} \, \phi_{j'} = \phi_{j'} $ if $j>j'$ and 
$\phi_{j-6} (\xi) \phi_{j-6} (\eta)$ vanishes 
on $\supp \Psi_j$, $j \geq 1$, we have
\begin{align*}
\Psi_j ( \xi, \eta )  
&=
\Psi_j ( \xi, \eta ) \phi_{j-6} (\xi) \,
\phi_{j+1} (\eta) (1- \phi_{j-6} (\eta) )
\\&\quad+
\Psi_j ( \xi, \eta ) \phi_{j+1} (\xi) (1- \phi_{j-6} (\xi) ) \,
\phi_{j-6} (\eta) 
\\&\quad+
\Psi_j ( \xi, \eta ) \phi_{j+1} (\xi) (1- \phi_{j-6} (\xi) ) \, 
\phi_{j+1} (\eta) (1- \phi_{j-6} (\eta) ) .
\end{align*}
We further decompose the first factor above as follows:
\[
\Psi_j ( \xi, \eta ) \phi_{j-6} (\xi) \,
\phi_{j+1} (\eta) (1- \phi_{j-6} (\eta) )
\left\{
	\phi_{j-4} (\eta) + (1- \phi_{j-4} (\eta) )
\right\}.
\]
Then, this is equal to
$\Psi_j ( \xi, \eta ) \phi_{j-6} (\xi) \,
\phi_{j+1} (\eta) (1- \phi_{j-4} (\eta) ) $, since 
$\phi_{j-6} (\xi) \phi_{j-4} (\eta)$ vanishes 
on $\supp \Psi_j$, $j \geq 1$, 
and $(1-\phi_{j}) (1-\phi_{j'}) = (1-\phi_{j})$ if $j>j'$.
The second factor can be expressed similarly
because of symmetry.
Therefore, we have
\begin{align*}
\Psi_j ( \xi, \eta )  
&=
\Psi_j ( \xi, \eta ) \phi_{j-6} (\xi) \,
\phi_{j+1} (\eta) \big(1- \phi_{j-4} (\eta) \big)
\\&\quad+
\Psi_j ( \xi, \eta ) \phi_{j+1} (\xi)\big(1- \phi_{j-4} (\xi) \big) \,
\phi_{j-6} (\eta) 
\\&\quad+
\Psi_j ( \xi, \eta ) \phi_{j+1} (\xi)\big(1- \phi_{j-6} (\xi) \big) \, 
\phi_{j+1} (\eta) \big(1- \phi_{j-6} (\eta) \big) 
\\
&=:
\Psi_j ( \xi, \eta )
\, \phi_{j}' (\xi) \psi_{j}' (\eta)
+
\Psi_j ( \xi, \eta ) 
\, \psi_{j}' (\xi) \phi_{j}' (\eta)
+
\Psi_j ( \xi, \eta ) 
\, \psi_{j}'' (\xi) \psi_{j}'' (\eta)
\end{align*}
with $\phi_{j}' = \phi' (\cdot/2^{j})$,
$\psi_{j}' = \psi' (\cdot/2^{j})$, and
$\psi_{j}'' = \psi'' (\cdot/2^{j})$,
and then we realize that 
\begin{align*}
\supp \phi' 
&\subset \left\{ \zeta \in\R^n : |\zeta| \leq 2^{-5}\right\},
\quad
\\
\supp \psi'
&\subset \left\{ \zeta \in\R^n : 2^{-4} \leq |\zeta| \leq 2^{2}\right\},
\\
\supp \psi''
&\subset \left\{ \zeta \in\R^n : 2^{-6} \leq |\zeta| \leq 2^{2}\right\}.
\end{align*}
Hence, 
we obtain the information \eqref{suppdecomphi'}, 
\eqref{suppdecompsi'}, and \eqref{suppdecompsi''} 
given in Section \ref{secchangeform}.



\section{Boundedness from $L^2 \times L^2$ to $L^1$}
\label{appl2l2l1bdd}

In this appendix, we shall prove the following boundedness
stated in Remark \ref{rems0critical}.
\begin{theorem}
Let $0 \leq \rho <1$, $m=-(1-\rho)n/2$, and
$\boldsymbol{s}=(s_0,s_1,s_2) \in [0,\infty)^3$ satisfy
$s_{0}, s_{1}, s_{2} \geq n/2$.
Then, if $\sigma \in BS_{\rho,\rho}^m (\boldsymbol{s};\R^{n})$,
the bilinear pseudo-differential operator $T_\sigma$ is bounded 
from $L^2 (\R^n) \times L^2 (\R^n) $ to $L^1 (\R^n) $.
\end{theorem}
The proof is much simpler than 
that of the boundedness to $h^1$
done in Section \ref{secl2l2h1}.

\begin{proof}
As in Section \ref{secchangeform} and Appendix \ref{appexistdecom},
we decompose a Littlewood--Paley partition 
$\{ \Psi_{j} \}_{j \in \N_0}$ on $(\R^n)^2$ 
into the following form:
For $j \geq 1$,
\begin{equation*}
\Psi_j ( \xi, \eta )
=
\Psi_j ( \xi, \eta )
\phi'_{j} (\xi) \psi'_{j} (\eta)
+
\Psi_j ( \xi, \eta ) 
\psi'_{j} (\xi) \phi'_{j} (\eta),
\end{equation*}
where $\phi'_{j} = \phi' (\cdot/2^{j})$,
$\psi'_{j} = \psi' (\cdot/2^{j})$,
$\supp \phi'
\subset \{ \zeta \in\R^n : |\zeta| \leq 2^{-2} \} $, and
$\supp \psi'
\subset \{ \zeta \in\R^n : 2^{-3} \leq |\zeta| \leq 2^{2} \} $.
Then, repeating the same lines as in Section \ref{secchangeform},
the dual form of $T_{\sigma}(f,g)$ can be expressed by the sum of
the forms $I_{0}$, $I_{1}$, and $I_{2}$ as follows.
The form $I_{0}$ is the same as in \eqref{i0}.
The form $I_{1}$ is the following:
\begin{align*}
I_{1}
=
\sum_{ j\gg1 }
\sum_{\boldsymbol{k} \in (\N_0)^3}
\sum_{\boldsymbol{\nu} \in (\Z^{n})^2} 2^{-j\rho n}
\int_{\R^n} 
T_{ \sigma_{ j, \boldsymbol{k}, \boldsymbol{\nu}}^\rho }
	(\square_{\nu_1} f_{j}', \square_{\nu_2} g_{j}' ) (x) \, {h_j(x)} \, dx,
\end{align*}
where
$f_{j}' = \phi'_{j}(D) f (2^{-j\rho} \cdot )$, 
$g_{j}' = \psi'_{j}(D) g (2^{-j\rho} \cdot )$, and 
$h_j = h(2^{-j\rho} \cdot )$.
Also, we have
\begin{align*}
\supp \widehat{f_{j}'}
&\subset \left\{ \xi \in\R^n : |\xi| \leq 2^{j(1-\rho)-2}\right\} ,
\\
\supp \widehat{g_{j}'}
&\subset \left\{ \eta \in\R^n : 2^{j(1-\rho)-3} \leq |\eta| \leq 2^{j(1-\rho)+2}\right\}.
\end{align*}
The form $I_{2}$ is in a symmetrical position with $I_{1}$,
and thus we omit stating it.

We shall consider the three forms above.
However, we only consider $I_{1}$,
since the proof for $I_{0}$ is exactly the same 
as in Section \ref{boundl1i0}
and the proof for $I_{2}$ is similar to
that for $I_{1}$ because of symmetry.
We take a Littlewood--Paley partition $\{ \psi_{\ell} \}$ on $\R^n$
and decompose the factor of $f$ as
\begin{align*}
I_{1}
=
\sum_{ j\gg1 }
\sum_{\boldsymbol{k} \in (\N_0)^3}
\sum_{\ell \in \N_0}
\sum_{\boldsymbol{\nu} \in (\Z^{n})^2} 
2^{-j\rho n}
\int_{\R^n} 
T_{ \sigma_{ j, \boldsymbol{k}, \boldsymbol{\nu}}^\rho }
	(\square_{\nu_1} \Delta_{\ell} f_{j}', \square_{\nu_2} g_{j}' ) (x) \, {h_j(x)} \, dx.
\end{align*}
Then, the sums over $\boldsymbol{\nu}$ and $\ell$ are restricted to
$\nu_1 \in \{ \nu_1 \in \Z^n : |\nu_1| \lesssim 2^{\ell} \}$,
$\nu_2 \in \{ \nu_2 \in \Z^n : |\nu_2| \lesssim 2^{j(1-\rho)} \}$, and
$\ell \leq j(1-\rho)$ (see Section \ref{sech1T1}).
Applying Lemma \ref{nu12lemma} (1) 
with $p=q=2$ and $r=\infty$
to the restricted sums, 
we have
\begin{align*}
\left| I_{1} \right|
\lesssim
\sum_{ j, \, \boldsymbol{k}}
\sum_{\ell :\, \ell \leq j(1-\rho)}
2^{ \ell n/2 }
\, 2^{ (k_0 + k_1 + k_2 ) n/2} 
\| \Delta_{\boldsymbol{k}} [ \sigma_j^\rho ] \|_{ L^{2}_{ul} }
\| \Delta_{\ell+j\rho} f \|_{L^2}
\| \psi'_{j}(D)g \|_{L^2}
\| h \|_{L^\infty}
\end{align*}
where we used the calculation as in \eqref{h1I11fg}
to the factors of $f$ and $g$.
Since we are not dividing the sum over $j$,
the right hand side above is simply bounded by
\begin{align*}
\| \sigma \|_{ BS_{\rho,\rho}^{m} (\boldsymbol{s};\R^{n}) }
\| h \|_{L^\infty}
\sum_{ j \gg1 }
\sum_{\ell :\, \ell \leq j(1-\rho)}
2^{ \ell n/2 } \, 2^{-j(1-\rho) n/2}
\| \Delta_{\ell+j\rho} f \|_{L^2}
\| \psi'_{j}(D)g \|_{L^2},
\end{align*}
where $m=-(1-\rho)n/2$ and 
$\boldsymbol{s} = (n/2,n/2,n/2)$.
The sums over $j$ and $\ell$ are
bounded by a constant times $\| f \|_{L^2} \| g \|_{L^2}$,
recalling the proof of \eqref{estofII}. 
Hence, we obtain
\begin{align*}
\left| I_{1} \right|
\lesssim
\| \sigma \|_{ BS_{\rho,\rho}^{m} (\boldsymbol{s};\R^{n}) }
\| f \|_{L^2} \| g \|_{L^2} \| h \|_{L^\infty},
\end{align*}
where $m=-(1-\rho)n/2$ and 
$\boldsymbol{s} = (n/2,n/2,n/2)$.
This completes the proof.
\end{proof}




\end{document}